\documentclass[10pt]{article}

\usepackage{amscd}      
\usepackage{amssymb}     
\usepackage[intlimits]{amsmath}     
\usepackage[all,v2]{xy}      
\xyoption{2cell}
\UseAllTwocells
\xyoption{frame}

\usepackage{theorem}

\newtheorem{prop}{Proposition}[section]  
\newtheorem{lem}[prop]{Lemma}

\newtheorem{cor}[prop]{Corollary}
\newtheorem{them}[prop]{Theorem}

\theorembodyfont{\upshape}

\newtheorem{defn}[prop]{Definition}

\newtheorem{numrmk}[prop]{Remark}

\newtheorem{numex}[prop]{Example}

\newtheorem{rmk}{Remark}

\newenvironment{pf}{\begin{trivlist}\item[]{\sc Proof.}}%
            {\nolinebreak $\Box$ \end{trivlist}}

\newcommand{\noprint}[1]{}

\newcommand{\comment}[1]{{\marginpar{\footnotesize #1}}}

\renewcommand{\tilde}{\widetilde}

\newcommand{\toto}{\rightrightarrows}
\newcommand{\qed}{{\nolinebreak $\,\,\Box$}}

\newcommand{\upst}{^{\ast}}

\newcommand{\lst}{_{\ast}}

\newcommand{\com}{^{\scriptscriptstyle\bullet}}
\newcommand{\lcom}{_{\scriptscriptstyle\bullet}}
\newcommand{\upcom}{^{\scriptscriptstyle\bullet}}
\newcommand{\argument}{{{\,\cdot\,}}}
\newcommand{\pt}{\mathop{pt}}
\newcommand{\cinf}{\mathop{C^\infty}}

\newcommand{\XX}{{\mathfrak X}}

\renewcommand{\SS}{{\mathfrak S}}

\newcommand{\YY}{{\mathfrak Y}}      
      
\newcommand{\PP}{{\mathfrak P}}

\newcommand{\MM}{{\mathfrak M}}

\newcommand{\RR}{{\mathfrak R}}

\newcommand{\zz}{{\mathbb Z}}
\newcommand{\hh}{{\mathbb H}}

\newcommand{\qq}{{\mathbb Q}}

\newcommand{\rr}{{\mathbb R}}

\renewcommand{\O}{{\cal O}}

\newcommand{\del}{\partial}

\newcommand{\resto}{{\,|\,}}
\newcommand{\st}{\mathrel{\mid}}

\newcommand{\tot}{\mathop{\rm tot}}

\newcommand{\pr}{\mathop{\rm pr}\nolimits}

\newcommand{\gr}{\mathop{\rm gr}\nolimits}

\newcommand{\Ad}{\mathop{\rm Ad}\nolimits}

\newcommand{\cosk}{\mathop{\rm cosk}\nolimits}

\newcommand{\id}{\mathop{\rm id}\nolimits}
\newcommand{\Hom}{\mathop{\rm Hom}\nolimits}

\newcommand{\bhom}{\mathop{\rm\bf Hom}\nolimits}

\newcommand{\Isomu}{\mathop{\underline{\rm Isom}}\nolimits}

\newcommand{\band}{\mathop{\underline{\rm Band}}\nolimits}
\newcommand{\End}{\mathop{\rm End}\nolimits}

\newcommand{\Aut}{\mathop{\rm Aut}\nolimits}
\newcommand{\Out}{\mathop{\rm Out}\nolimits}

\newcommand{\tto}{\longrightarrow}

\newcommand{\comp}{\mathbin{{\scriptstyle\circ}}}

\newcommand{\ldiag}[1]%
       {\makebox[0cm]{${\scriptstyle#1}\downarrow\phantom{\scriptstyle#1}$}}
\newcommand{\ldiagup}[1]%
       {\makebox[0cm]{${\scriptstyle#1}\uparrow\phantom{\scriptstyle#1}$}}
\newcommand{\rdiag}[1]%
       {\makebox[0cm]{$\phantom{\scriptstyle#1}\downarrow{\scriptstyle#1}$}}
\newcommand{\sediagr}[1]%
       {\makebox[0cm]{$\phantom{\scriptstyle#1}\searrow{\scriptstyle#1}$}}
\newcommand{\nediagr}[1]%
       {\makebox[0cm]{$\phantom{\scriptstyle#1}\nearrow{\scriptstyle#1}$}}
\newcommand{\rdiagup}[1]%
       {\makebox[0cm]{$\phantom{\scriptstyle#1}\uparrow{\scriptstyle#1}$}}
\newcommand{\swdiag}[1]%
       {\makebox[0cm]{$\phantom{\scriptstyle#1}\swarrow{\scriptstyle#1}$}}
\newcommand{\sediag}[1]%
       {\makebox[0cm]{${\scriptstyle#1}\searrow\phantom{\scriptstyle#1}$}}
\newcommand{\nediag}[1]%
       {\makebox[0cm]{${\scriptstyle#1}\nearrow\phantom{\scriptstyle#1}$}}

\newcommand{\longiso}{\stackrel{\textstyle\sim}{\longrightarrow}}

\newcommand{\doublearrowstack}[2]%
                      {{{{\scriptstyle#1}\atop{\textstyle\longrightarrow}}\atop{{\textstyle\longrightarrow}\atop{\scriptstyle#2}}}}
\newcommand{\rightleftarrowstack}[2]%
                      {{{{\scriptstyle#1}\atop{\textstyle\longrightarrow}}\atop{{\textstyle\longleftarrow}\atop{\scriptstyle#2}}}}
\newcommand{\leftrightarrowstack}[2]%
                      {{{{\scriptstyle#1}\atop{\textstyle\longleftarrow}}\atop{{\textstyle\longrightarrow}\atop{\scriptstyle#2}}}}

\newcommand{\overtoparrow}%
{\makebox[0cm]{\beginpicture
\setcoordinatesystem units <.8cm,.4cm> point at 0 0
\setplotarea x from -3 to 3, y from 0 to 1
\setquadratic
\plot -3 0 0 1 3 0 /
\put{\vector(3,-1){0}}[Bl] at 3 0
\endpicture}}

\newcommand{\underbottomarrow}%
{\makebox[0cm]{\beginpicture
\setcoordinatesystem units <.8cm,.4cm> point at 0 0
\setplotarea x from -3 to 3, y from 0 to 1
\setquadratic
\plot -3 1 0 0 3 1 /
\put{\vector(3,1){0}}[Bl] at 3 1
\endpicture}}

\newcommand{\ses}[5]%
{0\longrightarrow#1\stackrel{#2}{ \longrightarrow}#3\stackrel{#4}{
\longrightarrow}#5\longrightarrow0}

\newcommand{\dt}[6]%
{#1\stackrel{#2}{longrightarrow}#3 \stackrel{#4}{\longrightarrow}#5
\stackrel{#6}{\longrightarrow} #1[1]}  
 
\newcommand{\cat}[1]%
{(\mbox{\rm #1})}

\newcommand{\vfx}{V}
\newcommand{\vfy}{W}

\def\Label#1{\label{#1}{\tt [#1]}\phantom{h}}

\def\Label{\label}
\def\comment{\noprint}

\title{Differentiable Stacks and Gerbes}
\author{Kai Behrend\\ University of British Columbia\\
1984 Mathematics Road\\
Vancouver, B.C., Canada V6T 1Z2\\ \texttt{behrend@math.ubc.ca} \bigskip\\
  Ping Xu\thanks{Research partially supported by NSF grants DMS03-06665
 and DMS-0605725  \& NSA grant 03G-142.}  \\ 
Department of Mathematics \\ Pennsylvania State University \\ University Park, PA 16802 \\  \texttt{ping@math.psu.edu}}

\date{}
\begin{document}
\sloppy
\maketitle

\begin{abstract}
We introduce differentiable stacks and explain the relationship with
Lie groupoids.  Then we study $S^1$-bundles and $S^1$-gerbes over
differentiable stacks.  In particular, we establish the relationship
between $S^1$-gerbes and groupoid $S^1$-central extensions.
We define connections and curvings for groupoid $S^1$-central extensions
extending the corresponding notions of Brylinski,   Hitchin and Murray for
 $S^1$-gerbes over manifolds. We  develop a Chern-Weil
theory of characteristic classes in this general setting by
presenting  a  construction of  Chern classes and Dixmier-Douady classes
 in terms of analogues of connections and curvatures.
We also describe a   prequantization result for  both
$S^1$-bundles and $S^1$-gerbes extending the well-known
result of Weil and Kostant. In particular, we give
an explicit construction of $S^1$-central extensions with
prescribed curvature-like data.
\end{abstract}

\tableofcontents

\section{Introduction}

Grothendieck introduced stacks  to give geometric
meaning to higher non-commutative cohomology classes.
This is also the context in which {\em gerbes} first appeared \cite{Giraud}.
However  most of the work  on stacks so far remains
algebraic, though there is  increasing
evidence that  differentiable stacks
will find many useful applications.
One example of the notion of stack is that of orbifolds. In algebraic geometry,
these correspond to   Deligne-Mumford stacks \cite{LL}.
In differential geometry, orbifolds or $V$-manifolds
have been studied for many years using local charts.
Recently, it has been realized that viewing
orbifolds as  a  very special kind of
Lie groupoids, i.e. \'etale proper groupoids,  is quite useful \cite{Moer}.

The notion of a groupoid is a  generalization of the
concepts of {spaces} and {groups}. A groupoid
consists of a space of objects (units) $X_0$, and a space of arrows
$X_1$ with source and target maps $s, t: X_1\to X_0$.
 There is a  multiplication  defined only for
composable pairs $X_2=\{(x, y)\st \text{$t(x)=s(y)$, for $x,  y\in X_1$}\}
 \subset X_1\times X_1$. There is
also an inverse map. These structures satisfy
the usual axioms. Lie groupoids are groupoids where
both $X_0$ and $X_1$ are manifolds, $s$ and $t$ are surjective
submersions,  and all the structure maps are required to be smooth.
 A Lie groupoid
$X_1 \toto X_0$ is said to be {\em proper} if the map
$s\times t: X_1 \to X_0\times X_0$ is  proper (in algebraic geometry,
this would be called {\em separated or Hausdorff}).
In the theory of groupoids, spaces and groups are treated on equal footing.
Simplifying somewhat, one could say that a groupoid is a mixture of a space and
a group; it has space-like and group-like properties that interact in
a delicate way. In a certain sense, groupoids provide a uniform
framework for many different geometric objects.
For instance, when a Lie group acts on a manifold
properly, the corresponding equivariant cohomology theories, including
$K$-theory, can be treated using the transformation groupoid
$ M\rtimes G\toto M$. Here the structure maps are
$s(x, g)=x, \ t(x, g)=xg, \ \ (x, g)(y, h)=(x, gh)$.


There exists a  dictionary between differentiable stacks and
Lie groupoids. Roughly speaking,
differentiable  stacks  are
Lie groupoids {\em up to Morita equivalence}.
Any Lie groupoid $X_1\toto X_0$ defines a differentiable
stack $\XX$ of $X\lcom$-torsors.   Two differentiable stacks
 $\XX$ and  $\XX'$ are isomorphic if and only
if the Lie groupoids $X\lcom$ and $X\lcom'$ are Morita equivalent.
In a certain sense,   Lie groupoids are like ``local charts''
on   a differentiable stack. Establishing such a dictionary
consists of the first part of the paper. We note that
this viewpoint of connecting    stacks with groupoids is
somehow folklore 
 (see \cite{Duskin, Metzler, Noohi1, Noohi2}). However, we feel
that it is useful to spell it  out in detail in the 
differentiable geometry setting, which is of ultimate interest for
our purpose.  

Our main goal of this paper is to
develop the theory of  $S^1$-gerbes over differentiable  stacks.
Motivation  comes
from   string theory in which ``gerbes with connections'' appear
naturally \cite{Dijkgraaf, Freed, Kalkkinen, Sharpe}.

 For $S^1$-gerbes over manifolds, there has been extensive work
on this subject pioneered  by Brylinski \cite{Brylinski:book}, 
Chatterjee \cite{Chatterjee}, 
Hitchin \cite{Hitchin},  Murray \cite{Murray} and many others.
Also, there is  interesting work on equivariant
 $S^1$-gerbes, e.g., 
 by Brylinski \cite{Brylinski:gerbe},
 Meinrenken \cite{Meinrenken}, Gawedzki-Neis \cite{GN}, Stienon \cite{Stienon}
 and others, as well as on gerbes over orbifolds \cite{LU}.
These endeavors make the foundations of  gerbes over
 differentiable  stacks a very important  subject.
An important step is to  geometrically realize a class $H^2(\XX,  S^1)$
(or $H^3(\XX,  \zz)$ when $\XX$ is Hausdorff).
Such a   geometrical realization
is  crucial in applications to twisted $K$-theory \cite{TXL, TX1, TX2}.

Our method is to use the dictionary  mentioned above, under which we show that
 $S^1$-gerbes are in one-to-one correspondence with
Morita equivalence classes of groupoid $S^1$-central extensions.  Thus
it follows from a well-known theorem of Giraud \cite{Giraud}
that there is a  bijection
between $H^2(\XX, S^1)$ and  Morita equivalence 
classes of  Lie groupoid $S^1$-central extensions.
We  note that there are  several   independent investigations of
similar topics; see \cite{BSS, Moer:regular, Moerdijk, tu:05, WS}. 

An $S^1$-central extension of a Lie groupoid $X_1\toto X_0$
is a Lie groupoid $R_1  \toto X_0$ with a  groupoid morphism
$\pi: R_1\to X_1$ such that $\ker \pi \cong X_0\times S^1$  lies
in the center of $ R_1$.  It is easy to see
that $\pi: R_1\to X_1$  is then naturally an
$S^1$-principal bundle.
A standard example is an  $S^1$-central extension
 of a \v{C}ech groupoid:
Let $N$ be  a manifold and $\alpha \in H^3 (N, \zz )$, and
let  $\{U_i \}$ be  a  good  covering of $N$. Then  the groupoid
$\coprod_{ij} U_{ij}\toto \coprod_i U_i$, where $U_{ij}
=  U_i \cap U_j $,
 which is called the
\v{C}ech groupoid, is Morita equivalent to the manifold
$N$.
 Then the $S^1$-gerbe  corresponding to
the class $\alpha$ can be  realized as an $S^1$-central
extension of groupoids $\coprod_{ij} U_{ij} \times S^1\to \coprod_{ij} U_{ij}\toto
 \coprod_i U_i$, where the multiplication on $\coprod_{ij} U_{ij} \times S^1$
is given by
$(x_{ij}, \lambda_1) (x_{jk}, \lambda_2 )= (x_{ik}, \lambda_1  \lambda_2 c_{ijk})$,
where $x_{ij}, \ x_{jk}, x_{ik}$ are the same point $x$ in the
three-intersection $U_{ijk}$ considered as elements in the
two-intersections, and $c_{ijk}: U_{ijk}\to S^1$
 is a  2-cocycle which represents the
\v{C}ech class in $H^2 (N, {{S}}^1 ) \cong H^3(N, \zz )$
 corresponding to $\alpha$.
This is essentially the picture of an $S^1$-gerbe over a
manifold described by Hitchin \cite{Hitchin}.

The exponential sequence  $\zz\to \Omega^0\to S^1$ 
 induces a boundary map 
$H^2(X\lcom,S^1)\to H^3(X\lcom,\zz)$.  The image
in $H^3(X\lcom,\zz)$   of 
the class in $H^2(X\lcom,S^1)$ of  a groupoid $S^1$-central 
 extension  $R\lcom$  is called the {\em
Dixmier-Douady class }of $R\lcom$.  The Dixmier-Douady class behaves well
with respect to the  pullback and the tensor operations.
A fundamental   question is   to develop a Chern-Weil
characteristic class theory to construct the Dixmier-Douady classes
 geometrically.  For this purpose, we
need the  de Rham double complex of a Lie groupoid.
Let $X_1\toto X_0$ be a Lie groupoid.
By $X_{p}$, we denote the manifold of composable
sequences of $p$ arrows in the groupoid $X_1\toto X_0$. 
We have $p+1$ canonical maps, called face maps,
 $X_{p}\to X_{p-1}$ giving rise to a diagram
 \begin{equation}
\xymatrix{
\ldots X_{2}
\ar[r]\ar@<1ex>[r]\ar@<-1ex>[r] & X_{1}\ar@<-.5ex>[r]\ar@<.5ex>[r]
&X_{0}\,.}
\end{equation}
In fact, $X\lcom$ is a  simplicial manifold \cite{Dupont}.
Thus for any abelian sheaf $F$
  (e.g., $\zz$, $\rr$,  or $S^1$),
we have the cohomology groups $H^k(X\lcom,F)$.
Just like  for manifolds, $H^k(\XX, \rr)$ is canonically isomorphic
to the de Rham cohomology of $X_1\toto X_0$,
which is defined by  the double complex $\Omega\upcom(X\lcom)$, with
boundary maps $d:
\Omega^{k}( X_{p} ) \to \Omega^{k+1}( X_{p} )$, the usual exterior
derivative of differentiable forms, and $\partial
:\Omega^{k}( X_{p} ) \to \Omega^{k}( X_{p+1} )$,  the
alternating sum of the pull-back of face maps.
We denote the total differential by $\delta = (-1)^pd+\del$.
The  cohomology groups of the total complex
$H_{dR}^k(X\lcom)= H^k\big(\Omega\upcom(X\lcom)\big)$
are called the {\em de~Rham cohomology }groups of $X_1\toto X_0$.
When $X_1\toto X_0$ is the \v{C}ech groupoid associated
to an open covering of a manifold $N$, this
is isomorphic  to the usual de Rham cohomology of the manifold
$N$. On the other hand, when
 $X_1\toto X_0$ is a transformation groupoid $G\rtimes M\toto M$, then
$H_{dR}^k(X\lcom)$ is isomorphic to  the equivariant cohomology
$H^k_G (M)$.

Following Murray  \cite{Murray} and  Hitchin \cite{Hitchin}, for a given
groupoid  $S^1$-central extension,
one can also define the notions of connections,
curvings and $3$-curvatures. A {\em flat} gerbe is
one whose $3$-curvature vanishes. In this case,
there exists a holonomy map as well. However, a substantial
difference between $S^1$-gerbes over an arbitrary differential
stack and  $S^1$-gerbes over a manifold is that
connections and curvings may not always exist.
Therefore they may not be as useful as one expects in calculating
Dixmier-Douady classes. For this purpose, we need the
notion of so called pseudo-connections.
Given an $S^1$-central extension $R_1\to X_1 \toto X_0$,
a {\em pseudo-connection } consists of a pair
$(\theta,B)$, where  $\theta\in\Omega^1(R_1 )$ is
 a connection  1-form for the $S^1$-principal  bundle
$R_1\to X_1$, and $B\in\Omega^2(X_0)$ is a 2-form.
It is simple to check that $\delta(\theta+B) \in Z^3_{dR}(R\lcom)$
 descends to $Z^3_{dR}(X\lcom)$, i.e. there exist unique
$\eta\in\Omega^1(X_2)$, $\omega\in\Omega^2(X_1)$ and
$\Omega\in \Omega^3(X_0)$ such that
$  \delta(\theta+B)= \pi\upst(\eta+\omega+\Omega)$.
Then $\eta+\omega+\Omega$ is called the {\em  pseudo-curvature }of
the pseudo-connection $\theta+B$. It is simple to
check that the class $[\eta+\omega+\Omega]\in H^3_{dR}(X \lcom)$
is independent of the
choice of the pseudo-connection $\theta+B$.
One of the    main results of this paper is that
$[\eta+\omega+\Omega]$ is indeed  the Dixmier-Douady class, or more
precisely,  the image of the Dixmier-Douady class
  under the canonical
homomorphism $H^3(\XX,\zz)\to H^3(\XX,\rr)\cong H^3_{dR}(X\lcom)$.
Recently,  Ginot-Stienon found an
alternative proof of this result using $2$-group bundles \cite{GS}
(in fact they proved a more general result for central $G$-extensions).
We also describe a prequantization result, an analogue of the
Kostant-Weil \cite{Kostant, Weil} theorem for $S^1$-gerbes.
That is, given any integral   3-cocycle
 $\eta +\omega+\Omega\in Z^3_{dR}(X \lcom)$,  we describe a
sufficient condition that guarantees the 3-cocycle
as   the pseudo-curvature
of a groupoid  $S^1$-central extension  $R\lcom$
with a pseudo-connection
$\theta +B$,  and classify all such  pairs $(R\lcom, \theta+ B)$.

$S^1$-central extensions of  Lie groupoids
also appear naturally in Poisson geometry. It was proved in
\cite{WeinsteinXu} that  a certain  prequantization of a symplectic
groupoid  naturally becomes  an $S^1$-central extension of groupoids
with a connection, which is indeed  a contact groupoid.
The proof utilizes Lie algebroids as a tool.
Lie algebroids are infinitesimal versions  of Lie groupoids.
It is thus natural to  study Lie groupoid central extensions
via Lie algebroid central extensions in a general framework.
More precisely, let $X_1\toto X_0$ be an $s$-connected Lie groupoid 
with Lie algebroid $A$, and let  $\eta +\omega
\in Z^3_{dR} (X\lcom )$ be  a de-Rham 3-cocycle,
 where $\eta \in \Omega^1(X_2)$ and $\omega \in \Omega^2 (X_1 )$. Then
$\omega -d\eta^r\in \Omega^2 (X_1^t)$ is
a right invariant $t$-fiberwise closed two-form
on $X_1$, and therefore  defines a Lie algebroid two-cocycle
of $A$, which in turn defines 
a Lie algebroid  central extension $\tilde{A}=A\oplus (X_0\times \rr )$ 
of $A$.  Here $\eta^r$ is a $t$-fiberwise
   one-form     on $X_1$ given  by
$\eta^r  (\delta_{x})=\eta (r_{x^{-1}*}\delta_{x} , \  0_{x} ),
\forall \delta_x \in T_x X_1^t$,
and  $r_{x^{-1}}$  denotes the right  translation.
A natural  question is: under what condition
does this Lie algebroid central extension give rise to
a Lie groupoid central extensione? 
The last part of the paper is devoted to investigating this
question. Our method is  to adapt  the method
of characteristics developed  by Coste-Dazord-Weinstein \cite{CDW}.
As a consequence,  we obtain a geometrical characterization of the integrality
 condition of a de Rham $3$-cocycle $\eta +\omega
\in Z^3_{dR} (X\lcom )$ of a Lie  groupoid $X_1\toto X_0$.

The results of this paper were announced in \cite{BX}. See
also \cite{BXZ} for a construction of $S^1$-gerbes over
the quotient stack  $[G/G]$ ($G$ is
a compact simple Lie group and $G$ acts on
$G$ by conjugations) as an example.

{\bf Acknowledgments.}
We  would like to thank several institutions
for their hospitality while work on this project was being done:
 RIMS and IHP (Behrend and Xu), University of British Columbia and
 Universit\'e Pierre et Marie Curie (Xu).
 We  also  wish to thank many people for  useful discussions and comments,
including  Camille Laurent-Gengoux, Jim Stasheff, Mathieu Stienon,
Jean-Louis Tu and Alan Weinstein.

\section{Differentiable Stacks}

Our goal in this section is to define the notion of differentiable
stack and establish a dictionary between differentiable stacks and
Lie groupoids.  Roughly speaking, differentiable stacks are
Lie groupoids {\em up to Morita equivalence}. 

Our differentiable manifolds will {\em not }be assumed to necessarily
be Hausdorff. We use the words $\cinf$ and {\em smooth
}interchangeably. The manifold consisting of one point is denoted by
$\ast$ or $\pt$. 

Let us start by recalling some terminology.
A $\cinf$-map $f:U\to X$ of $\cinf$ manifolds is 
a {\em submersion}, if for all $u\in U$ the derivative
$f\lst:T_uU\to T_{f(u)}X$ is surjective.  The {\em relative
dimension }of 
the submersion $f$ is the (locally constant on $U$) dimension of the
kernel of $f\lst$.  A submersive map of relative dimension $0$ is called {\em
\'etale}.  Thus $f$ is \'etale if and only if it is a local
diffeomorphism.

Let $\SS$ be the category of all $\cinf$-manifolds with $\cinf$-maps as
morphisms. Note that not all fiber products exist in $\SS$, but if at
least one of the two morphisms $U\to X$ or $V\to X$ is submersive, then
the fiber product $U\times_X V$ exists in $\SS$. 
In general, the fiber product $U\times_X V$ exists 
if $U\to X$ and $V\to X$ satisfy the transversality condition.

We endow $\SS$ with the Grothendieck topology given by the following
notion of covering family.  Call a family $\{U_i\to X\}$ of morphisms
in $\SS$ with target $X$ a {\em covering family of $X$}, if all maps
$U_i\to X$ are \'etale and the total map $\coprod_i U_i\to X$ is
surjective.

One checks that the conditions for a Grothendieck topology (see
Expos\'e~II in \cite{sga4}) are satisfied.  (Note that in the terminology of
loc. cit.  we have actually defined a {\em pretopology}. This
pretopology gives rise to a Grothendieck topology, as explained in
loc.\ cit..) We call this topology the {\em \'etale }topology on
$\SS$.

One can also work with the topology of open
covers. In this topology, all covering families are open covers $\{U_i\to
X\}$, in the usual topological sense.  The notion of sheaf or stack
over $\SS$ obtained using this topology is the same as using the \'etale
topology.

A {\em site }is just a category endowed with a Grothendieck
topology. So if we refer to $\SS$ as a site, we emphasize that we
think of $\SS$ together with its \'etale topology.

A {\em Lie groupoid }is a groupoid in $\SS$, whose source and target
maps are submersions.

\subsection{Groupoid fibrations}

 A {\em category fibered in groupoids }$\XX\to\SS$ is a category $\XX$,
together with a functor $\pi:\XX\to\SS$, such that the following two
fibration axioms are satisfied:

(i) for every arrow $V\to U$ in $\SS$, and every object $x$ of $\XX$
lying over $U$ (i.e. $\pi(x)=U$), there exists an arrow $y\to x$ in
$\XX$ lying over $V\to U$,

(ii) for every commutative triangle $W\to V\to U$ in $\SS$ and arrows
$z\to x$ lying over $W\to U$ and $y\to x$ lying over $V\to U$, there
exists a unique arrow $z\to y$ lying over $W\to V$, such that the
composition $z\to y\to x$ equals $z\to x$.

The object $y$ over $V$, whose existence is asserted in (i), is unique
up to a unique isomorphism by (ii).  Any choice of such a $y$ is called
a {\em pullback }of $x$ via $V\to U$, notation $y=x\resto V$, or
$y=f\upst x$, if the morphism $V\to U$ is called $f$.
Often it is convenient to choose pullbacks for all $x$ and all $V\to
U$ (where $U=\pi(x)$).

Given a category fibered in groupoids $\XX\to\SS$ and an object $U$ of
$\SS$, the category of all objects of $\XX$ lying over $U$ and all
morphisms of $\XX$ lying over $\id_U$ is called the {\em fiber }of
$\XX$ over $U$, notation $\XX_U$, sometimes $\XX(U)$.  Note that
all fibers $\XX_U$ are (set-theoretic) groupoids.  This follows from
Property~(ii), above.

We call {\em categories fibered in groupoids }over $\SS$ also simply
{\em groupoid fibrations}. The groupoid fibrations over $\SS$ (see
\cite{sga1}) form a 2-category.
Fibered products exist. They satisfy a 2-categorical version of the
universal mapping property for fibered products  (see \cite{LL}).

The notion of groupoid fibration is  the mathematical formalization
of the notion of {\em moduli problem}. 
Let $\XX\to\SS$ be  a groupoid fibration.  If we consider $\XX$ as a
moduli problem, then we think of an object $x\in\XX$ lying
over $S\in\SS$ as an {\em $\XX$-family parametrized
by $S$}.
The objects we wish to classify are the objects of the
category $\XX(\pt)$. 

Standard examples of categories fibered in groupoids over $\SS$ are:

\begin{numex}\label{bg}
Let $G$ be a Lie group. Let $\XX=BG$ be the category of pairs $(S,P)$,
where $S\in\SS$ is a $\cinf$-manifold and $P$ is a principal
$G$-bundle over $S$. A morphism from $(S,P)$ to $(T,Q)$ is a
commutative diagram
$$\xymatrix{
P\dto\rto & Q\dto\\
S\rto & T}$$
where $P\to Q$ is $G$-equivariant.  The functor $\pi:BG\to\SS$ is
defined by $(S,P)\mapsto S$.
\end{numex}

\begin{numex}
Every manifold $X$ defines a groupoid fibration $F_X$ over $\SS$.
The objects of $F_X$ are pairs $(U, f)$, where $U$ is a
$C^\infty$-manifold and $f: U\to X$ is a smooth map.
Morphisms in $F_X$ are the  commutative triangles
$$\xymatrix{U\rto \drto& V\dto\\&X}$$
The  functor $F_X\to \SS$ is the  projection onto the first component.
The groupoid fibration  $F_X$ satisfies
$$F_X(U)=\Hom_{\SS}(U,X)\,.$$
By abuse of notation, we identify $F_X$ with $X$ in the sequel.
\end{numex}

\begin{numex}
Let $\MM_g$ be the following groupoid fibration: objects
are fiber bundles $X\to S$ endowed with
a smoothly varying fiberwise complex structure, 
such that all fibers are Riemann surfaces of genus $g$.
Morphisms are commutative diagrams
$$\xymatrix{X\rto \dto&Y\dto\\S\rto&T}$$
such that $X\to Y\times_T S$ is  a conformal isomorphism.
This   groupoid fibration is the  moduli stack of Riemann surfaces of
genus $g$. An object $X\to S$ of $\MM_g$ is a family
of Riemann surfaces parametrized by $S$. The functor
$\MM_g\to \SS$ maps $X\to S$ to $S$.
\end{numex}

\begin{numex}\label{cf}
Any contravariant functor $F:\SS\to(\text{sets})$ gives rise to a
category fibered in groupoids $\XX\to\SS$ defined as follows: objects
of $\XX$ are pairs $(U,x)$, where $U$ is a $\cinf$-manifold and $x\in
F(U)$.  A morphism $(U,x)\to(V,y)$ is a $C^\infty$ map $a: U\to V$
such that $F(a)(y) =x$. The functor $\pi:\XX\to\SS$ is defined by
$(U,x)\mapsto U$.

In particular, a sheaf over $\SS$ defines a groupoid fibration
over $\SS$ in a canonical way.
\end{numex}

\begin{defn}
A groupoid fibration $\XX$ over $\SS$ is {\bf representable}, if there exists a
manifold $X$ such that $X\cong\XX$ (as groupoid fibrations over $\SS$).
\end{defn}

\begin{defn}
A morphism of groupoid fibrations $\XX\to\YY$ is called a {\bf representable
submersion}, if for 
every manifold $U$ and every morphism $U\to\YY$ the fibered product
$V=\XX\times_{\YY}U$ is representable and the induced morphism
of manifolds $V\to U$ is a submersion.

If the relative dimension $V\to U$ is always equal to $n\in\zz$, then
we call $n$ the {\bf relative dimension }of $\XX\to\YY$.
\end{defn}

\begin{numex}\label{falg}
For a Lie group $G$, the canonical morphism $\ast\to BG$ 
is a representable submersion. Here the functor
assigns to any smooth manifold $U$  the trivial
$G$-bundle over $U$. We can think of $\ast\to BG$ as the universal
$G$-bundle, because every $G$-bundle $P\to S$ gives rise to a
2-fibered product
$$\xymatrix{
P\rto\dto & S\dto\\
\ast\rto& BG}$$
\end{numex}

The following lemma will be useful in the future.

\begin{lem}[Descent]
\label{descent}
Let $F$ be a sheaf over $\SS$. Let $X$ be a manifold and $F\to X$
a morphism. Suppose that $\{U_i\to X\}$ is a covering family of $X$
and that for every $i$ the sheaf $F_i=U_i\times_X F$ is
representable. Then $F$ is representable.
\end{lem}
\begin{pf}
First note that we can choose a refinement of the covering $\{U_i\to
X\}$ consisting of open subsets of $X$. Replacing the covering by such
a refinement, we reduce to the case of a cover $\{U_i\to X\}$
consisting of open subsets.

Let, as usual, $U_{ij}=U_i\times_X U_j=U_i\cap U_j$. Define
$F_{ij}=U_{ij}\times_XF$. Then all $F_{ij}$ are
representable.
 Moreover, all
maps $F_{ij}\to F_i$ and $F_{ij}\to F_j$
are (isomorphic to) embeddings of open subsets.  Thus we can glue the
manifolds $F_i$ along the open submanifolds $F_{ij}$ to obtain a
manifold representing $F$.
\end{pf}

\begin{defn}
\label{epi}
A morphism of groupoid fibrations $\XX\to \YY$ is an {\bf epimorphism}
if for every $U\to \YY$, where $U$ is a manifold, there exists a
surjective submersion $V\to U$ and a 2-commutative diagram
$$\xymatrix{V\rto\dto\drtwocell\omit & U\dto\\ \XX\rto & \YY}$$
Equivalently, $V$ may be replaced by an open cover of $U$, in this
statement.
\end{defn}

\begin{rmk}\label{stup.ex}
Let $\XX$ be a category fibered in groupoids over $\SS$.
Given a manifold $U\in\SS$ and an object $x\in\XX_U$
(we write $x|U$), 
 the choice of
pullbacks of $x$ for all maps $V\to U$ defines a morphism $U\to \XX$.
Conversely, given a morphism $U\to \XX$, the image of $\id_U$ is an
object in the fiber $\XX_U$.  In this way we identify morphisms $U\to
\XX$ with objects in the fiber $\XX_U$.
\end{rmk}

\subsection{Stacks}
Recall the definition of stack \cite{LL}:

\begin{defn}
Let $\XX\to\SS$ be a category fibered in groupoids. 
  We call $\XX$
a {\bf stack }over $\SS$, if the following three  axioms are
satisfied:

(i) for any $\cinf$-manifold $X\in\SS$, any two objects $x,y\in
\XX$ lying over $X$, and any two isomorphisms $\phi,\psi:x\to y$ over
$X$, such $\phi\resto U_i=\psi\resto U_i$, for all $U_i$ in a covering
family $U_i\to X$, we have that $\phi=\psi$;

(ii) for any $\cinf$-manifold $X\in\SS$, any two objects $x,y\in
\XX$ lying over $X$, a covering family $U_i\to X$ and, for every $i$,
an isomorphism $\phi_i:x\resto U_i\to y\resto U_i$, such that
$\phi_i\resto U_{ij}=\phi_j\resto U_{ij}$, for all $i,j$, there exists
an isomorphism $\phi:x\to y$, such that $\phi\resto U_i=\phi_i$, for
all $i$;

(iii) for every $\cinf$-manifold $X$, every covering family $\{U_i\}$
of $X$, every family $\{x_i\}$ of objects $x_i$ in the fiber
$\XX_{U_i}$ and every family of morphisms $\{\phi_{ij}\}$,
$\phi_{ij}:x_i\resto U_{ij}\to x_j\resto U_{ij}$, satisfying the
cocycle condition $\phi_{jk}\comp\phi_{ij}=\phi_{ik}$ (which is an
equation in the fiber $\XX_{U_{ijk}}$), there exists an object $x$ over
$X$, together with isomorphisms $\phi_i:x\resto U_i\to x_i$ such that
$\phi_{ij}\comp\phi_i=\phi_j$ (over $U_{ij}$).
\end{defn}

Note that the isomorphism $\phi$, whose existence is asserted in (ii)
is unique, by (i). Similarly, the object $x$, whose existence is
asserted in (iii), is unique up to a unique isomorphism, because of
(i) and (ii).  The object $x$ is said to be obtained by {\em gluing
}the objects $x_i$ according to the gluing data $\phi_{ij}$.

Note also that there are choices to be made for all the pullbacks
mentioned in the definition of stacks, but no property depends on any
of these choices.

\begin{rmk}
To any covering family $U_i$ of $X$, we  can associate a  
groupoid fibration $R$, together with a monomorphism $R\subset X$ of
groupoid fibrations, the {\em covering sieve }given by $U_i$. The
stack axioms may be reformulated in terms of covering sieves: thus, a
groupoid fibration $\XX$ is a stack if and only if for every covering
sieve $R\subset X$, of every object $X\in\SS$, the functor
\begin{equation}\Label{sieve}
\bhom_\SS(X,\XX)\longrightarrow\bhom_\SS(R,\XX)
\end{equation}
is an equivalence of groupoids.  More precisely, $\XX$ satisfies Stack
Axiom~(i) if and only if (\ref{sieve}) is always faithful, $\XX$
satisfies Stack Axiom~(ii) if and only if (\ref{sieve}) is always
full and $\XX$ satisfies Stack Axiom~(iii) if and only if
(\ref{sieve}) is always essentially surjective.
\end{rmk}

The following  lemma is useful in practice.

\begin{lem}\label{critrep}
Let $f:\XX\to\YY$ be a morphism of stacks over $\SS$. Suppose given
a manifold $U$ and a morphism $U\to \YY$ which is an {\em epimorphism}. If
the fibered product $V=\XX\times_{\YY}U$ is representable and $V\to U$
is a submersion, then $f$ is a representable submersion.

If $V\to U$ has relative dimension $n$, then so does $f$.
\end{lem}
\begin{pf}
Let $W\to\YY$ be an arbitrary morphism, where $W$ is a manifold. First
we have to show that the fibered product $F=\XX\times_{\YY}W$ is
representable. By the fact that $U\to\YY$ is an epimorphism, we can
choose a covering family $\{W_i\to W\}$ of $W$ and morphisms
$\phi_i:W_i\to U$ making the diagram
$$\xymatrix{
W_i\dto\rto^{\phi_i} \drtwocell\omit{^}& U \dto\\
W\rto & \YY}$$
commute (which involves, of course, also a choice of a 2-arrow, for
every $i$). By Lemma~\ref{descent}, it suffices to prove that
$F_i=W_i\times_WF$ is representable, for all $i$. But
$F_i=W_i\times_UV$, as can be seen from the cartesian cube
$$\xymatrix@=10pt{
F_i\ddto\rrto\drto&& V\ddto\drto & \\
& W_i\ddto\rrto && U\ddto\\
F\drto\rrto &&\XX \drto&\\
& W \rrto &&\YY}$$
and so, indeed, $F_i$, and hence $F$ is representable.

Now the fact that $F\to W$ is a submersion, follows from the fact that
for every $i$ the map $F_i\to W_i$ is a submersion, because being a
submersion is a local property. But $F_i\to W_i$ is a submersion as a
pull back of the submersion $V\to U$.
\end{pf}

\begin{numex}
Let $G$ be a Lie group and $H$ a closed Lie subgroup. The induced morphism
$BH\to BG$ is a representable submersion. To see this, let us apply
Lemma~\ref{critrep}.  Note that $\ast\to BG$ is an epimorphism, 
because every $G$-bundle is locally trivial. Note also that we have a
cartesian diagram 
$$\xymatrix{
G/H\dto\rto & \ast\dto\\
BH\rto & BG}$$
because the reductions of structure group from $G$ to $H$ of the
trivial $G$-bundle over a manifold $U$ are classified by the maps
$U\to G/H$. Since $G/H$ is a manifold, $G/H\to \ast$ is a submersion,
which finishes the proof.
 The relative dimension of $BH\to BG$ is
equal to $\dim G-\dim H$.
\end{numex}

Two stacks $\XX$ and $\YY$ over $\SS$ are said
to be  {\em isomorphic }if they are
equivalent as categories over $\SS$. This means that there exist
morphisms $f:\XX\to \YY$ and $g:\YY\to\XX$ and 2-isomorphisms
$\theta:f\comp g\Rightarrow \id_{\YY}$ and $\eta:g\comp
f\Rightarrow\id_{\XX}$.

\begin{prop}\label{crit.iso}
For stacks $\XX$ and $\YY$ over $\SS$ to be isomorphic, it
suffices that
 there  exists a morphism $f:\XX\to \YY$ satisfying the two
conditions:

(i) for any two objects $x,x'$ of $\XX$, lying over the same object $U$ of
$\SS$, and any arrow $\eta:f(x)\to f(x')$ in $\YY_U$, there exists a unique
arrow $x\to x'$ mapping to $\eta$ under $f$ (we say $f$ is {\em fully
faithful }or a {\em monomorphism\/});

(ii) for every object $y$ of $\YY$, lying over $S\in\SS$, there exists
a covering family $\{U_i\}$ of $S$ and objects $x_i$ of $\XX$ lying
over $U_i$, such that $f(x_i)\cong y\resto U_i$, for all $i$ (we say
that $f$ is an {\em epimorphism\/}).
\end{prop}

A morphism satisfying both these conditions is called an {\em
isomorphism }of stacks.

\subsection{Differentiable stacks}

Let $\XX$ be a  groupoid fibration over $\SS$.
Recall that we may think of $x/S$ equivalently as a
morphism of groupoid fibrations $x:S\to\XX$.

For $x/S$ and $y/T$, 
 consider the fibered product
$$\xymatrix{{\Isomu(x,y)}\rto\dto\drtwocell\omit{^} & T\dto^y\\
S\rto^x & \XX}$$
For an $\XX$-family $x$, parametrized by $S$, we call
$\Isomu(x,x)\toto S$ the {\em symmetry groupoid }of $x$. A priori,
$\Isomu(x,x)$ is just a groupoid fibration over $\SS$, but it may be
hoped that it is (represented by) a Lie groupoid. Note that we have a
cartesian diagram
$$\xymatrix{{\Isomu(x,x)}\rto\dto\drtwocell\omit{^} & {S\times
      S}\dto\\
{\XX}\rto^{\Delta} & {\XX\times\XX}}$$
Thus, ultimately, properties of the diagonal $\Delta:\XX\to\XX\times\XX$,
will assure that the symmetry groupoids $\Isomu(x,x)$ are manifolds, at
least if $S\to\XX$ is sufficiently well-behaved.


\begin{lem}
Let $f:\XX\to\YY$ be a representable submersion of stacks over
$\SS$. Then the following are equivalent:

(i) $f$ is an epimorphism;

(ii) for every manifold $U\to\YY$ the submersion $V\to U$, where $V$
is the fibered product $V=\XX\times_{\YY}U$, is surjective;

(iii) for some manifold $U\to\YY$, where $U\to\YY$ is an epimorphism,
the submersion $V\to U$ is surjective.

A representable submersion satisfying these conditions is called a
{\em surjective }representable submersion.
\end{lem}
\begin{pf}
This follows from the fact that a submersion between manifolds is an
epimorphism of stacks if and only if it is surjective. We also use
that to be an epimorphism is a local property.
\end{pf}

\begin{defn}
\label{diff.stack}
A stack $\XX$ over $\SS$ is called {\bf differentiable }or a {\bf
$\cinf$-stack}, if there exists a manifold $X$ and a surjective
representable submersion $x: X\to\XX$.  Such a manifold $X$, together
with the structure morphism $X\to\XX$ is called a {\bf presentation of
$\XX$ }or an {\bf atlas for $\XX$},  and
such a family $x/X$ is called a {\em versal family}.
\end{defn}

Alternatively, one can describe a  differentiable stack
in a slightly weaker condition.

\begin{prop}
\label{prop:equiv}
A stack over $\SS$ is a  differentiable stack, if there exists
an $\XX$-family $x/X$, such that 

(i) the symmetry groupoid $\Isomu(x,x)$ is 
representable, and the projections $\Isomu(x,x)\to X$ are  submersions;

(ii) the morphism $x: X\to \XX$ is an epimorphism. I.e., 
for every $\XX$-family $y/S$, there exists a covering family
$U_i$ of $S$, and morphisms $\phi_i:U_i\to X$, such that $y\resto
U_i\cong\phi_i\upst x$. 
\end{prop}
\begin{pf}
Given such an $\XX$-family $x/X$, it
suffices to show that $x: X\to \XX$ is  representable submersion.
This follows from Lemma \ref{critrep} since
$x: X\to \XX$ is epimorphism,  $X\times_{\XX} X$ is  
 representable and $X\times_{\XX} X\to X$ is a submersion.

The converse is obvious.
\end{pf}

The 2-category of differentiable stacks is the full sub-2-category of
the 2-category of groupoid fibrations over $\SS$ consisting of
differentiable stacks.

Given a  differentiable stack $\XX$, a versal family $x/X$
 gives rise to a  Lie groupoid $\Isomu(x,x)\toto X$ in a canonical
way.
 The points of $\Isomu(x,x)$ are by definition 
triples $(y,\phi,y')$, where $y$ and $y'$ are
points of $X$ and $\phi: x|y \to x|y'$ is a morphism in the groupoid
$\XX_\ast$ (the fiber of $\XX$ over $\ast\in\SS$).  So it is clear how
to define the composition: 
\begin{equation}\label{form1}
(y,\phi,y')\comp(y',\psi,y'')=(y,\psi\comp\phi,y'')\,.
\end{equation}
To see that this, indeed, defines the structure of a Lie groupoid on
$\Isomu(x,x)\toto X$, the quickest way is to note that for every manifold $U$,
evaluating at $U$ we get a (set-theoretic) groupoid $\Isomu(x,x)(U)\toto X(U)$,
defined by the same formula~(\ref{form1}) and compatible with all maps
$V\to U$.

A morphism  $\XX\to \YY$ of differentiable stack is 
{\bf representable} if one  of the following
equivalent conditions is satisfied:
\begin{enumerate}
\item there is a presentation $Y\to \YY$ such that 
$\XX\times_\YY Y$ is representable;
\item for any  representable submersion $Y\to \YY$,
$\XX\times_\YY Y$ is representable.
\end{enumerate}
For instance, the diagonal map $\XX\to \XX\times \XX$ of
a differentiable stack $\XX$ is always representable.

A representable morphism  $\XX\to \YY$ is called
{\bf proper} if there exists a 
presentation $Y\to \YY$ such that the base change
$X\to Y$ is proper. If this is the case, $X\to Y$ is
proper for all representable submersion $Y\to \YY$.
\begin{numex}
Let $\XX_g$ be the following groupoid fibration: objects
are fiber bundles $X\to S$ with fibers being isomorphic to 
a fixed  connected surface $Y$ of genus $g$.
Morphisms are commutative diagrams
$$\xymatrix{X\rto \dto&Y\dto\\S\rto&T}$$
such that $X\to Y\times_T S$ is  an  isomorphism. 
Consider the constant  family $Y\to *$. Then 
$*\stackrel{Y}{\to} \XX_g$ is an epimorphism, because every family of surfaces
is locally trivial. But $*\stackrel{Y}{\to}
\XX_g$ is not a representable submersion since the 
symmetry groupoid of this family is the diffeomorphism
group of $Y$, which is not a finite dimensional manifold.
 So $\XX_g$ is not  a differentiable stack.
\end{numex}

\subsection{Torsors for Lie groupoids}

Next,  we show  how to get a
differentiable stack starting from a Lie groupoid.  (This is, in fact,
 a generalization of passing from $G$ to $BG$.)

\begin{defn}
Let $\Gamma\toto M$ be a Lie groupoid and $S$ a manifold. A {\bf
$\Gamma$-torsor over $S$} is a manifold $P$, together with a
surjective submersion $\pi:P\to S$ and a (right) action of $\Gamma$ on
$P$, such that for all $p,p'\in P$, such that $\pi(p)=\pi(p')$, there
exists a unique $\gamma\in\Gamma$, such that $p\cdot\gamma$ is defined
and $p\cdot\gamma=p'$. 
\end{defn}

We call  the  map $P\to M$ of the $\Gamma$-torsor $P$ the {\em
anchor map }and denote it by  $a:P\to M$. 
(In the theory of symplectic
groupoids the anchor map is also called the  ``momentum map" \cite{MW}.)
And the  surjective submersion $\pi:P\to S$ is called the
{\em structure map}.

\begin{rmk}
Think of a $\Gamma$-torsor as follows.  View an element $p\in P$ as an
arrow eminating at $\pi(p)$ and terminating at $a(p)$.  Then view the
action of $\Gamma$ on $P$ as composing arrows. 
\end{rmk}

\begin{defn}
Let $\pi:P\to S$ and $\rho:Q\to T$ be $\Gamma$-torsors. A {\bf
morphism }of $\Gamma$-torsors from $Q$ to $P$ is given by a
commutative diagram of differentiable maps
\begin{equation}\label{morph}
\vcenter{\xymatrix{ 
Q\dto\rto^\phi & P\dto \\
T\rto & S}}
\end{equation}
such that $\phi$ is $\Gamma$-equivariant.
\end{defn}

Note that for a morphism of $\Gamma$-torsors the diagram (\ref{morph})
is necessarily a pullback diagram.

The $\Gamma$-torsors form a category with respect to this notion of
morphism. In particular, we now know what it means for two
$\Gamma$-torsors to be isomorphic.

\begin{numex}[trivial torsors]
Let $f:S\to M$ be a smooth map. Given $f$, we can induce over $S$ in a
canonical way a $\Gamma$-torsor, which we call the {\em trivial
}$\Gamma$-torsor given by $f$. 

Simply define $P$ to be the fibered product
$P=S\times_{f,M,s}\Gamma$. The structure map $\pi:P\to S$ is the first
projection. The anchor  map of the $\Gamma$-action is the second
projection followed by the target map $t$. The action is then defined
by 
$$(s,\gamma)\cdot\delta=(s,\gamma\cdot\delta)\,.$$
One checks that this is, indeed, a $\Gamma$-torsor over $S$.

Of course, we can take $S=M$ and $f$ the identity map of $M$. Then we
get the {\em universal }trivial $\Gamma$-torsor, whose base is
$M$. The structure morphism and the anchor map of the universal
$\Gamma$-torsor are, respectively,  $t, s:\Gamma\to M$.

Let $\pi:P\to S$ be an arbitrary $\Gamma$-torsor over the manifold
$S$. One checks that every section $s:S\to P$ of $\pi$ can be used to
construct an isomorphism between the $\Gamma$-torsor $P$ and the
trivial $\Gamma$-torsor over $S$ given by $a\comp s$, where $a:P\to M$
is the anchor map of $P$. 

Since every surjective submersion admits local sections, we see that
every $\Gamma$-torsor is {\em locally trivial}. 
\end{numex}

Let us denote the category of $\Gamma$-torsors by $B\Gamma$. There is a
canonical functor $B\Gamma\to\SS$ given by mapping a torsor $P\to S$
to the underlying manifold $S$.

The following proposition provides us with plenty of examples of
differentiable stacks. Theorem~\ref{torsors=stack} below
indeed  shows that
it provides us with {\em all }examples of differentiable stacks.

\begin{prop}
For every Lie groupoid $\Gamma\toto M$, the category of
$\Gamma$-torsors $B\Gamma$ is a differentiable stack.
\end{prop}
\begin{pf}
The fact that $B\Gamma$ is fibered in groupoids over $\SS$ follows from
the fact that diagrams such as (\ref{morph}) are always cartesian.
Note that given a $\Gamma$-torsor $P\to S$ and a
morphism of manifolds $T\to S$, $T\times_S P\to T$
is naturally a $\Gamma$-torsor over $T$.

To check the stack axioms, one has to prove that one can glue together
$\Gamma$-torsors and morphisms of $\Gamma$-torsors.  This is
 rather standard and will be omitted.

Finally we need to prove that $B\Gamma$ admits a presentation.
For this, we take the universal trivial torsor. 
 We shall construct a morphism $M\to B\Gamma$.
This means defining for every manifold $S$ a map $M(S)\to B\Gamma(S)$.
This we do by assigning to any smooth map $a:S\to M$ (i.e. object of
$M(S)$) the trivial $\Gamma$-torsor over $S$, which is an object of
$B\Gamma(S)$. Alternatively, we can use the universal trivial
$\Gamma$-torsor, which gives rise to the morphism $M\to B\Gamma$
directly, via the correspondence between objects of the fiber $ B\Gamma (U)$
and morphisms $U\to B\Gamma (U)$ (see the remark following
Definition \ref{epi}).

Now that we have a morphism $M\to B\Gamma$ from a manifold $M$, it
remains to prove that this morphism is a surjective representable
submersion.
To prove that $M\to B\Gamma$ is an epimorphism, means proving that
every $\Gamma$-torsor is locally trivial. This we have done already.
By Proposition \ref{prop:equiv},
 it now suffices to prove that the fibered
product
$$\xymatrix{
X\dto\rto\drtwocell\omit{^}& M\dto\\
M\rto& B\Gamma}$$
is representable and that the maps $X\toto M$ are submersions.

Let $S$ be an arbitrary manifold. Then $X(S)$ is the set of
triples $(a,\gamma,b)$, where $a,b:S\to M$ are
$C^\infty$  maps and $\gamma:Q_a\to Q_b$ is a morphism of
$\Gamma$-torsors over $S$, where $Q_a$ and $Q_b$ are the trivial
$\Gamma$-torsors over $S$ given by $a$ and $b$, respectively. One
checks that this set is canonically identified with $\Gamma(S)$, the
set of $C^\infty$-maps from $S$ to $\Gamma$.

Thus we have that $X\cong\Gamma$ as stacks over $\SS$, and so $X$ is
representable.  To check that the two projections $X\to M$ are
submersions, note that they are identified   with $s,t:\Gamma\to M$,
under this isomorphism.  Since $s$ and $t$ are submersions, we are done.
\end{pf}


\begin{rmk}
(1) Note that in the course of the proof we have seen that we have a
cartesian diagram
$$\xymatrix{
\Gamma\dto_s\rto^t\drtwocell\omit{^}& M\dto\\
M\rto& B\Gamma}$$
Thus $\Gamma\toto M$ is (isomorphic to) the Lie groupoid
arising from the atlas $M\to B\Gamma$.

(2) From the above proof, we see that for a given
$a:S \to M$ the corresponding trivial torsor over $S$
corresponds to the composition of morphisms
$$S\stackrel{a}{\longrightarrow}M\stackrel{\pi}{\longrightarrow} B\Gamma.$$
 
\end{rmk}

\begin{them}
\label{torsors=stack}
Let $\XX$ be a differentiable stack and $x/X$ a versal family for
$\XX$. Then 
$$\XX\cong B\Isomu(x,x)\,,$$
as groupoid fibrations over $\SS$.
\end{them}
\begin{pf}
We shall prove that the functor $f$:
\begin{align}\label{isom}
\XX&\longrightarrow B\Isomu(x,x)\\
y&\longmapsto \Isomu(x,y)\,,\nonumber
\end{align}
provides us with the required isomorphism of groupoid fibrations.

Since $x: X\to \XX$ is a representable submersion,
 it follows that $\Isomu(x,y)$
is representable.
  The fact
that $\Isomu(x,x)$ acts simply transitively on $\Isomu(x,y)$ is
clear. Thus, $\Isomu(x,y)$ is, in fact, an $\Isomu(x,x)$-torsor.

It remains to prove that (\ref{isom}) is an equivalence of
categories. Since both groupoid fibrations are stacks, we can use the
local criterion: Proposition~\ref{crit.iso},  i.e.
to prove  that $f$ is a monomorphism and an
epimorphism.

For the monomorphism property, let $y,y':S\to \XX$ be two objects of
$\XX$ lying over $S$. Let $Q$ and $Q'$ be the $\Isomu(x,x)$-torsors induced by
$y$ and $y'$ over $S$. We need to show that any isomorphism of
torsors $\phi:Q\to Q'$ comes from a 2-isomorphism $\theta:y\to
y'$. This follows from the fact that $\Isomu(y,y')$ is a sheaf: choose
a covering $\{U_i\}$ of $S$ trivializing the torsor $Q$. Then $\phi$
gives rise to isomorphisms $\theta_i:y\resto U_i\to y'\resto U_i$.
One checks that the $\theta_i$ glue together, giving rise to
$\theta$.

For the epimorphism property, suppose $Q\to S$ is an $\Isomu(x,x)$-torsor over
$S$. Then there exists a cover $\{U_i\}$ of $S$ and sections
$s_i:U_i\to Q$, trivializing $Q$ over $\{U_i\}$. The sections $s_i$
induce morphisms $x_i:U_i\to\XX$ (which are the 
compositions $U_i \to X\stackrel{x}{\to} \XX$) identifying $Q\resto U_i$ with
$x_i\upst X$. Thus we see that every $\Isomu(x,x)$-torsor over $S$ comes
locally from objects of $\XX$, proving that $f$ is an epimorphism.
\end{pf}

\begin{defn}
For a differentiable stack $\XX$, if the
 diagonal $\XX\to\XX\times\XX$ is proper,
 we call $\XX$ {\bf
  separated }or {\bf Hausdorff}. 
\end{defn}

For a differentiable stack $\XX$, 
the diagonal $\XX\to\XX\times\XX$ is always representable.
Indeed if $X_1\toto X_0$ is a Lie groupoid
representing   $\XX$, then $\XX_{\XX\times\XX} (X_0 \times X_0) \cong X_1$
and the base change map is $s\times t: X_1\to X_0\times X_0$.
Hence  $\XX$ is separated if and only if 
  $X_1\toto X_0$ is a proper groupoid.

 In the definition of Metzler~\cite{Metzler},   all differentiable stacks are
required to be separated. We believe that this
 is too restrictive.  Many interesting
differentiable stacks are not separated.

\subsection{Morita equivalence}

We have now established procedures to go back and forth between
Lie groupoids and differentiable stacks.  Given a
differentiable stack $\XX$, we choose a presentation $X_0\to\XX$ and
form the associated Lie groupoid $X_1\toto X_0$ by taking the
fibered product.  Conversely, starting with a Lie groupoid
$\Gamma\toto M$, we construct the differentiable stack $B\Gamma$ of
$\Gamma$-torsors, which comes with a canonical presentation, giving
back the groupoid $\Gamma\toto M$ we started with (up to
isomorphism).  It remains to see when exactly two different
Lie groupoids give rise to isomorphic differentiable stacks,
or put another way, what relationship there is between  various
Lie groupoids arising from various  presentations of a
differentiable stack.

\begin{defn}
Let $X\lcom$ and $Y\lcom$ be Lie groupoids. A morphism
$\phi\lcom:X\lcom\to Y\lcom$ is called a {\bf Morita morphism}, if 

(i) $\phi_0:X_0\to Y_0$ is a surjective submersion;

(ii) the diagram
$$\xymatrix{
X_1\dto\rto &X_0\times X_0\dto\\
Y_1\rto & Y_0\times Y_0}$$
is cartesian.
\end{defn}

\begin{defn}
Two Lie groupoids $X\lcom$ and $Y\lcom$ are called {\bf
Morita equivalent}, if there exists a third Lie groupoid
$Z\lcom$ and Morita morphisms $Z\lcom\to X\lcom$ and $Z\lcom \to Y\lcom$
\end{defn}

\begin{them}\label{morita-them}
Let $X\lcom$ and $Y\lcom$ be Lie groupoids. Let $\XX$ and
$\YY$ be the associated differentiable stacks, i.e., $\XX$ is the
stack of $X\lcom$-torsors and $\YY$ the stack of
$Y\lcom$-torsors. Then the following are equivalent:

(i) the differentiable stacks $\XX$ and $\YY$ are isomorphic;

(ii) the Lie groupoids $X\lcom$ and $Y\lcom$ are Morita
equivalent;
 
(iii) there exists a manifold $Q$ together with two
$C^\infty$ maps $f:Q\to X_0$ and $g:Q\to Y_0$ and (commuting) actions of $X_1$
and $Y_1$ (the 
action of $X_1$ comprising $f$ and the action of $Y_1$ comprising
$g$), in such a way that $Q$ is at the same time an $X\lcom$-torsor
over $Y_0$ (via $g$) and a $Y\lcom$-torsor over $X_0$ (via $f$). We
call such a $Q$ an {\em $X\lcom$-$Y\lcom$-bitorsor}.
\end{them}
\begin{pf}
Let us start by proving that (i) implies (iii). Choose an isomorphism
identifying $\XX$ with $\YY$. Then let $Q$ be the fibered product
$Q=Y_0\times_{\XX}X_0$. One checks that $Q$ is a bitorsor. 

To prove that (iii) implies (ii), choose a bitorsor $Q$. Let $Q_1$ be
the fibered product $Q_1=Y_1\times_{Y_0}Q\times_{X_0}X_1$. There is a
canonical way to define a Lie groupoid $Q_1\toto Q$,
together with Morita equivalences $Q\lcom\to Y\lcom$ and $Q\lcom\to
X\lcom$. 

One also proves that (ii) implies (iii). This follows from the
following two facts: (1) if $\phi: X\lcom \to Y\lcom$
is a Morita morphism, then $Q=X_0\times_{Y_{0}, s} Y_1$ is
naturally an  $X\lcom$-$Y\lcom$-bitorsor;
(2) if $Q$ is  a  $X\lcom$-$Y\lcom$-bitorsor, and 
$Q'$ is  an  $Y\lcom$-$Z\lcom$-bitorsor, then $(Q\times_{Y_0}Q')/Y_{1}$
is an $X\lcom$-$Z\lcom$-bitorsor.

Finally, we need to prove that (iii) implies (i). 
Given an $X\lcom$-torsor $F$ over $U$, let  $E=X_1\backslash (Q\times_{X_0}
F)$. Then $E$ is a $Y\lcom$-torsor over $U$, where
 the anchor map   $E\to Y_0$ is $a( [q, f])=g (q)$ and the
$Y_1$-action is $y \cdot [q, f]=[q\cdot y^{-1}, f]$.
Also it is clear that  a morphism of $X\lcom$-torsors $F_1 \to F_2$
induces a morphism of $Y\lcom$-torsor $E_1 \to E_2$ in a canonical
way.  Thus one obtains a functor, which can be easily seen to be
an equivalence of categories.
\end{pf}

\begin{rmk}
Note that (iii) is the definition of Morita equivalence  used 
in a lot of literature on operator algebras \cite{hilsum-skandalis87, muhly-renault-williams87}.
\end{rmk}

\begin{defn}
If $\XX$ is a differentiable stack and there exists a Lie groupoid
$X\lcom$ presenting $\XX$, such that $X_0$ and $X_1$ both have  constant
dimensions, then we call $\dim\XX=2\dim X_0-\dim X_1$ the {\bf
  dimension }of $\XX$.
\end{defn}

We see that,  from Theorem \ref{morita-them}, 
$\dim\XX$ is independent of the presentation of $\XX$, and
therefore is well-defined.

\begin{rmk}
Note that $\dim\XX$ can also be written as
the base dimension minus the fibre dimension
of the representing groupoid $X_1\toto X_0$, which
is also the orbit space dimension minus the
isotropy group dimension.
Also $\dim\XX$ can be negative. In particular,
if $G$ is a Lie group of dimension $n$,  the stack $BG$
 is of dimension $-n$.
\end{rmk}

\subsection{Dictionary}

Theorem~\ref{morita-them} is only the beginning of a dictionary between
differentiable stacks and Lie groupoids.  We will now list a few
propositions that give more precise information, in particular with respect
to morphisms and 2-isomorphisms.

All these results are standard in stack theory. Proofs are elementary, but
usually tedious, and we omit them. 

\subsubsection{The 2-category of Lie groupoids}

Recall the notion of {\em natural equivalence }between groupoid morphisms:

\begin{defn}
Let $\phi:X\lcom\to Y\lcom$ and $\psi:X\lcom\to Y\lcom$ be two morphisms of
Lie groupoids. A {\bf natural equivalence }from $\phi$ to $\psi$, notation
$\theta:\phi\Rightarrow\psi$, is a $C^\infty$ map $\theta:X_0\to Y_1$
such that for every $x\in X_1$ we have
$$\theta\big(s(x)\big)\ast\psi(x)=\phi(x)\ast\theta\big(t(x)\big)\,.$$
\end{defn}

Fixing the Lie groupoids $X\lcom$ and $Y\lcom$, the morphisms and natural
equivalences form a category $\Hom (X\lcom,Y\lcom)$, which is a
(set-theoretic) groupoid. With this notion of morphism groupoid, the Lie
groupoids form a 2-category. 

\subsubsection{The  Dictionary Lemmas}

We consider two Lie groupoids $X\lcom$ and $Y\lcom$ with associated
differentiable stacks $\XX$ and $\YY$, respectively.  The dictionary lemmas
relate groupoid morphisms $X\lcom\to Y\lcom$ to stack morphisms
$\XX\to\YY$. 

The first Dictionary Lemma says that a morphism of Lie groupoids induces a
morphism of associated differentiable stacks, unique up to unique
2-isomorphism: 

\begin{lem}[First Dictionary Lemma]
Let $\phi:X\lcom\to Y\lcom$ be a morphism of Lie groupoids. Let $\XX$ and
$\YY$ be differentiable stacks associated to $X\lcom$ and $Y\lcom$,
respectively. Then there exists a morphism of stacks $f:\XX\to\YY$ and a
2-isomorphism 
\begin{equation}\Label{eta-di}
\vcenter{\xymatrix{
X_0\dto\rto^{\phi_0}\drtwocell\omit{^\eta} & Y_0\dto\\
\XX\rto^f&\YY}}
\end{equation}
such that the cube
\begin{equation}\Label{cube}
\vcenter{\xymatrix@=.5pc{
& X_1\ddlto\drto \ar[rrr]^{\phi_1} &&& Y_1\ddlto\drto &\\
&& X_0\ddlto\ar[rrr]^(.3){\phi_0} &&& Y_0\ddlto \\
X_0\drto\ar[rrr]^(.3){\phi_0} &&& Y_0\drto &&\\
& \XX\ar[rrr]^f &&& \YY &}}
\end{equation}
2-commutes.  If $(f',\eta')$ is another pair satisfying these properties,
then there is a unique 2-isomorphism $\theta:f\Rightarrow f'$ such that
$\theta\ast\eta'=\eta$. \qed
\end{lem}

The second and third Dictionary Lemmas treat the converse:

\begin{lem}[Second Dictionary Lemma]\Label{sdl}
Let $f:\XX\to\YY$ be a morphism of stacks, $\phi_0:X_0\to Y_0$ a morphism of
manifolds and $\eta$ a 2-isomorphism as in~(\ref{eta-di}). 
Then there exists a unique morphism of Lie groupoids
$\phi_1:X_1\to Y_1$ covering $\phi_0$ and making the cube~(\ref{cube})
2-commutative. \qed
\end{lem}

\begin{lem}[Third Dictionary Lemma]
Let $f:\XX\to\YY$ be a morphism of stacks. Let $\phi:X\lcom\to Y\lcom$ and
$\psi:X\lcom\to Y\lcom$ be two morphisms of Lie groupoids. Let $\eta$ and
$\eta'$ be 2-isomorphisms, where $(\phi,\eta)$ and $(\psi,\eta')$ both form
2-commutative cubes such as~(\ref{eta-di}).  Then there exists a unique
natural equivalence $\theta:\phi\Rightarrow\psi$ such that the diagram
$$\xymatrix@=.5pc{
&&&& Y_1\ddlto\drto &\\
X_0\ar[urrrr]^\theta\ddrto\ar[drrr]_{\phi_0}\ar[rrrrr]^{\psi_0} &&&&&
Y_0\ddlto\\ 
&&& Y_0\drto &&\\
&\XX\ar[rrr]^f &&&\YY &}$$
2-commutes. \qed
\end{lem}

\subsection{Differentiable Spaces}


Differentiable spaces are generalizations of manifolds. The are
differentiable stacks whose isotropy groups are trivial.  They occur
when one tries to define the quotient of an equivalence relation which
is ``of Lie type'' (i.e. is given by a Lie groupoid) but the usual
quotient has bad properties (i.e. is not a manifold or not a
principal bundle quotient).  Differentiable spaces have slightly
better properties than manifolds.
  The main advantage is that
Lemma~\ref{sub.des} holds for them.

\begin{defn}
A sheaf over 
 $\SS$, which, considered as a stack over $\SS$ is
differentiable, is called a {\bf differentiable space}. 
\end{defn}

Thus a  sheaf $F$ is a differentiable space if there
exists a manifold $X$ and a surjective representable submersion $X\to
F$. 


\begin{numex}
If a Lie group acts on a manifold freely, but not properly, we get a
differentiable space.
\end{numex}

\begin{prop}
The differentiable stack $\XX$ defined by a Lie groupoid $X\lcom$ is
(isomorphic to) a differentiable space if and only if $X\lcom$ is a Lie
equivalence relation (i.e. $X_1\to X_0\times X_0$ is injective). 

In particular, if $X\lcom$ and $Y\lcom$ are Morita
equivalent Lie groupoids, then $X\lcom$ is an equivalence relation if
and only if $Y\lcom$ is. 
\end{prop}

Thus we may think of differentiable spaces as Lie equivalence
relations up to Morita equivalence.

\begin{lem}[Submersive descent for differentiable
spaces]\Label{sub.des}
Let $F$ be a sheaf over $\SS$. Let $X$ be a manifold and $F\to X$
a morphism. Suppose that $U\to X$ is a surjective submersion of
manifolds
and that the sheaf $G=U\times_X F$ is a differentiable space. Then $F$
is a differentiable space.\qed
\end{lem}

Note that there is no corresponding statement for manifolds.  For
manifolds we only have \'etale descent (Lemma~\ref{descent}). 

\begin{defn}
Let $\XX$ and $\YY$ be stacks over $\SS$. We call a morphism
$f:\XX\to\YY$ {\bf weakly representable}, if for every 
representable submersion $U\to\YY$, where $U$ is a manifold, the
fibered product 
$V=\XX\times_{\YY}U$ is isomorphic to a differentiable space.
\end{defn}

\begin{prop}
Let $\XX$ and $\YY$ be differentiable stacks. The morphism
$f:\XX\to\YY$ is weakly representable if there exists a presentation
$Y\to\YY$ such that $X=\XX\times_{\YY}Y$ is isomorphic to a
differentiable space.
\end{prop}
\begin{pf}
The proof is very similar to the proof of Lemma~\ref{critrep}. 
We need to use submersive descent for differentiable spaces.
\end{pf}

\begin{numex}
Representable morphisms are weakly representable.  In particular, the
diagonal $\XX\to\XX\times\XX$ of a differentiable stack is 
weakly representable, and any $\cinf$-map of manifolds is weakly
representable. 

Moreover, any morphism
from a differentiable space to a differentiable stack is weakly
representable, and any morphism of differentiable stacks which is faithful
is weakly representable. 
\end{numex}

\begin{rmk}
We get a weaker notion of differentiable stack if we work with
groupoids where $X_0$ and $X_1$ are differentiable spaces rather than
manifolds.  Equivalently, we can relax the condition that the diagonal
$\XX\to\XX\times \XX$ be representable to it being weakly representable.
We could call these stacks {\em weakly differentiable stacks}. 

For example, the quotient $\rr/\qq$ is a differentiable space but not
a manifold.  It is also a group.  The associated stack $B(\rr/\qq)$ is
weakly differentiable but not differentiable.
\end{rmk}


\begin{rmk}
It would be interesting to investigate the relationship
between differentiable spaces and Souriau's
 diffeology structures \cite{Souriau}.
\end{rmk}

\section{Homology  and cohomology}

Here our goal is to define the cohomology of a differentiable stack
with values in a sheaf (or a complex of sheaves) of abelian groups.
Of particular interest is the {\em de Rham complex}, which gives rise
to de Rham cohomology.  


Then we pass to Lie groupoids and define 
 the cohomology of a Lie
groupoid with values in a sheaf of abelian groups. This cohomology is
Morita invariant.
 For any complex of sheaves of abelian groups, we also define
  a double complex and
its associated cohomology groups. These cohomology groups are not
necessarily Morita invariant, but they will be if all component
sheaves of the complex are {\em acyclic on manifolds}. An example of
this is the de Rham complex.  Thus de Rham cohomology of a groupoid is
also Morita invariant.


\subsection{Sheaves over stacks and their cohomology}

Let $\XX$ be a differentiable stack. We endow the category
$\XX$ with a Grothendieck topology defined as follows: call a family
$\{x_i\to x\}$ of morphisms in $\XX$ a {\em  covering family }of the
object $x\in\XX$, if the image family $\{U_i\to U\}$ in $\SS$ is a
covering family, i.e. is a family of \'etale maps such that $\coprod
U_i\to U$ is surjective.  One checks that, indeed, the axioms of a
topology are satisfied. 
Thus we may now speak of sheaves over $\XX$: i.e.
contravariant functors
$\XX\to(\text{sets})$ satisfying the sheaf axioms.  We get the
category $(\text{sheaves/$\XX$})$ of sheaves over $\XX$.

\begin{numrmk}
Let $F$ be a sheaf over the stack $\XX$. Consider  $F$ as a
category fibered in groupoids $F\to\XX$.  Then, by composing with
$\XX\to\SS$, we may turn $F$ into a category fibered in groupoids over
$\SS$. One checks that $F$ is then a stack over $\SS$ and that
$F\to\XX$ is a morphism of stacks.  Moreover, $F\to\XX$ is faithful.

Conversely, if $f:\YY\to\XX$ is a faithful morphism of stacks over
$\SS$, we may associate a sheaf $F$ over $\XX$ defined by
$F(x)=\{(y,\phi)\st \phi:x\to f(y)\}/\sim$, where
$(y,\phi)\sim(y',\phi')$ if there exists $\eta:y\to y'$ such that
$f(\eta)\comp\phi=\phi'$.  We call $F$ the {\em sheaf of sections }of
$f:\YY \to \XX$.

Thus we get an equivalence of categories between stacks which are
faithful over $\XX$ and sheaves over $\XX$.
\end{numrmk}

We define the global section functor
$$\Gamma(\XX,\argument):(\text{sheaves/$\XX$})
\longrightarrow(\text{sets})$$ 
by $\Gamma(\XX,F)=\Hom_\XX(\XX,F)$, the set of morphisms of sheaves
over $\XX$ from the trivial sheaf $\XX$ (whose set of sections is
always the one point set $\{\ast\}$) to the sheaf $F$. 

\begin{numrmk}
If $\XX$ is differentiable, $X\to\XX$ a presentation and $X_1\toto X$
the associated Lie groupoid, then for any sheaf $F$ on $\XX$ we have a
short exact sequence of sets
$$\xymatrix{
\Gamma(\XX,F)\ar[r]& F(X) \ar@<-.5ex>[r]\ar@<.5ex>[r] &
F(X_1)}\,.$$
In other words $\Gamma(\XX,F)$ is the equalizer of the two restriction
maps $F(X)\to F(X_1)$. (Note that there are two canonical morphisms
$X_1\to\XX$, so that $F(X_1)$ is ambiguous notation.  But since both
morphisms $X_1\to\XX$ are canonically isomorphic, it is irrelevant
which choice one makes for $F(X_1)$.)
\end{numrmk}

Restricting $\Gamma(\XX,\argument)$ to the category of sheaves of
abelian groups over $\XX$ we get the functor
$$\Gamma(\XX,\argument):(\text{abelian
sheaves/$\XX$})\longrightarrow(\text{abelian groups})\,.$$
This functor is left exact, and $(\text{abelian sheaves/$\XX$})$ has
sufficiently many injectives, so we may derive this functor to get the
functors 
$$H^i(\XX,\argument):(\text{abelian
sheaves/$\XX$})\longrightarrow(\text{abelian groups})\,.$$

Passing to the derived category of complexes of abelian sheaves over
$\XX$, we get the total derived functor
$$R\Gamma(\XX,\argument):D^+(\XX)\longrightarrow D^+(\text{abelian
groups})\,.$$
For a complex $M\com\in D^+(\XX)$ of abelian sheaves on $\XX$, the
homology groups of the complex $R\Gamma(\XX,M\com)$ are denoted by 
$$\hh^i(\XX,M\com)=h^i\big(R\Gamma(\XX,M\com)\big)$$
and called the {\em hypercohomology }groups of $\XX$ with values in
$M\com$. 

Of course, if $M\com\to N\com$ is a quasi-isomorphism of complexes of
abelian sheaves over $\XX$, we get induced isomorphisms 
$\hh^i(\XX,M\com)\to\hh^i(\XX,N\com)$ on hypercohomology.

\begin{defn}
Let $U$ be a manifold.  A sheaf in the usual sense (defined only on
open subsets of $U$) is called a {\bf small }sheaf on $U$.  This is to
distinguish such sheaves from sheaves over the stack over $\SS$
obtained from $U$.
\end{defn}

Let $\XX$ be a stack over $\SS$ and $F$ a sheaf over $\XX$. Let $x$ be
an object of $\XX$ lying over the manifold $U\in\SS$.  The small sheaf
on $U$, which maps the open subset $V\subset U$ to $F(x\resto V)$ is
called the small sheaf {\em induced }by $F$ via $x:U\to \XX$ on
$U$. Notation: $F_{x,U}$, or simply $F_U$, if there is no risk of
confusion.
Given a morphism $\theta:y\to x$ in $\XX$ lying over $f:V\to U$ in
$\SS$, there is an induced morphism of small sheaves over $V$ called
$\theta\upst:f^{-1}F_{x,U}\longrightarrow F_{y,V}$.
This induced morphism is contravariantly functorial in $\theta$.

\begin{lem}\label{additional}
If $\XX$ is representable, represented by the manifold $X$, then for any big sheaf $F$ over $\XX$, we have
$H^i(\XX,F)=H^i(X,F_X)$, for all $i$.
\end{lem}
\begin{pf}
This follows from the fact that for a manifold $\XX=X$ the functor
$F\mapsto F_X$, which maps a big sheaf to its associated small sheaf is
exact. This property fails for stacks.
\end{pf}

\subsection{Differential forms}

Let $\XX$ be a differentiable stack.
Define the sheaf $\Omega^i_\XX$ of differential  forms of degree $i$
on $\XX$ as follows: for an object $x\in\XX$ lying over $U\in\SS$, we
let $\Omega^i_\XX(x)=\Omega^i(U)$ be the $\rr$-vector space of
($\rr$-valued) differentiable $i$-forms on $U$.  For a morphism $y\to
x$ in $\XX$ lying over the $\cinf$-map $V\to U$, we define the
restriction map $\Omega^i_\XX(x)\to\Omega^i_\XX(y)$ to be the
pullback map $\Omega^i(U)\to\Omega^i(V)$. The sheaf axioms are easily
verified. Note that the $\Omega^i_\XX$ are sheaves of $\rr$-vector
spaces, i.e. they take values in $(\text{$\rr$-vector
spaces})\subset(\text{sets})$. 

The sheaf $\Omega^0_\XX$ is also called the {\em structure sheaf }of
$\XX$, notation $\O_\XX$.
It is isomorphic to the sheaf of sections of the projection
$\XX\times\rr\to\XX$. 

Note that none of the
$\Omega^i_\XX$, for $i>0$, are isomorphic to sheaves of sections of any
morphism of differentiable stacks $\YY\to\XX$.

The exterior derivative $d:\Omega^i(U)\to\Omega^{i+1}(U)$, where $U$
is any manifold, commutes with the pullback of forms via any
$\cinf$-map $V\to U$.  Thus $d$ induces a homomorphism of sheaves
$d:\Omega^i_\XX\to \Omega^{i+1}_\XX$, for all $i\geq0$.  Clearly,
$d^2=0$, and so we have defined a complex $\Omega\com_\XX$ of sheaves
of $\rr$-vector spaces over $\XX$. We call $\Omega\com_\XX$ the {\em
de Rham }complex of $\XX$.  Its hypercohomology is called the
{\em de Rham }cohomology of $\XX$:
$$H^i_{DR}(\XX)=\hh^i(\XX,\Omega\com_\XX)\,.$$
If there is any danger of confusion (for example if $\XX$ is a
manifold), we refer to $\Omega_\XX\com$ as the {\em big }de Rham
complex of $\XX$.


Let $\rr_\XX$ denote the sheaf over $\XX$
 defined by
$$\rr_\XX(x)=\{f:U\to\rr\st\text{$f$ is locally constant}\}\,,$$
for any object $x$ of $\XX$ lying over $U\in\SS$. The sheaf $\rr_\XX$ is
in a natural way a subsheaf of the structure sheaf $\Omega^0_\XX$. If
we let $\tilde{\rr}$ denote the manifold with the same underlying set as
$\rr$, but the discrete differentiable structure, then we may identify
$\rr_\XX$ with the sheaf of sections of the projection
$\XX\times\tilde{\rr}\to\XX$. 

The usual Poincar\'e Lemma proves that the big de Rham complex
$\Omega\com_\XX$ is a resolution of $\rr_\XX$. Thus we conclude that 
$$H^i_{DR}(\XX)=H^i(\XX,\rr_\XX)\,,$$
for all $i$.

\subsection{Groupoid cohomology}

Let $\XX$ be a differentiable stack over $\SS$,
and $X\to \XX$ an atlas for $\XX$.
  Define for all $p\geq0$
$$X_p=\underbrace{X\times_\XX\ldots\times_\XX X}_{\text{$p+1$
times}}\,.$$ 
(hence $X_0=X$). 
Since $X\to \XX$ is a representable submersion, all $X_p$ are
manifolds.  Of course, $X_1\toto X_0$ is the Lie groupoid associated
to the atlas $X\to\XX$. Furthermore, we assume
that $X_1\toto X_0$  is  a Hausdorff and second
countable Lie groupoid \cite{Mackenzie, MoerdijkM}.

We have $p+1$ projection maps $X_p\to X_{p-1}$, giving rise to a
diagram
\begin{equation}\label{sim.ma}
\xymatrix{
\ldots\ar@<-1.5ex>[r]\ar@<-.5ex>[r]\ar@<1.5ex>[r]\ar@<.5ex>[r] & X_2
\ar[r]\ar@<1ex>[r]\ar@<-1ex>[r] & X_1\ar@<-.5ex>[r]\ar@<.5ex>[r]
&X_0\,.}
\end{equation}
In fact more is true: $X\lcom$ is a simplicial manifold. 

Every $X_p$
has $p+1$ canonical projections $X_p\to\XX$.  They are all canonically
isomorphic to each other.  Choose any one of them and call it
$\pi_p:X_p\to\XX$.  As usual, we identify $\pi_p$ with an object of
$\XX$ lying over $X_p$.
Let $F$ be a sheaf over $\XX$.
 We denote the set $F(\pi_p)$ by $F(X_p)$. 
Let $F_p$ denote the small sheaf on $X_p$ induced by $F$. Then we have
$F(X_p)=\Gamma(X_p,F_p)$.

Diagram~(\ref{sim.ma}) induces a diagram
\begin{equation}\label{cosim}
\xymatrix{
F(X_0)\ar@<-.5ex>[r]\ar@<.5ex>[r] &
F(X_1)\ar[r]\ar@<1ex>[r]\ar@<-1ex>[r] &
F(X_2)\ar@<-1.5ex>[r]\ar@<-.5ex>[r]\ar@<1.5ex>[r]\ar@<.5ex>[r]
&\ldots}
\end{equation}
which can, in fact, be refined to a cosimplicial set.

Now assume that $F$ is a sheaf of abelian groups.  Let
$\del:F(X_p)\to F(X_{p+1})$ be the alternating sum of the maps of
Diagram~(\ref{cosim}).
We obtain a complex of abelian groups
$$\xymatrix{
F(X_0)\rto^\del&F(X_1)\rto^\del&F(X_2)\rto^\del & \ldots}$$
The homology groups of this complex are denoted by
$${\check H}^i(X\lcom,F)=h^i\big(F(X\lcom)\big)$$
and called the {\em \v{C}ech }cohomology groups of $F$ with respect to
the covering $X\to \XX$.

Note that when $F$ is the sheaf $\Omega^0$,
${\check H}^i(X\lcom, \Omega^0 )$ is also called 
{\em groupoid cohomology } with trivial coefficients \cite{Mackenzie, WeinsteinXu}.

Lemma~\ref{additional} gives us the following lemma and proposition.

\begin{lem}
There is an $E_1$ spectral sequence
$$H^q(X_p,F_p)\Longrightarrow H^{p+q}(\XX,F)\,.$$
\end{lem}

\begin{prop}
Assume that for every $p$ the induced small sheaf $F_p$ is acyclic,
i.e. satisfies $H^i(X_p,F_p)=0$, for all $i>0$.
Then we have
$${\check H}^i(X\lcom,F)=H^i(\XX,F)\,.$$
\end{prop}

\begin{cor}\label{wfo}
We have, for all $i,j\geq0$,
$${\check H}^j(X\lcom,\Omega^i_\XX)=H^j(\XX,\Omega^i)\,.$$
\end{cor}

In particular we see that the sheaf cohomology $H^j(\XX,\Omega^0 )$
is isomorphic to the groupoid cohomology.



In the sequel, when $X\lcom$ is a Lie groupoid, a {\em sheaf over
}$X\lcom$ is defined to be a sheaf over the associated stack
$\XX$.  Moreover, we define groupoid cohomology 
$H^i(X\lcom, F)$ to be equal to $H^i(\XX,F)$. This is in line with
sheaf cohomology of simplicial manifolds.

Now let $M\com$ be a complex of abelian sheaves over $\XX$, bounded
below.  Denote the differential on $M\com$ by $d$. Let $X\to\XX$ be an
atlas.  For every $i$ we get a \v{C}ech complex $M^i(X\lcom)$, with
differential $\del$. 
Because $d$ and $\del$ commute, we obtain, in fact, a double complex
$$\big(\{M^i(X_p)\}_{i,p},d,\del\big)\,.$$
Our convention will always be that $\del$ is the vertical
differential, $d$ the horizontal differential.
The homology groups of the associated total complex are denoted by
$$\check{\hh}^i(X\lcom,M\com)=h^i\big(\tot M\com(X\lcom)\big)$$
and called the \v{C}ech hypercohomology groups of $M\com$ with respect
to the covering $X\to \XX$.

\begin{prop}
Assume that for every $i$ and every $p$ the small sheaf $M^i_p$
induced by $M^i$ on $X_p$ is acyclic.  Then we have
$$\check{\hh}^i(X\lcom,M\com)=\hh^i(\XX,M\com)\,,$$
for all $i\geq0$.
\end{prop}

\begin{cor}
We have, for every atlas $X\to\XX$,
$$H_{DR}^i(\XX)=\check{\hh}^i(X\lcom,\Omega\com_\XX)=
h^i\big(\tot \Omega\com(X\lcom)\big)\,.$$
\end{cor}

In the sequel, we also use $H^i_{DR}(X\lcom)$ to denote the above
cohomology group. We also write $\Omega^k
(X\lcom):=\oplus_{i+j=k}\Omega^i (X_j)$, $Z^k (X\lcom)$ and $B^k
(X\lcom)$ to denote the spaces of $k$-cochains, $k$-cocycles and
$k$-coboundaries of $\tot \Omega\com(X\lcom)$, respectively.

In particular, if $\XX$ is 
representable, represented by a manifold $X$, then
$H_{DR}^i(\XX)$ coincides
with  the $i$-th homology group of the usual de Rham complex of $X$. Thus
our definition of de Rham cohomology gives the usual de Rham
cohomology in the case of manifolds.

\begin{rmk}
The wedge product of differential forms turns $\Omega\com_\XX$ into
a sheaf of differential graded $\rr$-algebras.  From general
principles
it follows that $R\Gamma(\XX,\Omega\com)$ is also a differential
graded $\rr$-algebra. Thus the hypercohomology
$\hh^\ast(\XX,\Omega\com)$ is a graded $\rr$-algebra.


This multiplicative structure can be described 
explicitly  at the level of the
double complex $\Omega\com(X\lcom)$ associated to an atlas $X\to\XX$
of $\XX$.

Let $a \in \Omega^{k}(X_{p})$ and $b\in \Omega^{l}(X_{q})$. Define
$a\cup b\in \Omega^{k+l}(X_{p+q})$ by

\begin{equation}
\label{eq:cup1}
a\cup b = (-1)^{kq} p_{1}^*a\wedge p_{2}^* b,
\end{equation}
where $p_1, \ p_2$ are the natural  projections from $X_{p+q}$ to $X_{p}$ and
$X_q$  given by $(x_1, \cdots , x_{p}, x_{p+1}, \cdots , x_{p+q})\to
(x_1, \cdots , x_{p})$
and $(x_1, \cdots , x_{p}, x_{p+1}, \cdots , x_{p+q})\to
(x_{p+1}, \cdots , x_{p+q})$ respectively. One checks that
for  any  $a, b, c\in  \Omega\com(X\lcom)$  the following identities hold:
\begin{eqnarray}
&&  (a\cup b)\cup c=a\cup (b\cup c)\\
&&\delta (a\cup b)=\delta a \cup b+ (-1)^{|a|}a\cup \delta b,
\end{eqnarray}
where   $\delta$ 
is the total differential of the double complex $\Omega\com(X\lcom)$,
and  $|a|$ denotes the total degree of $a$ in the double complex
$\Omega\com(X\lcom)$.
Thus   $(\Omega\com(X\lcom), \cup , \delta )$ is a graded
differential algebra. One can prove that  the
induced graded algebra structure on its  cohomology groups
coincides with the one on $\hh^\ast(\XX,\Omega\com)$.
\end{rmk}

\section{$S^1$-bundles and $S^1$-gerbes}

In this section we study connections on bundles and gerbes. We
often restrict to the case of $S^1$ as structure group.

\subsection{$S^1$-bundles}
In this subsection, we study  differential geometry,
including characteristic classes, of $S^1$-bundles over a 
differentiable stack in terms of Lie groupoids.

Let $\XX$ be a differentiable stack and $X_1\toto X_0$ a Lie groupoid
presenting $\XX$. By an $X\lcom$-space, we mean 
a manifold $P_0\to X_0$ with a left $X\lcom$-action.

\begin{defn}
An $S^1$-bundle over $\XX$ is a 2-commutative diagram
$$\xymatrix{
\PP\times S^1\rto^\sigma\dto_{\text{\rm proj.}}\drtwocell\omit &
\PP\dto^\pi\\
\PP\rto^\pi & \XX}$$
such that the pullback via $U\to \XX$, for any submersion from a
manifold $U$, defines an $S^1$-bundle over $U$.
(Note 
that this implies that $\PP\to\XX$ is a representable surjective
submersion, and hence that $\PP\resto
U$ is an $S^1$-bundle for {\em every } $U\to\XX$, submersive or not.)
\end{defn}

\begin{defn}
Let $X_1\toto X_0$ be a Lie groupoid. A (right) {\em $S^1$-bundle }over
$X_1\toto X_0$ is a (right) $S^1$-bundle $P_0$ over $X_0$,
 together with a (left) action of $X\lcom$ on $P_0$,
 which respects the $S^1$-action, i.e. we have
\begin{equation}
\label{eq:bundle}
(\gamma \cdot x) \cdot t=\gamma \cdot (x \cdot t ) ,
\end{equation}
for all $t\in S^1$ and  all compatible pairs $(\gamma,x)\in
\Gamma\times_{t,X_0} P_0$.
\end{defn}

\begin{prop}
 There is a canonical equivalence of categories
$$(\text{\rm $S^1$-bundles over $\XX$})\longrightarrow(\text{\rm
  $S^1$-bundles over $X_1\toto X_0$})\,.$$
\end{prop}
\begin{pf}
Let $\PP\to\XX$ be an $S^1$-bundle.
Denote the pullback of $X_0$ via $\PP\to\XX$ by $P_0$.
Thus $P_0\to X_0$ is an $S^1$-bundle by assumption.
And $P_0\to \PP$ is a representable
submersion.  Let $P_1\toto P_0$ be the associated groupoid.
 We get an induced morphism of groupoids $P\lcom\to X\lcom$,
 which is {\em cartesian}, i.e. the diagram
$$\xymatrix{P_1\rto\dto & X_1\dto \\ P_0\rto & X_0}$$
is a pullback diagram of manifolds, where the vertical maps are source
maps (or, equivalently, target maps). 
Therefore $P_1\to X_1$ is an $S^1$-bundle and the
vertical maps in the diagram above are
$S^1$-bundle maps. As a consequence,  $X_1$ acts on $P_0$ and
Eq. \eqref{eq:bundle} is satisfied. The
functor in the proposition is  $\PP\mapsto P_0$.

Conversely, given an $S^1$-bundle $P_0$ over
$X_1\toto X_0$, 
let $P_1=X_1\times_{t,X_0} P_0$.
  Action and projection form a diagram
$P_1\toto P_0$, and it is easy to check that
$P_1\toto P_0$ is  naturally  a groupoid (called the transformation
groupoid of the $X_1$-action).  
 It is clear that $P_1$ is an $S^1$-bundle over $X_1$.
Moreover, there is a natural morphism of groupoids
 $\pi$ from $P_1\toto P_0$ to $X_1\toto X_0$, which respects
the $S^1$-bundle structures $P_1\to X_1$ and $P_0\to X_0$. 
 Let $\PP$ be the corresponding
stack of $P\lcom$-torsors. The groupoid morphism $P\lcom\to
X\lcom$ induces a morphism of stacks $\PP\to\XX$, which is
representable, as its pullback to $X_0$ equals $P_0\to X_0$.
It is also simple to see that for
any morphism  $U\to \XX$,
$\PP\resto U$ is an $S^1$-bundle. Therefore
$\PP$ is an $S^1$-bundle over $\XX$.
 The backwards functor is given by $P_0\mapsto\PP$.
\end{pf}

As a consequence,
$S^1$-bundles over a given Lie groupoid
$X_1\toto X_0$ are classified by     
$H^1(X\lcom,S^1)$.  The exponential
sequence $\zz\to \Omega^0\to S^1$ induces a boundary map
$H^1(X\lcom,S^1)\to
H^2(X\lcom,\zz)$; the image of the class of 
  an $S^1$-bundle under this
boundary map is called its {\em Chern class}.

Let $P_0 \to X_0$ be a principal $S^1$-bundle over $X_1\toto X_0$. Let
$\theta\in\Omega^1 (P_0 )$ be a connection 1-form 
on $P_0$. One checks that $\delta  \theta\in \Omega^2_{DR}(P\lcom)$
descends to $ \Omega^2_{DR}(X\lcom)$. In other words, there exist
unique $\omega\in \Omega^1(X_1)$ and $\Omega\in
\Omega^2(X_0)$ such that $\pi\upst(\omega+\Omega)=\delta\theta$. 

\begin{prop}
\label{pro:chern}
The class $[\omega+\Omega]\in H^2_{DR}(X\lcom)$ is independent of the
choice of the connection $\theta$ on $P_0\to X_0$. Under the canonical
homomorphism $H^2(X\lcom,\zz)\to  H^2(X\lcom,\rr )\cong H^2_{DR}(X\lcom)$,
 the Chern class
of $P$ maps to $[\omega+\Omega]$.
\end{prop}
\begin{pf}
The proof of independence of choice of connection is a direct
calculation. See \cite{LTX}.
  Thus we concentrate on the second statement.

Let $\XX$ be the differentiable stack represented by $X\lcom$ and $\PP$
the $S^1$-bundle on $\XX$ defined by $P\lcom$. 

Consider on $\XX$ the diagram of abelian sheaves
\begin{equation}\Label{mor.com}
\vcenter{\xymatrix{
0\rto & \zz\dto\rto & \Omega^0\dto\rto^{\exp} & S^1\dto^{d\log}\rto
&0\dto\rto & \ldots\\ 
0\rto & \rr\rto & \Omega^0\rto^d & \Omega^1\rto^d & \Omega^2\rto &
\ldots}}
\end{equation}
The upper row is a resolution $[\Omega^0\to S^1]$ of $\zz$, and the
lower row is the de Rham resolution $\Omega\com$ of $\rr$, and the
whole diagram is a morphism of resolutions.

It follows that we have a commutative diagram
\begin{equation}\Label{mor.coh}
\vcenter{\xymatrix{
{H^i(\XX,\zz)}\rrto^\sim\dto && {\hh^i(\XX,[\Omega^0\to S^1])}\dto\\
{H^i(\XX,\rr)}\rrto^\sim && {H^i_{DR}(\XX)}}}
\end{equation}
This diagram gives us a way to calculate the image of the Chern class
of $\PP$ in de Rham cohomology. 

In fact, consider now the \v Cech resolution (simplicial manifold)
$X\lcom$. For any complex of sheaves $F\com$, we get an associated
double complex of abelian groups $F\com(X\lcom)$ and canonical maps
$$h^i\big(\tot F\com(X\lcom)\big)\tto \hh^i(\XX,F\com)\,.$$
(If all $F^q$ are acyclic on every $X_p$, then these are isomorphisms.)

If we assume that $P_0$ admits a section $\sigma$ over $X_0$,
then $\rho=t\upst(\sigma)-s\upst(\sigma)\in S^1(X_1)$ is a 2-cocycle
in $\tot [\Omega^0\to S^1](X\lcom)$ and the associated cohomology
class $[\rho]\in\hh^2(\XX,[\Omega\to S^1])$ is the image of the Chern
class of $\PP$ under the upper row of~(\ref{mor.coh}).

The morphism of complexes of sheaves~(\ref{mor.com}) induces a
morphism of double complexes of abelian groups
$$[\Omega^0\to S^1](X\lcom)\tto \Omega\com(X\com)\,.$$
The image of $\rho$ under this morphism is 
$$d\log\big(t\upst(\sigma)-s\upst(\sigma)\big)\in\Omega^1(X_1)\,.$$

Now, $\sigma$ also induces a flat connection
on the $S^1$-bundle $P_0\to X_0$ (ignoring the groupoid action),
hence a connection $1$-form $\theta\in \Omega^1(P_0)$. We
have $\delta(\theta)=s\upst(\theta)-t\upst(\theta)$, and thus, all we
need to prove is that 
$$s\upst(\theta)-t\upst(\theta) = \pi\upst
\big(d\log(t\upst\sigma-s\upst\sigma)\big) .$$
This can  be easily checked.

Let us now prove the general case, i.e. the case where $P_0$ is not
necessarily trivial over $X_0$.

Choose an open covering $(U_i)$ of $X_0$ such that each $U_i$ is
contractible. Let $Y_0=\coprod U_i$ and
$Y_1=\coprod U_i \times_{s} X_1 \times_{t} U_j$.
Then $Y_1\toto Y_0$ is a Lie groupoid Morita equivalent
to $X_1\toto X_0$. In fact, the projection 
 $f:Y\lcom\to X\lcom$ is a Morita morphism.
 Let $Q\lcom$ be the
$S^1$-bundle over $Y\lcom$ induced from $P\lcom$.
Since $Q_0$ can be trivialized over $Y_0$, we
  choose $[\omega'+\Omega']$
coming from a trivialization of $Q_0$, as above (in which case
$\Omega'=0$). 

Choose $(\omega,\Omega)$ for $P\lcom$ as in the statement of the
proposition. By Morita invariance of de Rham cohomology of Lie
groupoids, to prove that $[\omega+\Omega]$ is the Chern class of
$\PP$, it suffices to do this for $f\upst[\omega+\Omega]$. Thus we
reduce to proving that $f\upst[\omega+\Omega]=[\omega'+\Omega']\in
H^2_{DR}(Y\lcom)$, which is just the invariance under choice of
connection.
\end{pf}

Note that the class $[\omega+\Omega]$ in the above proposition
 is an integer class in $H^2_{DR}(X\lcom)$, and $\omega+\Omega$
an integer 2-cocycle in $Z^2_{DR}(X\lcom)$. In general
a $k$-cocycle in $Z^k_{DR}(X\lcom)$ is said to be
an {\em integer $k$-cocycle} if it defines an {\em integer class}
 in $H^k_{DR}(X\lcom)$,
i.e. a  class  in the image of the homomorphism
 $H^k(X\lcom,\zz)\to  H^k(X\lcom,\rr )\cong H^k_{DR}(X\lcom)$.

\begin{defn}
A connection 1-form $\theta$ on $P_0$ is called a
 {\em pseudo-connection} on
$P\lcom$. The de~Rham cocycle $\omega+\Omega\in Z^2_{DR}(X\lcom)$ such
that $\pi\upst(\omega+\Omega)=\delta\theta$ is called the
 {\em pseudo-curvature }of $\theta$.

A pseudo-connection $\theta$ is a {\em connection }if
$\del\theta=0$.

A {\em flat connection }is a pseudo-connection $\theta$ 
whose pseudo-curvature vanishes.       
\end{defn}

\begin{rmk}
Unlike in the manifold case, connections do not always exist.
Thus connections are not necessarily as useful
to compute characteristic classes as in the manifold case. 
For instance, the universal $S^1$-bundle $\ast\to BS^1$, which
corresponds to $S^1\to \ast$ considered as an
$S^1$-bundle over the groupoid $S^1\toto \ast$ (where the groupoid
$S^1\toto \ast$ acts on $S^1$ by left translation), does not admit any
connections.  (Any connection on the universal bundle would
necessarily be flat, and the existence of a flat connection on the
universal bundle would imply that all connections on all
bundles over all manifolds were flat.)
\end{rmk}

A {\em flat} $S^1$-bundle is an $S^1$-bundle
with a flat connection.  It is simple to see that
a flat $S^1$-bundle over $X\lcom$ is equivalent to
a $\rr/\zz$-bundle over $X\lcom$. Therefore, the
equivalence classes of  flat $S^1$-bundles
are classified by $H^1(X\lcom , \rr/\zz )$. The functor
$$\{ \text{flat $S^1$-bundles over $X\lcom$}\}
\longrightarrow H^1(X\lcom , \rr/\zz )$$
is called the {\em holonomy} map. When $X\lcom$ is a
manifold, this reduces to the usual holonomy map
for flat bundles.

We are now ready to prove the following proposition, which generalizes
the prequantization theorem of Kostant and Weil \cite{Kostant, Weil}.

\begin{prop}
\label{prop:prequantization}
Assume that $\check 	H^1 (X\lcom, \Omega^0 )=0$. 
Let $\omega+\Omega\in Z^2_{DR}(X\lcom)$ be an integer 2-cocycle. Then
there exists an $S^1$-bundle $P\lcom$ over $X_1\toto X_0$ and a
pseudo-connection $\theta$ whose pseudo-curvature is
$\omega+\Omega$.  

Moreover, the set of isomorphism classes of all such pairs
$(P\lcom,\theta)$ is a simply transitive $H^1(X\lcom,\rr/\zz)$-set.
 Here
$(P\lcom,\theta)$ and $(P'\lcom, \theta')$ are isomorphic 
if $P_1$ and $P'_1$ are
isomorphic as $S^1$-bundles over $X_1\toto X_0$ and under such an
isomorphism $\theta$ is identified with $\theta'$.       
\end{prop}
\begin{pf}
Consider the exact sequence
$$\to H^1(X\lcom,S^1)\stackrel{\phi}{\to}
  H^2(X\lcom,\zz)\to 
H^2(X\lcom,\Omega^0)\to$$
induced by the  exponential sequence $\zz\to \Omega^0\to S^1$.
The map $H^2(X\lcom,\zz)\to 
H^2(X\lcom,\Omega^0)$ factors through $H^2 (X\lcom , \rr  )
\cong H^2_{DR}(X\lcom )$,
i.e. we have the following commutative diagram
$$\xymatrix@C=1pc{H^2 (X\lcom , \zz )\rrto\drto && H^2_{DR}(X\lcom
  )\dlto^p\\ &  H^2 (X\lcom,\Omega^0) \cong {\check H}^2(X\lcom,
  \Omega^0 )  &}$$   
where  $p$ is the natural projection. It is clear that
$p([\omega +\Omega])=0$.  Thus
there is an  $S^1$-bundle  $P\lcom $
over $X_1\toto X_0$, whose Chern class equals  $[\omega+\Omega]$.     
Let $\theta'\in \Omega^1 (P_0 )$ be a pseudo-connection, and
$\delta\theta'= \pi\upst(\omega'+\Omega')$.
According to Proposition \ref{pro:chern}, $\omega+\Omega$ and
$\omega'+\Omega'$  are cohomologous. Hence, $\omega+\Omega-
(\omega'+\Omega')=\delta (f+\alpha)$ for some $f\in \Omega^0 (X_1)$
and $\alpha \in \Omega^1 (X_0)$. It thus follows that
$\partial f=0$, which implies that $f=\partial g$ for
$g\in  \Omega^0 (X_0)$ since $\check H^1 (X\lcom , \Omega^0 )=0$. Thus
$\delta f=\delta \partial g=\delta dg $. Let
$\theta =\theta'+\pi^*(\alpha +dg)\in \Omega^1 (P_0 )$.
It is clear that $\theta$ is the desired pseudo-connection.

If $(P, \theta )$ and $(P', \theta')$ are two such  $S^1$-bundles,
then $(P\otimes \overline{P'}, \pr_1^* \theta +\pr^*_2  \overline{\theta'})$
is a flat bundle, whose isomorphism class  is  
classified by   $H^1(X\lcom,\rr/\zz)$. Here
$P\otimes \overline{P'}$ denotes the $S^1$-bundle $(P\times_{X_0} P')/S^1$,
and $\pr_1: P\otimes \overline{P'}\to P$ and $\pr_2: P\otimes \overline{P'}
\to P'$ are projections.
This completes the proof.
\end{pf}

\begin{rmk}
The condition $\check     H^1 (X\lcom, \Omega^0 )=0$
always holds for a proper Lie groupoid $X_1\toto X_0$
according to Crainic \cite{Crainic}. In particular,
when $X_1\toto X_0$ is a manifold $M\toto M$ (which is clearly a 
proper Lie groupoid), an integer 2-cocycle in
$Z^2_{DR} (X\lcom)$ corresponds to an  integer
closed two-form  on $M$. Thus  Proposition
\ref{prop:prequantization} reduces to the usual
prequantization theorem of
Kostant and Weil \cite{Kostant, Weil}.
\end{rmk}

\subsection{$S^1$-gerbes and $S^1$-central extensions}
Let us first  recall the definition of gerbes.
Let $\XX$ be the differentiable stack associated to the Lie groupoid
$X_1\toto X_0$. Thus $\XX$ is the stack of $X\lcom$-torsors.

\begin{defn}\Label{defn.gerbe}
An $\SS$-stack  $\RR$, endowed with a morphism $\RR\to\XX$ is called a
{\bf gerbe }over $\XX$, if both $\RR\to\XX$ and $\RR\to\RR\times_\XX\RR$
are epimorphisms.
\end{defn}

\begin{rmk}
Under  the  correspondence between $\SS$-stacks
equipped  with  morphisms to
$\XX$ and $\XX$-stacks, the gerbes over $\XX$, according to our
definition,  correspond to gerbes over  the site $\XX$ in the usual sense,
i.e. in the sense of Giraud~\cite{Giraud}, Chapter~III.2.
\end{rmk}


$BS^1\times \XX\to \XX$ is an example of a gerbe over
$\XX$. We will study gerbes that locally look like this example.

The groupoid of automorphisms of $BS^1$ is equal
to the transformation groupoid of $S^1$ on $\Aut S^1\cong \zz_2$.
This action is by ``inner automorphisms'' and hence trivial,
as $S^1$ is abelian. The group of automorphism
classes of $BS^1$ is therefore equal to $\zz_2$. The sheaf
of automorphism classes of $BS^1\times \XX$ over $\XX$, which takes
$U/\XX$ to the 2-isomorphism classes of diagrams
$$ \xymatrix{BS^1 \times U\rto^\sim \drto\drtwocell\omit{<-2>} &BS^1 \times U\dto \\&U}$$
is therefore equal to $\zz_2\times \XX\to \XX$.
So if the gerbe $\RR \to \XX$ is locally isomorphic to
$BS^1\times \XX$, then the sheaf of automorphism
classes of $\RR$ over $\XX$, which maps $U/\XX$
to the 2-isomorphism classes of diagrams 
$$ \xymatrix{\RR| U\rto^\sim \drto\drtwocell\omit{<-2>} &\RR| U\dto \\&U}$$
  is a 2-sheeted covering $\band (\RR) \to \XX$, called the {\em band}
of $\RR$.

\begin{defn}
An {\bf $S^1$-gerbe }over $\XX$ is a gerbe $\RR\to \XX$
which is locally isomorphic to $BS^1\times \XX$ and is
endowed with a trivialization of its band (the 2-sheeted covering $\band (\RR) \to \XX$).
%
\end{defn}


%

The following is a well-known theorem of Giraud  \cite{Giraud}.

\begin{them}[Giraud]
\label{thm:giraud}
Isomorphism classes of $S^1$-gerbes over $\XX$ are in one-to-one
correspondence with $H^2\big(\XX, S^1 \big)$.
\end{them}

Now we recall $S^1$-central extensions of Lie groupoids \cite{WeinsteinXu}.

\begin{defn}
Let $X_1\toto X_0$ be a Lie groupoid. An {\em $S^1$-central
extension }of  $X_1\toto X_0$ consists of

1. a  Lie groupoid ${R_1}\toto  X_0$, together with a morphism of Lie
groupoids $(\pi,\id):[R_1\toto  X_0]\to[X_1\toto X_0]$,

2. a left $S^1$-action on $R_1$, making $\pi:R_1\to X_1$ a (left)
principal $S^1$-bundle.

\noindent These two structures are compatible in the sense that
$(s\cdot x) (t\cdot y)=st \cdot (xy )$,
for all  $ s,t \in S^{1}$ and $(x, y) \in R_1\times_{X_0} R_1 $.
\end{defn}

The proposition below gives  an equivalent definition.

\begin{prop}
Let $X_1\toto X_0$ be a Lie groupoid. A Lie groupoid
$R_1\toto X_0$ is an $S^1$-central extension of
$X_1\toto X_0$ if and only if it is endowed with $\phi$ and $\pi$
forming an
exact sequences of groupoid  morphisms
$$1\to X_0\times S^1\stackrel{\phi}{\to} R_1\stackrel{\pi}{\to}X_1\to 1$$
over the identities on the unit spaces,
and the image of $\phi$ lies  in the center of $R_1$.
\end{prop}
\begin{pf}
The proof is straightforward and is left to  the reader.
\end{pf}

The following result describes the precise
 connection between
$S^1$-gerbes and  $S^1$-central extensions.

\begin{prop}
\label{pro:M-trivial}
Let $X_1\toto X_0$ be  a Lie groupoid and $\XX$ its corresponding 
differentiable stack of $X\lcom$-torsors.
There is 
a  one-to-one correspondence
 between isomorphism classes of 
$S^1$-central extensions of $X_1\toto X_0$ and isomorphism classes of 
$S^1$-gerbes $\RR$ over $\XX$ endowed with a trivialization of the
restriction of $\RR$ to $X_0$.
\end{prop}
\begin{pf}
Given an $S^1$-central extension $R_1\toto X_0$ of $X_1\toto X_0$, let
$\RR$ be the stack of $R\lcom$-torsors and $\XX$ the stack of
$X\lcom$-torsors. Then the groupoid morphism $\pi$ induces a morphism
of stacks $\RR\to\XX$, via which we think of $\RR$ as a stack over
$\XX$.

The groupoid morphism $S^1\times X_0\stackrel{\phi}{\to} R_1$ induces the morphism of
stacks $BS^1\times X_0\to\RR$. Consider the diagram
$$\xymatrix{
R_1\rto\dto\ar@{}[dr]|{(1)} & X_0\dto  &\\
X_1\rto\dto\ar@{}[dr]|{(2)} & BS^1\times X_0\rto\dto\ar@{}[dr]|{(3)}
& X_0\dto\\
X_0\rto & \RR\rto & \XX}$$
The square (1) is cartesian, because $R_1|X_1$ is an
$S^1$-torsor. The combination of squares (1) and (2) is cartesian by
definition of $\RR$. Hence, by descent, (2) is cartesian. The
combination of squares (2) and (3) is cartesian by definition of
$\XX$, and so, again by descent, square (3) is cartesian. This proves
that $\RR$ restricted to $X_0$ is isomorphic to $BS^1\times X_0$, and
in particular, $\RR\to\XX$ satisfies the first condition in the
definition of $S^1$-gerbe.

The band of $\RR\to\XX$ is an $\Out(S^1)$-torsor, trivialized by
$X_0$, so the band is given by a map $X_1\to\Out(S^1)$. It is given as
follows: $x\mapsto [Ad_{\tilde{x}}]$, $\forall x\in X_1$,
where $\tilde{x}\in R_1 $ is any point satisfying $\pi (\tilde{x} )=x$
and $Ad_{\tilde{x}} y=\tilde{x}y \tilde{x}^{-1}$.
Here $y\in \ker \pi_x\cong S^1$.
Then because $\ker \pi$ is central in $R_1$, the map $X_1\to\Out(S^1)$
is trivial, showing that the band of $\RR$ is trivial. 

Conversely, given such a gerbe $\RR$,  by taking the section
$X_0\to \RR\resto X_0$, one obtains a commutative diagram
of stacks:
\begin{equation}
\label{eq:xrx}
\xymatrix{
X_0\drto\rto & \RR\dto\\
& \XX .}
\end{equation}
 So $X_0 \to \RR $ is a presentation. 
Let $R_1=X_0\times_\RR X_0$. Thus we have a Lie groupoid
morphism $(\pi,\id):[R_1\toto  X_0]\to[X_1\toto X_0]$.
Moreover the kernel of $\pi$ is isomorphic to $X_0 \times S^1$ as
a bundle of groups, by assumption. Since $\band \RR$ is trivial,
it follows that   the conjugation action of $R_1$ on $\ker \pi$
must be trivial. Therefore $\ker \pi$ lies in the center of
$R_1$. This concludes the proof.
\end{pf}

\subsection{Morita equivalence of $S^1$-central extensions}

We now  introduce the definition of Morita equivalence of $S^1$-central
 extensions.

\begin{defn}
We say that  two $S^1$-central extensions
 $ R_1\to X_1\toto X_0$
 and $R'_1\to X'_1 \toto X'_0$ are Morita equivalent if
there exists an  $S^1$-equivariant $R\lcom$-$R'\lcom$-bitorsor
$Z$, by which we mean that $Z$ is an $R\lcom$-$R'\lcom$-bitorsor
endowed with an $S^1$-action such that
$$(\lambda r)\cdot z\cdot r'= r\cdot (\lambda z)\cdot r' = r\cdot
z \cdot
(\lambda r')$$
 whenever $(\lambda, r, r', z)\in S^1\times R\times R'\times Z$ and
the products make sense.
\end{defn}

The following result is immediate.

\begin{lem}
\label{lem:eq}
Let $R_1\to X_1\toto X_0$
 and $ R'_1\to X'_1 \toto X'_0$ be Morita 
equivalent $S^1$-central extensions, and $Z$ an
$S^1$-equivariant $R\lcom$-$R'\lcom$-bitorsor. Then
the $S^1$-action on $Z$ must be free and $Z/S^1$
is a $X\lcom$-$X'\lcom$-bitorsor.
As a consequence, $X\lcom$ and $X'\lcom$ are Morita
equivalent.
\end{lem}

\begin{prop}
 Let $R_1\to X_1\toto X_0$
 and $R'_1\to X'_1 \toto X'_0$ be $S^1$-central extensions
of Lie groupoids. Let $\RR$, $\RR'$, $\XX$, and $\XX'$
be their associated stacks. Then the following are equivalent:

(i) $R\lcom$ and $R'\lcom$ are Morita
equivalent $S^1$-central extensions.

(ii) there exists an $S^1$-central extension
$R''_1\to X''_1\toto X''_0$  and 
$S^1$-equivariant Morita morphisms  $R''\lcom\to R'\lcom$
and $R''\lcom\to R\lcom$.

(iii) $\XX\cong\XX'$, over which $\RR\cong \RR'$ as $S^1$-gerbes.
\end{prop}
\begin{pf}
To prove that (i) implies (ii), choose an $S^1$-equivariant
$R\lcom$-$R'\lcom$-bitorsor $Z$.  Take $R''_1 =R_1 \times_{X_0}
Z\times_{X'_0} R'_1$ and $X''_1 =X_1 \times_{X_0}
X\times_{X'_0} X'_1$, where $X=Z/S^1$. Then it is
simple to see that $R''_1\to X''_1 \toto X$ is the desired
$S^1$-central extension.

For (ii) to imply (iii), assuming that $R''_1\to X''_1 \toto X''_0$
is such an $S^1$-central extension,
  using Theorem \ref{morita-them},
one has the   commutative diagram
$$\xymatrix{\RR'\rto\dto\drtwocell\omit & \RR\dto\\ \XX'\rto & \XX},$$
where the horizontal  maps are isomorphism of stacks. (iii) thus follows.

Finally, we prove that (iii) implies (i).  By identifying
$\XX$ with $\XX'$ and $\RR$ with $\RR'$, we may
think  $ R_1\to X_1\toto X_0$
 and $R'_1\to X'_1 \toto X'_0$ as the
$S^1$-central extensions corresponding to the presentations
$X_0\to \RR$ and $X'_0 \to \RR$ respectively.
Take $Z=X_0\times_\RR X'_0$. Then  $Z$ is an
$S^1$-equivariant $R\lcom$-$R'\lcom$-bitorsor.
\end{pf}

We end this subsection by the   following exact sequences:

\begin{prop}
\label{pro:ex1}
There is a natural exact sequence
\begin{multline*}
H^1(X\lcom,S^1)\stackrel{\tau_1}{\longrightarrow} H^1(X_0 , S^1)
\stackrel{\tau_2}{\longrightarrow} \\
\{\text{$S^1$-central extensions of $X_1\toto X_0 $}\}
\stackrel{\tau_3}{\longrightarrow}
 H^2(X\lcom,S^1) \stackrel{\tau_4}{\longrightarrow} H^2(X_0 ,S^1)\,.
\end{multline*}
\end{prop}
\begin{pf}
Let $\XX$ be the stack of $X\lcom$-torsors.
 Note that $\tau_1, \cdots , \tau_4$ can be geometrically described as
follows:
\begin{enumerate}
\item $\tau_1 $ is the map sending an $S^1$-bundle $L\to \XX$
to its restriction to $X_0$, i.e. to $L\to X_0$ by forgetting the groupoid
 $X_1\toto X_0$-action.
\item $\tau_2 $ sends an  $S^1$-bundle $L\to X_0$
to the $S^1$-central extension $(s^*L\times_{X_1} {t^*L^{-1}})/S^1\to
X_1\toto X_0$. In stack language, $\tau_2$ maps the $S^1$-bundle $L$
to the stack of descent data (gluing data) for $L$ over the groupoid $X\lcom$.
\item $\tau_3 $ sends an $S^1$-central extension to
the class in $H^2(X\lcom  , S^1 )$ of its corresponding
gerbe.
\item $\tau_4$ is the pull back map under the map $X_0\to \XX$.
\end{enumerate}

Let $\phi :L\to X_0$ be an $S^1$-bundle over $X_1\toto X_0$. Define
a map $ X_1 \times S^1
 \to (s^*L\times_{X_1} {t^*\bar{L}})/S^1$
by $(r, \lambda )\to [( (r l) \lambda^{-1} , l)]$,  where $l \in L$ 
is any point satisfying $\phi (l )=t(r)$. One checks that
this is an isomorphism of $S^1$-central extensions.
Conversely, if  $\phi :L\to X_0$ is an $S^1$-bundle over 
$X_0$ such that $(s^*L\times_{X_1} {t^*\bar{L}})/S^1$ is a
 trivial central extension, then $(s^*L\times_{X_1} {t^*\bar{L}})/S^1\to
X_1$ admits a  section $\sigma $  which is 
a groupoid   homomorphism. Then the equation $\sigma (r)=[(r\cdot l, l)]$,
where $l\in L$ such that $\phi (l )=t(r)$, defines 
a groupoid action of $X_1\toto X_0$ on $L$.
This shows that the sequence is exact at 
$H^1(X_0 , S^1)$.

Let $R$ be     the $S^1$-central extension
 $(s^*L\times_{X_1} {t^*\bar{L}})/S^1$, where $L\to X_0$
is an $S^1$-bundle. One checks that the pullback groupoid
$R[L]\toto L$ is a trivial extension of $X_1 [L]\toto L$,
which means that $R$ defines the zero class in $H^2(X\lcom , S^1 )$.
Conversely, if $R$ is an $S^1$-central extension 
defining the zero class in $H^2(X\lcom , S^1 )$, then $R$ is
Morita equivariant to the trivial central extension 
$X_1 \times S^1 \toto X_0$ via an $S^1$-equivariant
bimodule $Y$: 
then $ L=Y/X_1\stackrel{\rho}{\to}X_0$ is an $S^1$-bundle.
One checks easily that $R$ is isomorphic to 
$(s^*L\times_{X_1} {t^*\bar{L}})/S^1$.
This shows that the sequence is exact at  
$\{\text{$S^1$-central extensions of $X_1\toto X_0 $}\}$.

Finally the exactness at $H^2(X\lcom, S^1 )$  follows from
Theorem \ref{thm:giraud} and Proposition
\ref{pro:M-trivial}.
\end{pf}

\subsection{Dixmier-Douady  classes} 
Let $R$ be an $S^1$-central extension of $X_1\toto X_0$.  Write the
underlying Lie groupoid of $R$ as $R_1\toto R_0$. Call the structure
morphism $\pi:R\lcom\to X\lcom$. Since $R_1\toto R_0$
defines an $S^1$-gerbe over $\XX$, it defines a class in
$H^2(X\lcom,S^1)$ according to Theorem \ref{thm:giraud}.
The exponential sequence gives rise to a homomorphism
$H^2(X\lcom,S^1)\to H^3(X\lcom,\zz)$.  The image of
$[R]\in H^2(X\lcom,S^1)$ in $H^3(X\lcom,\zz)$ is called the
{\em Dixmier-Douady class }of $R$ and denoted by $DD (R)$.
 The Dixmier-Douady class behaves well
with respect to pullbacks and the tensor operation.

Let $f: Y\lcom \to X\lcom$ be a Lie  groupoid homomorphism.
 Then  the pullback $S^1$-bundle $f^*R_1\to Y_1$ is 
an  $S^1$-central extension over $Y_1\toto Y_0$, called the pullback
central extension.

Assume that  $R'$  and  $R''$
are two $S^1$-central extensions of $X_1\toto X_0$.
Let $R_1=(R'\times_{X_1}R'' )/S^1$, where $S^1$ acts on
$R'\times_{X_1}R''$ by
$t\cdot (r_1, r_2 )=(t\cdot r_1, t^{-1}\cdot  r_2 ),
\  \forall t\in S^1, \ \ (r_1, r_2)\in R'\times_{X_1}R'' $.
It is clear that $R_1$ with the natural projection
to $X_1$ is still an $S^1$-principal bundle, where the
$S^1$-action is given by $t\cdot [ (r_1, r_2 )]=[(t\cdot r_1 , r_2)]$.
The groupoid structures on $R'$ and $R''$ induce
a groupoid structure on $R_1$ in a natural way, which
in fact makes $R_1$ into a groupoid  $S^1$-central extension,
called the tensor product of $R'$ and $R''$ and is  denoted,
by $R'\otimes R''$.

The following proposition can be easily verified.

\begin{prop}
\begin{enumerate}
\item $DD(f^*R)=f^* DD(R)$; and 
\item $DD(R'\otimes R'')=DD(R')+DD(R'').$
\end{enumerate}
\end{prop} 

\begin{defn}
Let $\theta\in\Omega^1(R_1)$ be a connection  1-form for the 
$S^1$-principal bundle
$R_1\to X_1$, and $B\in\Omega^2(R_0)$ be any  2-form. Any such pair
$(\theta,B)$ is called a {\em pseudo-connection }for the central
extension $R$.
\end{defn}

It is simple to check that $\delta(\theta+B) \in Z^3_{DR}(R\lcom)$
 descends to $Z^3_{DR}(X\lcom)$, i.e. there exist unique
$\eta\in\Omega^1(X_2)$, $\omega\in\Omega^2(X_1)$ and
$\Omega\in \Omega^3(X_0)$ such that
$$  \delta(\theta+B)= \pi\upst(\eta+\omega+\Omega).$$

Then $\eta+\omega+\Omega$ is called the
{\em  pseudo-curvature }of the pseudo-connection $\theta+B$.

We will now show that pseudo-connections can be used to calculate
Dixmier-Douady classes.

\begin{them}
\label{pro:DD}
The class $[\eta+\omega+\Omega]\in H^3_{DR}(X\lcom)$ is independent of the
choice of the pseudo-connection $\theta+B$. Under the canonical
homomorphism $H^3(X\lcom,\zz)\to H^3_{DR}(X\lcom)$, the Dixmier-Douady
class of $R$ maps to $[\eta+\omega+\Omega]$.
\end{them}
\begin{pf}
One checks directly that the class $[\eta+\omega+\Omega]\in H^3_{DR}(X\lcom)$ is independent of the
choice of the pseudo-connection $\theta+B$.

We prove the second part of the theorem.
Let $\XX$ be the stack given by $X_1\toto X_0$ and $\RR\to\XX$ the
$S^1$-gerbe over $\XX$ defined by $R_1\toto R_0$.

We will construct a hypercovering in the site $\XX$. Note that
$R\lcom$ is a simplicial object in $\XX$. The hypercovering we shall
use is the 1-coskeleton of $R\lcom$:
$$Y\lcom=\cosk_\XX R\lcom$$
This is a hypercovering because
$R_1\to R_0\times_\XX R_0=X_1$ and $R_0\to\XX$ are surjective
submersions. (For the theory of hypercoverings, see \cite{sga4, art.maz, deligne}.
 In the generality we need them, the necessary results are
 proved in~\cite{sga4}.)

Intuitively, $Y\lcom$ is the set of all (1-skeleta of) simplices in
$R\lcom$, whose image in $X\lcom$ commutes. 
  Explicitly, $Y_p$ is the
fibered product
$$\xymatrix{
Y_p\rrto\dto && {\displaystyle\prod_{0\leq i<j\leq p}^{\phantom x}  R_1}
\dto\\ 
X_p\rrto & &    {\displaystyle\prod_{0\leq i<j\leq p}^{\phantom x}
  X_1}}$$ 
Here the horizontal arrow at the bottom is the map which sends
$X_p$ to the edges of a commutative $p$-simplex, i.e. the
product of $\big({p+1\atop 2}\big)$   maps $f_{ij}: X_p\to X_1,
\  0\leq i<j\leq p$
$$ (x_1, x_2, \cdots , x_p) \to x_{i+1}\cdots x_j\,.$$

Since $Y\lcom$ is a hypercovering of $\XX$, we have a canonical
homomorphism
\begin{equation}
\label{eq:7}
f:\check H^2(Y\lcom,S^1)\tto  H^2(\XX,S^1).
\end{equation}
Since $\Omega\com$ consists of soft sheaves,
we also have an isomorphism
\begin{equation}
\check H^2(Y\lcom, \Omega\com )\longiso  H^2(\XX, \Omega\com).
\end{equation}
We will see that the class of $\RR$ in $H^2(\XX,S^1)$ is in the image
of the homomorphism \eqref{eq:7}.

In fact, 
$$Y_2=\{(\alpha,\beta,\gamma)\in R_1\times R_1\times R_1\st
\pi(\alpha ) \pi(\gamma)=\pi(\beta)\}\,,$$
so we have a $\cinf$-map
\begin{align*}
c:Y_2&\tto S^1\\
(\alpha,\beta,\gamma)&\longmapsto (\alpha \gamma)\beta^{-1}\,.
\end{align*}
Recall that a composition in $R_1$ makes sense if and only if
the composition of its image in $X_1$ makes sense and that we have $\ker\pi=S^1$.
One checks that the  coboundary of $c$ vanishes, and so $c$ defines
a \v Cech cohomology class $[c]\in \check H^2(Y\lcom,S^1)$. 
It is simple to see that
 $f([c])$ is the cohomology class of $\RR$.

Now consider the diagram
$$\xymatrix{
& {\check H^2(Y\lcom,S^1)} \rto^f\dto & {H^2(\XX,S^1)}\drto^{\del}\dto & \\ 
& {\check H^3\big(Y\lcom,[\Omega^0\to S^1]\big)} \rto\dto^{d\log} &
  {H^3\big(\XX,[\Omega^0\to S^1]\big)}\dto^{d\log} &
  H^3(\XX,\zz)\dto\lto_-\sim\\ 
{\check H^3(X\lcom,\Omega\com)}\rto^\sim_{\rho\upst} & {\check
  H^3(Y\lcom,\Omega\com)} \rto^\sim & 
  {H^3_{DR}(\XX)} & H^3(\XX,\rr)\lto_\sim}$$
which commutes. 
 The two vertical arrows in the first row are induced
by the trivial map 
$$S^1\tto[\Omega^0\to S^1][1]\,,$$
i.e. the map
$$\xymatrix{ 0\rto \dto &S^1 \dto\\
\Omega^0 \rto &S^1} $$ 

Considering this diagram, we see that we need to prove that
$$d\log([c])=\rho\upst([\eta+\omega+\Omega])\in {\check
  H^3(Y\lcom,\Omega\com)}\,,$$  
where we have denoted the canonical projection by $\rho:Y\lcom\to X\lcom$
and its induced map on \v{C}ech  cohomology by $\rho\upst$. We have also
committed the abuse of denoting $[c]$ and its induced class in ${\check
  H^3\big(Y\lcom,[\Omega^0\to S^1]\big)}$ by the same letter.

First we may assume that $B=0$ (thus $\Omega=0$) for simplicity since the 
class $[\eta+\omega+\Omega]$ is independent of the pseudo-connection.
Thus we have
$$\partial \theta =\pi^* \eta , \ \ \ d\theta =-\pi^* \omega .
$$
We have the following commutative diagram:
$$\xymatrix@C=1pc{Y_2\rrto^p\drto && R_2\dlto\\ & X_2 &}$$
where $p: Y_2\to R_2$ is the natural projection.
We have
$$\rho\upst\eta= p\upst \pi\upst \eta 
=p\upst \partial \theta =\alpha^* \theta-(\alpha \gamma )^* \theta
+\gamma^* \theta,
$$
where, by abuse of notation, we denote by $\alpha, \ \gamma $ and
$\alpha \gamma$ the maps $Y_2\to Y_1$ sending
$(\alpha,\beta,\gamma)$ to $\alpha, \ \gamma $ and $\alpha \gamma \in Y_1$,
respectively.

Since $\rho=\pi$  on $ Y_1$,  we have $\rho^* \omega 
=\pi^*\omega =-d\theta
\in \Omega^2 (Y_1)$, which is cohomologous to $-\partial_{Y\lcom} \theta$
in $\check H^3 (Y\lcom,\Omega\com)$. The latter is
equal to $-(\alpha^* \theta -\beta^* \theta +\gamma^* \theta)\in
\Omega^1 (Y_2)$. Thus it follows that
$$ \rho\upst([\eta+\omega ])=\beta^* [\theta-(\alpha \gamma )^* \theta ].$$
Now it suffices to prove that
$$d\log([c])=\beta^* \theta-(\alpha \gamma )^* \theta \in \Omega^1 (Y_2). $$
Let  $\psi : R_2\times S^1 \to Y_2$ be the diffeomorphism
given by $(\alpha, \gamma, t)\to (\alpha, t(\alpha \gamma),\gamma)$.
Then   $\psi^\ast ( d\log([c]) )$ is the Maurer-Cartan form $dt$ on $S^1$,
while $ \beta^* \theta-(\alpha \gamma )^* \theta $ is easily seen
to be equal to $dt$ as well.

This completes the proof.
\end{pf}

\subsection{Prequantization}

\begin{defn}
Given an $S^1$-central extension $R_1 \to X_1\toto X_0$,

(i) a connection $1$-form $\theta\in \Omega^1(R_1)$ 
for the $S^1$-principal  bundle $R_1\to X_1$, such that
$\del\theta=0$ is a {\em connection};

(ii) Given $\theta$, a 2-form $B\in\Omega^2(X_0 )$, such that $d\theta=\del B$
is a {\em curving}; 

(iii) and given $(\theta,B)$, the 3-form 
$\Omega=dB\in H^0(X\lcom,\Omega^3)\subset \Omega^3(X_0 )$ is called
the {\em $3$-curvature }of $(\theta,B)$;

(iv)  If $\Omega=0$, then $R_1 \to X_1\toto X_0$ together with
$(\theta,B)$ is called a {\em flat } $S^1$-central extension of
$X_1\toto X_0$.  Note that the flat central extensions form an abelian
group.
\end{defn}

In other words, a flat $S^1$-central extension of
$X_1\toto X_0$ is an $S^1$-central extension with a
pseudo-connection whose pseudo-curvature vanishes.
The following proposition is  immediate.

\begin{prop}
Let $R_1 \to X_1\toto X_0$ be an $S^1$-central extension. Then

(i)  $H^2 (X\lcom , \Omega^1)$ contains the obstruction to the
existence of a connection;

(ii) if we assume the existence of a connection, 
$H^1 (X\lcom , \Omega^2)$ contains the obstruction to the
existence of a curving.
\end{prop}

According to Theorem \ref{pro:DD}, we have
the following

\begin{prop}
 If an $S^1$-central extension  $R_1 \to X_1\toto X_0$ 
admits a connection and a curving with 3-curvature $\Omega$, then
$[\Omega]\in H^3_{DR}(X\lcom)$ is the image of  its Dixmier-Douady class
under the canonical
homomorphism $H^3(X\lcom,\zz)\to H^3_{DR}(X\lcom)$.
\end{prop}

\begin{rmk}
Given a manifold $M$, and a surjective submersion
$X_0\to M$,  $X_1 (=X_0\times_M X_0) \toto X_0$ is a Lie
groupoid Morita equivalent to $M$. An $S^1$-central extension $R_{1}\to X_1\toto X_0$
defines a bundle gerbe over $M$ in the terminology of Murray \cite{Murray, MurrayS}.
Since $\Omega^1$ and $\Omega^2$ are soft sheaves over $M$,\comment{check tem}
we have  $H^2(X\lcom , \Omega^1)\cong H^2(M, \Omega^1)=0$
and $H^1(X\lcom , \Omega^2)\cong H^1(M, \Omega^2)=0$.
As a consequence, connections and curvings always exist
for bundle gerbes. This result was due to Murray \cite{Murray}.
Moreover, in this case,  the $3$-curvature $\Omega\in
\Omega^3(X_0 )$ descends to a closed 3-form
on $M$ since $\partial \Omega=0$.

In particular, for an open cover $\{U_i\}$ of $M$, one can take
$X_0=\coprod U_i$.
Then $X_1\cong \coprod U_{ij}$.
An $S^1$-central extension $R_{1}\to X_1\toto X_0$
corresponds in this case  to a family of line bundles $L_{ij}\to U_{ij}$
satisfying all the axioms of bundle gerbes as  in \cite{Hitchin}.
This is the case of Chatterjee-Hitchin bundle gerbes \cite{Chatterjee, Hitchin}
\end{rmk}

\begin{prop}
Assume that $H^2(X_0, \rr )=0$.
There is a natural exact sequence
\begin{multline*}
H^1(X\lcom,\rr/\zz)\longrightarrow H^1(X_0 ,\rr/\zz)\longrightarrow
\\ \{\text{flat $S^1$-central extensions of $X_1\toto X_0$}\}
\longrightarrow H^2(X\lcom,\rr/\zz) \longrightarrow
H^2(X_0,\rr/\zz)\,.
\end{multline*}
\end{prop}
The proof is similar to that of Proposition  \ref{pro:ex1}
via replacing $S^1$ by $\rr/\zz$, and using the following

\begin{lem}
Let $X_1\toto X_0$ be a Lie groupoid. Assume that
 $H^2(X_0, \rr )=0$. Then there is a canonical 
one-to-one correspondence between flat $S^1$-central
extensions of $X_1\toto X_0$ and $\rr/\zz$-central
extensions of $X_1\toto X_0$.  
\end{lem}
\begin{pf}
Let $(R_1\to X_1\toto X_0, \theta , B)$ be a  flat $S^1$-central
extension. Then, in particular, $dB=0$.
 Since $H^2(X_0, \rr )=0$, we can write $B=dA$, where
$A\in \Omega^1 (X_0 )$. Set $\theta'=\theta +\partial A \in 
\Omega^1 (R_1 )$. Then $\theta'$ is again a connection 1-form
for the principal $S^1$-bundle $R_1\to X_1$, which satisfies $d\theta'
=0$ and $\partial \theta'=0$. 
The condition $d\theta' =0$ implies that
 $R_1\to X_1$ is flat, and can therefore equivalently be considered as
an $\rr/\zz$-bundle. Moreover, $\partial \theta'=0$
implies that under this new  differentiable  structure,
$R_1\to X_1$ is still a smooth groupoid homomorphism, and therefore
an $\rr/\zz$-central extension.

Conversely, given an $\rr/\zz$-central extension $R_1\to X_1\toto X_0$,
then $R_1\to X_1$ is a flat $S^1$-bundle. Let $\theta \in
\Omega^1 (R_1)$ be a flat connection one-form, i.e. $d\theta =0$.
Locally, if we write $R_1\cong X_1 \times \rr/\zz$, then  we may choose
$\theta =dt$ where $t$ is the coordinate on $\rr/\zz$.
Moreover locally the groupoid multiplication on $R_1$ is written
as 
$$(x, t)\cdot (y, s)=(x\cdot y, t+s +\omega (x, y)), \ \ \forall (x, y )\in
X_2, \ \ t, s \in \rr/\zz.$$
It is easy to see that $\omega (x, y)$ must be locally 
constant. Therefore  it follows that $\partial \theta =0$.
Hence $R_1\to X_1\toto X_0$ is a flat $S^1$-central extension.
\end{pf}

Following Hitchin \cite{Hitchin}, we call
the map
$$ \{\text{flat $S^1$-central extensions of $X_1\toto X_0$}\}
\longrightarrow H^2(X\lcom,\rr/\zz)$$
{\em the holonomy map}.

Next we give the following prequantization result, which can be
considered as an analogue, in the degree $3$-context, of
the well known theorem of Weil and Kostant \cite{Kostant, Weil}.

\begin{them}
Assume that $\check  H^2 (X\lcom , \Omega^0 )=0$.
Given any  3-cocycle  $\eta +\omega+\Omega\in Z^3_{DR}(X\lcom)$
as above, satisfying

1. $\eta+\omega+\Omega$ is  an integer 3-cocycle, and

2. $\Omega$ is exact,

\noindent there exists a groupoid $S^1$-central extension $R_1\toto X_0$ of the
groupoid $X_1 \toto X_0$ and  a  pseudo-connection
 $\theta +B$ such that its pseudo-curvature
is $\eta+\omega+\Omega$.
The pairs $(\theta, B)$ up to isomorphism form a simply transitive
set under the group of flat $S^1$-central extensions.
\end{them}
\begin{pf}
Consider the exact sequence
$$\to H^2(X\lcom,S^1)\stackrel{\phi}{\to}
  H^3(X\lcom,\zz)\to
H^3(X\lcom,\Omega^0)\to$$
induced by the  exponential sequence $\zz\to \Omega^0\to S^1$.
Since we have the following commutative diagram
$$\xymatrix@C=1pc{H^3 (X\lcom , \zz )\rrto\drto
 && H^3_{DR}(X\lcom )\dlto^p\\ &  H^3 (X\lcom,\Omega^0) \cong
{\check H}^3 (X\lcom,\Omega^0)  &}$$
where  $p$ is the natural projection, it is clear that
$[\eta+\omega +\Omega ]$ is in the image of $\phi$. Thus
there is an $S^1$-gerbe   $\RR \in H^2(X\lcom , S^1 )$ 
 whose  Dixmier-Douady   class equals   $[\eta +\omega+\Omega]$.     
Note that the image of $\RR$ under the map
$H^2(X\lcom , S^1 )\to H^2(X_0, S^1 )$ is zero since $\Omega$
is exact. This follows from the commutative diagram
$$\xymatrix{ H^2(X\lcom , S^1 )\rto \dto &
H^3(X\lcom, \zz ) \dto\\
H^2(X_0, S^1 )  \rto^\sim & H^3(X_0 , \zz )  } $$  
From Proposition  \ref{pro:ex1} it follows that
$\RR$ can be represented by  an $S^1$-central extension $R_1\toto R_0$
over $X_1\toto X_0$,    whose Dixmier-Douady
 class is   $[\eta +\omega+\Omega]$. 

Let $\theta' +B'$ be any pseudo-connection on
the $S^1$-central extension $R_1\toto R_0$  and 
$\eta' +\omega'+\Omega'$ its pseudo-curvature.
 Proposition \ref{pro:DD}   implies that $\eta +\omega+\Omega$
and $\eta' +\omega'+\Omega'$  are cohomologous.
Therefore 
$$(\eta +\omega+\Omega)-(\eta' +\omega'+\Omega' )=
\delta (f+\alpha + B''), $$
where $f\in  \Omega^0 (X_2 ), \ \alpha \in \Omega^1 (X_1 )$, and $
 B\in  \Omega^2 (X_0 )$.
It thus follows that
$\partial f=0$, which implies that $f=\partial g$ for
$g\in  \Omega^0 (X_1)$ since $\check H^2 (X\lcom , \Omega^0 )=0$. Thus
$\delta f=\delta \partial g=\delta dg $. Let
$\theta =\theta'+\pi^*(\alpha +dg)\in \Omega^1 (R_1 )$ and
$ B=B'+B''$.
It is clear that $\theta +B$ is the desired pseudo-connection
on $R_1\toto R_0$.     
Finally note that if   $(R', \theta', B')$ and 
$(R'', \theta'', B'')$ are two such pairs, then 
 $(R'\otimes (R'')^{-1} ,\theta'-\theta'' , B'-B'' )$
is a flat gerbe. So such pairs, up to isomorphism,
are indeed parametrized by the group of flat $S^1$-central
extensions. 
\end{pf}

\begin{rmk}
Note again, that the condition $\check     H^2 (X\lcom, \Omega^0 )=0$
always holds for a proper Lie groupoid $X_1\toto X_0$,
according to Crainic \cite{Crainic}. So prequantization
always works for a proper Lie groupoids. 
\end{rmk}

\def\nin{\noindent}
\def\eg{e.g.\ }
\def\ie{i.e.\ }
\def\pt#1#2{{\partial #1\over \partial #2}}

\newcommand{\alp }{\alpha }
\newcommand{\bet }{\beta }
\newcommand{\gm }{\Gamma }
\newcommand{\lon }{\longrightarrow }
\newcommand{\be }{\begin{eqnarray*}}
\newcommand{\ee }{\end{eqnarray*}}
\newcommand{\gt }{ {\frakg}\otimes  {\frakg} }
\newcommand{\per }{\backl }
\newcommand{\te }{\otimes  }
\newcommand{\poiddd }[3]{ (#1\gpd #2, s, t)}
\newcommand{\tu}{\tilde{u}}
\newcommand{\tm}{\tilde{m}}
\newcommand{\tv}{\tilde{v}}
\newcommand{\tg}{\tilde{g}}
\newcommand{\tee}{\tilde{\epsilon}}
\newcommand{\dress}{\lambda_{\psi (v^{-1})}u}
\newcommand{\dre}[1]{\lambda_{#1}}
\newcommand{\gmstar}{\gm_{G^{*}}}
\newcommand{\complex}{{\Bbb      C}}
\newcommand{\reals}{{\Bbb      R}}
\newcommand{\integers}{{\Bbb   Z}}
\newcommand{\ctn}{\complex^{2n}}
\newcommand{\rtn}{\reals^{2n}}
\newcommand{\rn}{\reals^{n}}
\newcommand{\rtm}{\reals^{2m}}
\newcommand{\frakb}{{\frak b}}
\newcommand{\frakd}{{\frak d}}
\newcommand{\frakg}{{\frak g}}
\newcommand{\frakh}{{\frak h}}
\newcommand{\frakl}{{\frak l}}
\newcommand{\fraks}{{\frak s}}
\newcommand{\bstar}{{\frak b}^{*}}
\newcommand{\gstar}{{\frak g}^{*}}
\newcommand{\hstar}{{\frak h}^{*}}
\newcommand{\lstar}{{\frak l}^{*}}
\newcommand{\sstar}{{\frak s}^{*}}
\newcommand{\starh}{*_{\hbar}}
\newcommand{\zbar}{\overline{z}}
\newcommand{\torus}{{\Bbb T}}
\newcommand{\boldphi}{\mbox{\boldmath $\phi$}}
\newcommand{\boldpi}{\mbox{\boldmath $\pi$}}
\newcommand{\inner}{\backl}
\newcommand{\half}{\frac{1}{2}}
\newcommand{\quarter}{\frac{1}{4}}
\newcommand{\canoncoords}{q_{1},\dots,q_{n},p_{1},\dots ,p_{n}}
\newcommand{\canonbracket}{\sum _{i=1}^{n}
[\frac{\del f}{\del q_{i}}\frac{\del g}{\del p_{i}}-
\frac{\del g}{\del q_{i}}\frac{\del f}{\del p_{i}}]}
\newcommand{\symplecticform}{\sum_{i=1}^{n}dq_{i}\wedge dp_{i}}
\newcommand{\poissend}[1]{#1\times\overline{#1}}
\newcommand{\oppositepoissend}[1]{\overline{#1}\times #1}
\newcommand{\threetimes}[1]{#1 \times \overline{#1}\times\overline{#1}}
\newcommand{\fourtimes}[1]{#1\times\overline{#1}\times\overline{#1}\times
#1}
\newcommand{\inverse}{^{-1}}
\newcommand{\boldc}{{\bf c}}
\newcommand{\boldi}{{\bf i}}
\newcommand{\boldj}{{\bf j}}
\newcommand{\boldk}{{\bf k}}
\newcommand{\boldR}{{\bf R}}
\newcommand{\boldS}{{\bf S}}
\newcommand{\boldx}{{\bf x}}
\newcommand{\boldy}{{\bf y}}
\newcommand{\D}{{\mathcal F}} 
\newcommand{\trace}{\mbox{\rm{Tr}}}
\newcommand{\area}{\mbox{\rm{Area}}}
\newcommand{\backl}{\mathbin{\vrule width1.5ex height.4pt\vrule height1.5ex}}
\newcommand{\cala}{{\cal A}}
\newcommand{\calb}{{\cal B}}
\newcommand{\calc}{{\cal C}}
\newcommand{\cald}{{\cal D}}
\newcommand{\cale}{{\cal E}}
\newcommand{\calf}{{\cal F}}
\newcommand{\calg}{{\cal G}}
\newcommand{\calh}{{\cal H}}
\newcommand{\cali}{{\cal I}}
\newcommand{\calj}{{\cal J}}
\newcommand{\calk}{{\cal K}}
\newcommand{\call}{{\cal L}}
\newcommand{\calm}{{\cal M}}
\newcommand{\caln}{{\cal N}}
\newcommand{\calo}{{\cal O}}
\newcommand{\calp}{{\cal P}}
\newcommand{\calq}{{\cal Q}}
\newcommand{\calr}{{\cal R}}
\newcommand{\calu}{{\cal U}}
\newcommand{\calx}{\mathfrak{X}}
\newcommand{\caly}{{\cal Y}}
\newcommand{\defequal}{\stackrel{\mbox {\tiny {def}}}{=}}
\newcommand{\ti}{\tilde{\iota}}
\newcommand{\ts}{{\tilde{s}}}
\newcommand{\tlt}{{\tilde{t}}}
\newcommand{\x}{\tilde{x}}
\newcommand{\y}{\tilde{y}}
\newcommand{\z}{\tilde{z}}
\newcommand{\w}{\tilde{w}}
\newcommand{\dda}{\ddagger}

\newcommand{\smalcirc}{\mbox{\,\tiny{$\circ $}\,}}     

\newcounter{bean}
\newcounter{bacon}
\newcounter{in}
\def\description label#1{\hfil\bf[#1]\hfil}
\parskip 5pt plus 1pt
\topmargin 4pt
\newcommand{\M}{\cal M}
\newcommand{\g}{\frak g}
\newcommand{\h}{\frak h}
\newcommand{\slm}{\boldS_{\lambda}\otimes \boldS_{\mu} }
\newcommand{\rlm}{\calr_{\lambda \mu}}
\newcommand{\rml}{\calr_{\mu \lambda }}
\newcommand{\sla}{S_{\lambda}}
\newcommand{\sml}{\boldS_{\mu}\otimes \boldS_{\lambda} }
\newcommand{\smu}{S_{\mu}}
\newcommand{\gmdouble}{\gmstar \times \gmstar}
\newcommand{\jml}{J_{\mu}\otimes J_{\lambda}}
\newcommand{\jlm}{J_{\lambda}\otimes J_{\mu}}
\newcommand{\jm}{J_{\mu}}
\newcommand{\jl}{J_{\lambda }}
\newcommand{\ssss}{S\times S\lon S\times S}
\newcommand{\hcalr}{\hat{\calr}}
\newcommand{\oo}{o}
\newcommand{\gdouble}{\gr \oplus (\gr^{*})^{\oo}}
\newcommand{\hdouble}{\h \oplus (\h^{*})^{\oo}}
\newcommand{\su}[1]{SU(#1)}
\newcommand{\sbb}[1]{SB(#1, \complex )}
\newcommand{\sll}[1]{SL(#1, \complex )}
\newcommand{\sustar}[1]{\frak{su}^{*}(#1)}
\newcommand{\liesl}[1]{\frak{sl}(#1, \complex )}
\newcommand{\hhdouble}{H\bowtie H^{*}}
\newcommand{\ggdouble}{G\bowtie G^{*}}
\newcommand{\hr}{\hat{R}}
\newcommand{\sect}{\gm(A)}
\newcommand{\sectg}{\gm(AG)}                               
\newcommand{\secdual}{\gm(A^{*})}
\newcommand{\secdualg}{\gm(A^{*}G)}                        
\newcommand{\secc}[1]{\gm(\wedge^{#1}A)}
\newcommand{\seccg}[1]{\gm(\wedge^{#1}AG)}                 
\newcommand{\secd}[1]{\gm(\wedge^{#1}A^{*})}
\newcommand{\secdg}[1]{\gm(\wedge^{#1}A^{*}G)}             

\newcommand{\parr}[1]{\frac{\partial}{\partial x^{#1}}}
\newcommand{\dx}{\dot{x}}
\newcommand{\anc}{a}
\newcommand{\ancd}{a_{*} }
\newcommand{\ot}{\omega \wed \theta}
\newcommand{\xia}{\xi_{a}}
\newcommand{\pib}{\pi^{\#}}
\newcommand{\talp}{\widetilde{\alpha}}           
\newcommand{\tbet}{\widetilde{\beta}}            
\newcommand{\difft}{\frac{d}{dt}}
\newcommand{\btheta}{\bar{\theta}}
\newcommand{\bomega}{\bar{\omega}}
\newcommand{\TLambda}{\widetilde{\Lambda}}
\newcommand{\tpr}{\tilde{pr}}
\newcommand{\xt}{\tilde{x}(t)}
\newcommand{\yt}{\tilde{y}(t)}
\newcommand{\zt}{\tilde{z}(t)}
\newcommand{\ut}{\tilde{u}(t)}
\newcommand{\vt}{\tilde{v}(t)}
\newcommand{\hth}{\tilde{h_1}(t)}
\newcommand{\htt}{\tilde{h_2}(t)}
\newcommand{\dxt}{\dot{\tilde{x}}(t)}
\newcommand{\dyt}{\dot{\tilde{y}}(t)}
\newcommand{\dzt}{\dot{\tilde{z}}(t)}
\newcommand{\dut}{\dot{\tilde{u}}(t)}
\newcommand{\dft}{\dot{{f}}(t)}
\newcommand{\dvt}{\dot{\tilde{v}}(t)}
\newcommand{\dht}{\dot{\tilde{h}}_1 (t)}
\newcommand{\dhtt}{\dot{\tilde{h}}_2 (t)}

\newcommand{\poidd }[2]{#1\toto #2}

\let\Tilde=\widetilde
\let\Bar=\overline
\let\Vec=\overrightarrow
\let\ceV=\overleftarrow
\let\Hat=\widehat



\renewcommand{\theenumi}{{\rm{(\roman{enumi})}}}

\def\Clin{C^\infty_{\ell in}}

\def\mvf{multiplicative vector field}
\def\mmvf{multiplicative multivector field}
\def\cdo{covariant differential operator}
\def\LAgpd{${\cal LA}$-groupoid}
\def\VBLalgd{${\cal VB}$-Lie algebroid}
\def\VBgpd{${\cal VB}$-groupoid}
\def\CDO{\mathop{\rm CDO}}
\def\Ad{\mathop{\rm Ad}}
\def\End{\mathop{\rm End}}
\def\fms#1{\Omega^1(#1)}

\def\ST{\ \vert\ }

\def\tilalpha{\widetilde\alpha}
\def\tilbeta{\skew6\widetilde\beta}

\newcommand{\wed}{\mathbin{\lower1.5pt\hbox{$\scriptstyle{\wedge}$}}}

\let\Tilde=\widetilde
\let\Bar=\overline
\let\Ri=\overrightarrow
\let\Le=\overleftarrow
\def\ssri{{\stackrel{\rightarrow}{X}}}
\def\ssli{{\stackrel{\leftarrow}{X}}}
\def\ssriy{{\stackrel{\rightarrow}{Y}}}
\let\Hat=\widehat
\let\isom=\cong
\let\sol=\bullet
\def\chigh{{\raise1.5pt\hbox{$\chi$}}}
\let\phi=\varphi

\def\til0{\Tilde{0}}
\def\til1{\Tilde{1}}

\def\dpl{\mathbin{+\hskip-6pt +\hskip4pt}}
\def\dminus{\raise2pt\hbox{\vrule height1pt width 2ex}\hskip3pt}
\def\dtimes{\mathbin{\hbox{\huge.}}}

\def\llangle{\langle\!\langle}
\def\rrangle{\rangle\!\rangle}

\def\pback#1{\mathbin{{{\lower1.2ex\hbox{$\times$}}\atop #1}}}

\def\ddt#1{\left.\frac{d}{dt}#1\right|_0}

\def\brev{\:\breve{\rule{0pt}{8pt}}}
\def\hatt{\Hat{\phantom{X}}}

\def\vlra{\hbox{$\,-\!\!\!-\!\!\!-\!\!\!-\!\!\!-\!\!\!
-\!\!\!-\!\!\!-\!\!\!-\!\!\!-\!\!\!\longrightarrow\,$}}

\def\vleq{\hbox{$\,=\!\!\!=\!\!\!=\!\!\!=\!\!\!=\!\!\!
=\!\!\!=\!\!\!=\!\!\!=\!\!\!=\!\!\!=\!\!\!=\!\!\!=\!\!\!=\,$}}

\def\lrah{\hbox{$\,-\!\!\!-\!\!\!
-\!\!\!-\!\!\!-\!\!\!-\!\!\!-\!\!\!\longrightarrow\,$}}

\def\surj{-\!\!\!-\!\!\!-\!\!\!\gg}

\def\inj{>\!\!\!-\!\!\!-\!\!\!-\!\!\!>}

\section{$S^1$-central extensions with prescribed pseudo-curvature}

\subsection{Geometry of $S^1$-central extensions}

First we need a technical lemma concerning 
   $S^1$-principal bundles over a Lie groupoid
(not necessary a groupoid central extension).

Let  $\poidd{X_1}{X_0}$ be a Lie groupoid with  a
3-cocycle  $\eta +\omega \in Z_{DR}^3 (X\lcom )$,
where $\eta \in \Omega^1 (X_2 )$ and $\omega \in \Omega^2 (X_1 )$, and 
$R_1\stackrel{\pi}{\lon}X_1 $ an $S^1$-principal bundle. 
Assume that $\theta \in \Omega^{1} (R_1)$ is a
principal bundle connection one-form with curvature $-\omega$, i.e.
$$d\theta =-\pi^* \omega. $$
 Consider the $T^2$-action on $R_1\times R_1\times R_1$:
\begin{equation}
\label{eq:T2}
(s, t)\cdot (\x, \y, \z)=(s\cdot \x, t\cdot \y, (st)\cdot \z ),
 \ \ \forall s, t \in S^1, \ \x, \y, \z \in R_1.
\end{equation}
Then $p: (R_1\times R_1\times R_1)/T^2 \lon X_1 \times X_1 \times X_1$
is an $S^1$-principal bundle. Consider the  following
diagram of principal bundles

\begin{equation}                         \label{eq:A}
\xymatrix{     
& &T^3\dto&T^1\dto \\
&T^2\rto&\pi^{-1}(\Lambda)\rto^{\tau}\dto^{\pi}&\pi^{-1}(\Lambda)/T^2\dto^p\\
& &\Lambda &\Lambda}
\end{equation}
where $\Lambda =\{(x, y, z)|z=xy, \ \ \forall (x, y)\in X_2\}
\subset X_1\times X_1\times X_1$ is the graph of the groupoid
multiplication of $X_1\toto X_0$.
Let  $\Tilde{\Theta}$ be the one-form on $\pi^{-1}(\Lambda )
\subset R_1\times R_1\times R_1$ defined by
\begin{equation}
\Tilde{\Theta}= \Theta -\pi^*\pr_{12}^*\eta,
\end{equation}
where $\Theta =(\theta , \theta , -\theta )$ and $\pr_{12}:
 \Lambda \to  X_2$  is  the projection to the first two components.  Then 
\begin{equation}
\label{eq:dTheta}
d\Tilde{\Theta}=0.
\end{equation} 

By $\xi$ we denote the Euler vector field on $R_1$ 
generating the $S^1$-action.

\begin{lem}
\begin{enumerate}
\item $(\xi, \xi , \xi)\per \Tilde{\Theta}=1$;
\item $  \Tilde{\Theta}\in \Omega^1 (\pi^{-1}(\Lambda ))$
is basic with respect to the $T^2$-action as in Eq. \eqref{eq:T2},
so it descends to  a one-form $\Hat{\Theta}$ on $\pi^{-1}(\Lambda )/T^2$;
\item  $  \Hat{\Theta}$ defines  a flat connection 
on the $S^1$-principal bundle $\pi^{-1}(\Lambda )/T^2\stackrel{p}{\lon}\Lambda$.
\end{enumerate}    
\end{lem}
\begin{pf}
(i) is obvious.
For (ii)-(iii),
note that $\Tilde{\Theta}\in \Omega^1(\pi^{-1}(\Lambda ))$
is invariant under the natural  $T^3$-action induced
from the one on $R_1\times R_1 \times R_1$.
It is also quite clear that $\xi_1\per \Tilde{\Theta}
=\xi_2\per \Tilde{\Theta} =0$, where 
\begin{equation}
\label{eq:xi12}
\xi_1 =(\xi, 0, \xi ), \  \mbox{and } \xi_2 =(0, \xi , \xi )
\end{equation} 
are the generating vector fields of the $T^2$-action as
in Eq. \eqref{eq:T2}. 
Hence $  \Tilde{\Theta}$ is basic with respect to this action,
and descends to a  one-form $\Hat{\Theta}$ on  $\pi^{-1}(\Lambda )/T^2$,
which is easily seen to be  a flat connection for the
$S^1$-bundle $\pi^{-1}(\Lambda )/T^2\stackrel{p}{\lon}\Lambda$.
\end{pf}

Now assume that $R_1 \toto R_0$ is a Lie  groupoid  $S^1$-central extension
over $X_1\toto X_0$. Then $R_1\to X_1$ is a principal $S^1$-bundle.
Assume, moreover, that $\theta \in \Omega^1 (R_1)$  is a pseudo-connection
of the  extension whose corresponding pseudo-curvature equals 
$\eta +\omega\in Z^3_{DR}(X\lcom )$. That is,
\begin{equation}
\partial \theta =\pi^* \eta, \ \ \ d\theta =-\pi^* \omega .
\end{equation}

The  proposition below describes the
relation between $\theta $ and the groupoid 
structure on $R_1\toto R_0$.  First, let us  fix some
notations as follows.

\begin{eqnarray}
&&\tee: R_0\lon R_1, \ \ u\to \tu\\
&& \tee_2: R_0\lon R_2, \ \ u\to (\tu, \tu )\\
&&\epsilon_2: X_0\lon X_2,  \ \ u\to (u, u)\\
&&\kappa: X_1\lon X_2,  \ \ x\lon (x, x^{-1}) \label{eq:alp}
\end{eqnarray}

Let $\eta_0$ be the one-form on $X_0$ given  by
\begin{equation}
\eta_{0}=\epsilon_2^* \eta .
\end{equation}

\begin{prop}
 Let  $R_1 \toto R_0$ be  a Lie groupoid  $S^1$-central extension
over $X_1\toto X_0$. 
Let $\theta \in \Omega^1 (R_1)$  be  a pseudo-connection
 whose corresponding pseudo-curvature equals 
$\eta +\omega\in Z^3_{DR}(X\lcom )$. Then
\begin{enumerate}
\item $\tee^* \theta =\eta_{0}$
\item $\ti^* \theta +\theta =\ts^* \eta_{0}+\pi^*\kappa^*\eta. $
\end{enumerate}

In particular, if $\theta $ is a  connection,
then $$ \tee^* \theta =0, \ \ \ \ti^* \theta=-\theta .$$
\end{prop}
\begin{pf}  
(i) It is clear that $\tee_2^* \partial  \theta =\tee^* \theta $. On the
other hand, we have $\tee_2^* \pi^* \eta =(\pi \smalcirc \tee_2 )^*\eta
=\epsilon_2^* \eta =\eta_{0}$. 
Thus we have $\tee^* \theta =\eta_{0}$.

(ii) Given any $\x\in R_1$, $\forall \delta_{\x}\in T_{\x}R_1$,
 consider the tangent vector $(\delta_{\x}, \ti_{*}\delta_{\x})$
of $R_2$ at the point $(\x, \x^{-1})$. It is clear that
$\tm_{*}(\delta_{\x}, \ti_{*}\delta_{\x} )=\tee_* \ts_{*}\delta_{\x} $. So
$$(\partial  \theta) (\delta_{\x}, \ti_{*}\delta_{\x} )=
\delta_{\x}\per \theta +\ti_{*}\delta_{\x}\per\theta 
-\tee_* \ts_{*}\delta_{\x}\per \theta
=\delta_{\x}\per (\theta +\ti^* \theta -\ts^* \eta_{0}).$$
On the other hand, $(\pi^* \eta )(\delta_{\x}, \ti_{*}\delta_{\x})=
\delta_{\x} \per \pi^*\kappa^*\eta$.  (ii) thus follows.
\end{pf}

\begin{numrmk}
{\em In the case of  an $S^1$-gerbe over a  manifold,
 the conditions that $ \tee^* \theta =0, \  \ti^* \theta=-\theta $
were included  in the definition of a connection 
\cite{Brylinski:book, Hitchin, Murray}.
From the above lemma, we see that they are 
easy consequences of the condition $\partial \theta=0$.
}
\end{numrmk}

\begin{prop}
\label{thm:prescribed}
Let  $R_1 \toto R_0$ be  a Lie groupoid  $S^1$-central extension
over $X_1\toto X_0$.
Let $\theta \in \Omega^1 (R_1)$  be  a pseudo-connection
 whose  pseudo-curvature equals 
$\eta +\omega\in Z^3_{DR}(X\lcom )$. 
Then the flat $S^1$-bundle  $p:  \pi^{-1}(\Lambda )/T^2 \lon \Lambda$
as in diagram \eqref{eq:A}  is holonomy free. 
\end{prop}  
\begin{pf} 
By $\TLambda \subset R_1\times R_1\times R_1$,
 we denote the graph of the groupoid  multiplication of $R_1\toto R_0$. 
 It is clear that $p (\TLambda/T^2 )=\Lambda$.
Given any $(\x_1 ,\y_1, \z_1), \ (\x_2 ,\y_2, \z_2)\in
\TLambda$, if $p[(\x_1 ,\y_1, \z_1)]=p[(\x_2 ,\y_2, \z_2)]$,
then $\pi (\x_1 ,\y_1, \z_1)=\pi (\x_2 ,\y_2, \z_2)$.
This implies that  $\x_1 =s \cdot \x_2$ and $\y_1=t\cdot  \y_2$.
Hence $\z_1=\x_1 \y_1=(s \cdot \x_2)(t\cdot  \y_2)=
(st)\cdot (\x_2 \y_2)=(st)\cdot\z_2 $,
and therefore $[(\x_1 ,\y_1, \z_1)]=[(\x_2 ,\y_2, \z_2)]$.
Hence $\TLambda/T^2$ is indeed a  section of the
$S^1$-bundle $p: \pi^{-1}(\Lambda )/T^2 \lon \Lambda$.
From the equation  $\partial  \theta =\pi^* \eta$, it
follows that $\Tilde{\Theta}$ vanishes on $\TLambda$.
So $\TLambda/T^2$ is indeed a horizontal section.  \end{pf}

\subsection{Sufficient condition}

In this subsection, we investigate the inverse question
to Proposition  \ref{thm:prescribed}. Namely,
given a Lie groupoid $X_1\toto X_0$ and
a 3-cocycle $\eta +\omega \in Z^3_{DR}(X\lcom )$,
if  $\pi: R_1\to X_1$ is an $S^1$-bundle
and $\theta \in \Omega^1 (R_1 )$
is a connection 1-form  of the bundle so that $d\theta =-\pi^*\omega$
and the corresponding  $S^1$-flat bundle 
$p: \pi^{-1}(\Lambda )/T^2 \lon \Lambda$
is  holonomy free, does $R_1$ always admit a structure
of groupoid $S^1$-central extension over $X_1\toto X_0$
so that $\theta $ is a  pseudo-connection with 
$\eta +\omega$ being  its pseudo-curvature?
Throughout this subsection, we will keep this assumption
and  all the  notations. Our  method is a 
modification of the one used in \cite{WeinsteinXu},
where a special case  was investigated.

Let  $\Lambda_1$ be  a  horizontal
section of the flat bundle $p:  \pi^{-1}(\Lambda )/T^2\lon \Lambda$.
Set $\TLambda =\tau^{-1}(\Lambda_1 )\subset \pi^{-1}(\Lambda )$,
which is clearly a  $T^2$-invariant submanifold.
It is also clear that $\mbox{dim}{\TLambda} =\mbox{dim}\Lambda +2
=\mbox{dim}X_2 +2$, and $\Tilde{\Theta}$ vanishes when
being restricted to $\TLambda$.

\begin{lem}
\begin{enumerate}
\item $\pi (\TLambda )=\Lambda $; and
\item $\TLambda$ is a graph over $R_2$
\end{enumerate}
\end{lem}
\begin{pf} (i) is  obvious.
 
(ii) Let $\tpr_{12} : R_1  \times    R_1\times R_1\to R_1\times R_1$
 be  the projection to its first two components.
Clearly $\tpr_{12} (\TLambda )\subseteq R_2$.
Let $(\x, \y )\in R_2$ be any point, and write $(x, y)\defequal \pi (\x, \y )$.
Then $(x, y, xy)\in \Lambda$. Assume that
$  (\x_1 , \y_1  , \z_1  )\in \TLambda$ such that
$\pi (\x_1 , \y_1  , \z_1  ) =(x, y, xy)$.
Then $\x=s\cdot \x_1$ and $ \y =t\cdot \y_1 $ for some 
$s, t\in S^1$. Since $\TLambda$ is $T^2$-invariant,
it thus follows that $(\x, \y, st\cdot \z_1 ) =(s, t)\cdot (\x_1 , \y_1  , \z_1  )\in \TLambda $.
This shows that $\tpr_{12} (\TLambda )=R_2$.

To show that $\TLambda$ is indeed a graph over $R_2$,
assume that $  (\x , \y  , \z  ),  \ (\x , \y  , \z_1  )
$ are  two points in $\TLambda$. Then it is 
clear that $\pi (\x , \y  , \z  ) =\pi (\x , \y  , \z_1  )$,
i.e.  $p\smalcirc \tau (\x , \y  , \z  ) =p\smalcirc \tau  (\x , \y  , \z_1  )$.
Since $\tau (\x , \y  , \z  )$ and $\tau (\x , \y  , \z_1  )\in 
\Lambda_1$ and $\Lambda_1$ is a section for $p$, it
follows that $\tau (\x , \y  , \z  )=\tau  (\x , \y  , \z_1  )$.
Hence $(\x , \y  , \z  )=(s, t)\cdot (\x , \y  , \z_1  )$
for some $(s, t)\in T^2$, which implies that
$s=t=1$ and $ \z = \z_1 $. \end{pf}

Now $\TLambda$ defines a smooth map $\tilde{m}': R_2 \lon R_1$,
$(\x , \y )\to \x *\y$. By construction, the operation  $*$ satisfies the
condition
\begin{equation}
\label{eq:st}
(s\cdot \x)*(t\cdot \y)=st \cdot (\x*\y )
\end{equation}
for all  $ s,t \in S^{1}$ and $(\x,\y) \in R_2$.

Obviously, $\tilde{m}'$ commutes with the projection
$\pi$.  Therefore for any  triple $(\x , \y , \z)\in R_3$, both elements
$(\x *\y)*\z$ and $\x *(\y*\z) \in R_1$ have the same image 
under the projection $\pi$, so they must  differ by a unique
 element in $S^1$. Hence,  we obtain a function $g: R_3\to S^1$.
  Note that  Eq. \eqref{eq:st}
implies that  $g$ descends to a   function
 on $X_3$. Hence,  symbolically,  we may write
$$g(x , y , z)=\frac{(\x *\y)*\z}{\x *(\y*\z)}, \ \ 
\forall  (x , y , z)\in X_3, $$
where $(\x , \y , \z)\in R_3$ is any point such that
$\pi (\x , \y , \z) =(x , y , z)$.
We  call  $g(x , y , z)$ the {\em modular function} of $\theta$.

Note  that $g(x , y , z)$  is independent of the choice of
the horizontal section $\Lambda_1$ of the
flat bundle  $p:  \pi^{-1}(\Lambda )/T^2\lon \Lambda$,
and therefore   depends solely  on $\theta$.

\begin{prop}
If the modular function $g(x , y , z)$ is equal to 1,
$\TLambda$ defines a Lie  groupoid structure on $R_1$,
which is an $S^1$-central extension of $X_1\toto X_0$  with
$\theta$ being  a  pseudo-connection  
and $\eta +\omega$ the corresponding
pseudo-curvature. 
\end{prop}
\begin{pf}
  By assumption, we know that
 $\x *\y$ is indeed associative.

Now we  need to show the existence of units.
For this  purpose, we show that there exists a unique
section for the principal $S^1$-bundle $R_1\stackrel{\pi}{\to }X_1$
 over  the unit space $\epsilon (X_0 )$,
  namely $\epsilon': \epsilon (X_0 )\lon R_1$,
 $\epsilon (u)\stackrel{\epsilon'}{\to} \tu$ such that $(\tu, \tu, \tu )\in \TLambda$ for 
any $u\in X_0 $. Let
$(\tu_1, \tu_2, \tu_3)\in \TLambda$ be any  point such that
$\pi (\tu_1, \tu_2, \tu_3 )=(u, u, u)$.
Then $\tu_2 =s \cdot \tu_1$ and $\tu_3 =t\cdot \tu_1$ for
some $s, t\in S^1$. Let $\tu =(st^{-1})\cdot \tu_1 $.
Then
$$ (\tu, \tu , \tu )
=(st^{-1}, t^{-1})\cdot (\tu_1, \tu_2, \tu_3 )\in \TLambda .$$
Assume that $(\tv, \tv ,\tv) $ is another point 
in $\TLambda$ such that $\pi (\tv, \tv,  \tv)
=(u, u, u)$. From the equation
 $(p\smalcirc \tau) (\tu , \tu , \tu)
=(p\smalcirc \tau ) (\tv, \tv ,\tv)$, we deduce that
$\tau (\tu , \tu , \tu)=\tau (\tv, \tv,  \tv)$.
Therefore $(\tu , \tu , \tu)=(s, t)\cdot (\tv, \tv,  \tv)$ for
some $s, t\in S^1$,
which means that $\tu =s \cdot \tv, \tu=t \cdot \tv, \tu
=st\cdot\tv$. This implies that $s=t=1$ and hence $\tu=\tv$.

Next we prove that $\tu  *\x =\x$ and $\x  *\tv =\x$ if
$\ts (\x )=u $ and $\tilde{t} (\x ) =v$.
 By construction, we have $\tu *\tu =\tu$.
From the associativity assumption, we have 
$$\tu * (\tu  *\y )=(\tu  *\tu ) *\y =\tu * \y, \mbox{ if }
\ts (\y)=u .$$
We must prove that $\x$ is of the form $\tu * \y $.
Let $x=\pi ( \x )$.
Since $(u, x, x)\in \Lambda$, there exists $(a, b, c)\in
\TLambda$ such that $\pi (a, b, c)=(u, x, x)$.
Thus $\tu =s\cdot a $ and $\x =t \cdot c$ for some
$s, t\in S^1$.
So $$(\tu , ts^{-1} \cdot b , \x)=(s\cdot a,  ts^{-1} \cdot b ,
 t\cdot c) =(s, ts^{-1})\cdot  (a, b, c)\in \TLambda .$$
Thus $\x =\tu  *( ts^{-1} \cdot b )$. In conclusion, we
have $\tu  *\x =\x$. Similarly, one proves that
$\x  *\tv =\x$.

Finally, we need to show the existence of inverse.
For any $\x \in R_1$, let $x=\pi (\x )$ and $\ts (\x )=v$.
Since $(x, x^{-1}, v)\in \Lambda$,
there exists $(\x_1 , \y_1, \z_1 )\in \TLambda$ such that
$\pi (\x_1 , \y_1, \z_1 )=(x, x^{-1}, v)$.
One may assume that $\x_1=\x$ by using the $T^2$-action.
Since $\pi \z_1 = \pi \tv$, we have $\tv=t\cdot \z_1$.
Thus $(\x, t\cdot \y_1 , \tv )=(1, t)\cdot (\x, \y_1, \z_1 )\in
\TLambda$. This shows that the right inverse of $\x$ exists.
Similarly, one shows that the left inverse exists as well.
It is then standard that the left and right inverses 
must coincide. 
This concludes the proof. \end{pf}

In general, the modular function is
not necessarily equal to $1$. Nevertheless,
we have the following characterization.

\begin{prop}
\label{pro:dg}
The modular function
$g: X_3\to S^1$ defines  an $\rr/\zz$-valued   groupoid  3-cocycle. I.e., 
$dg =0, \ \ \partial  g=1$.
\end{prop}
\begin{pf} 
Let $(x(t), y(t), z(t))$ be  any smooth curve
in $X_3$, and $(\xt , \yt , \zt )$
a smooth curve in $R_3$ such that
$\pi (\xt , \yt , \zt ) =(x(t), y(t), z(t))$.
Write  $\ut =\xt * \yt$, $\vt =\yt *\zt$, and
$\hth= (\xt * \yt )*\zt$, $\htt=\xt *(\yt *\zt )$.

Since $(\xt , \yt, \ut ), \ ( \ut, \zt , \hth)\in \TLambda$, 
we have
\be
&&\dxt \per \theta +\dyt  \per \theta -\dut \per \theta =\eta (\dxt
 * \dyt ), \ \mbox{and}\\
&&\dut \per \theta +\dzt  \per \theta -\dht  \per \theta =\eta (\dut 
* \dzt ),
\ee 
where, by abuse of notation, we use the same symbol $*$
to denote the  induced tangent map $TR_2\to TR_1$.
It follows that
\begin{equation}
\label{eq:h1}
\dht  \per \theta =\dxt \per \theta +\dyt  \per \theta+
\dzt \per \theta-\eta (\dxt * \dyt )-\eta (\dut * \dzt).
\end{equation}     

Similarly, one proves that
\begin{equation}
\label{eq:h2}
\dhtt\per \theta =\dxt \per \theta +\dyt  \per \theta+
\dzt \per \theta-\eta (\dxt * \dvt)-
\eta (\dyt *  \dzt).
\end{equation} 

Since $\partial  \eta =0$,  Eqs \eqref{eq:h1} and \eqref{eq:h2}
imply that $\dht  \per \theta=\dhtt\per \theta$.

Let $f(t)=\frac{\hth}{\htt}$. Therefore ${\hth}=f(t)\cdot {\htt}$,
which implies that
$$\dht=f(t)\cdot \dhtt+\dft \xi,$$
where $\xi$ is the Euler vector field on $R_1$.
Pairing with $\theta$ on both sides, one obtains that $\dft =0$.

Finally, the identity $\partial g=1$ can
be verified directly.
\end{pf}

\begin{cor}
\begin{enumerate}
\item For any $y\in X_1$, we have $g(s(y), y, t(y))=1$. In particular,
 $\forall u\in X_0$, we have $g(u, u, u)=1$;
\item If $X_1$ is $s$-connected, then $g(x, y, z)=1$.
\end{enumerate}
\end{cor}
\begin{pf}  Since $\partial g=1$,  we have
$$g(y, z , w)g(xy, z, w)^{-1}g(x, yz, w)g(x, y, zw)^{-1}g(x, y, z)=1. $$
By letting $x=s(y)$ and $z=t(y)$, we obtain that
$g(s(y), y, t(y))=1$.

For any $(x, y, z)\in X_3$, if $X_1$ is  $s$-connected,
then $x$ can be connected to $s(y)$ by a smooth path in the $t$-fiber
$t^{-1}(s(y))$ while $z$ can be connected to $t(y)$
by a smooth path in the $s$-fiber
$s^{-1}(t(y))$. In other words, $(x, y, z)$  and
 $(s(y), y, t(y))$ belong to the same connected
component of $X_3$. Thus $g(x, y, z)=1$ according to Proposition
\ref{pro:dg}.\end{pf}

An immediate consequence is the following

\begin{prop}
\label{pro:4.8}
Let  $\poidd{X_1}{X_0}$ be an $s$-connected Lie groupoid,
 and $\eta +\omega \in Z^3_{DR} (X\lcom )$ a 3-cocycle,
 where $\eta \in \Omega^1(X_2 )$ and $\omega \in \Omega^2 (X_1 )$.
  Assume that $\omega $ represents an integer cohomology class 
in $H^2_{DR}(X_1 )$, so that there exists an $S^1$-bundle 
$\pi:R_1\to X_1$ with a connection 1-form $\theta \in \Omega^1 (R_1)$,
 whose curvature is $-\omega$. 
If the associated $S^1$-bundle  $p:  \pi^{-1}(\Lambda )/T^2 \lon \Lambda$
as in diagram \eqref{eq:A}  is holonomy free,  then
$R_1\to X_1$ is a Lie groupoid  $S^1$-central extension with
$\theta$ being a pseudo-connection and $\eta +\omega$ being the
pseudo-curvature. In particular,
$\eta +\omega$ is of integer class  in $H^3 (X\lcom , \integers )$.
\end{prop}

%
%

\subsection{Properties of 3-cocycles}

In this subsection we study some  geometric properties of
3-cocycles of the De-Rham double complex of a 
Lie groupoid, which are important for our constructions
in the next section.  

Let  $\eta +\omega\in Z^3_{DR} (X\lcom )$  be a de-Rham
 three-cocycle, where $\eta \in \Omega^1 (X_2 )$ and $\omega \in \Omega^2 (X_1 )$.
Then 
\begin{equation}
\partial  \eta =0, \ \ \ \partial  \omega +d\eta   =0 , \ \ \ d\omega =0.
\end{equation}

By $X_1^s$ and  $X_1^t$ we denote the $s$- and $t$-fibrations of
$X_1\toto X_0$, respectively.
Define a  leafwise  one-form   $\lambda^{r} $ on $ X_1^t$  by
$$\lambda^r  (\delta_{x})=\eta (r_{x^{-1}*}\delta_{x} , \  0_{x} ), \ \ \ 
\forall \delta_x \in T_x X_1^t .$$
Similarly, let $\lambda^l$  be the leafwise  one-form on
$ X_1^s$ given by
$$\lambda^l  (\delta_{x})=\eta (0_{x}, \  l_{x^{-1} *}\delta_{x} ),
\ \ \ \forall \delta_x \in T_x X_1^s .$$
Here $r_{x^{-1}}$ and $l_{x^{-1}}$ denote the right and the left  translations,
respectively.

Note that  $\lambda^r $ (or $\lambda^l $)
is in general not right (left)-invariant.

By $A\to X_0$  we denote the Lie algebroid of $X_1\toto X_0$.
For any section $\vfx\in \gm (A)$, we denote, respectively, by
$\Vec{\vfx}$ and $\ceV{\vfx}$ the  right invariant and  the left invariant
vector fields   on $X_1$ corresponding to $\vfx$.

\begin{lem} 
\label{lem:3.1}
For any $\vfx\in \gm (A)$,
\begin{enumerate}
\item  $\eta (\Vec{\vfx}(x) , 0_y )=
\lambda^r (\Vec{\vfx}(xy))-\lambda^r (\Vec{\vfx}(x))$, $\forall (x, y)\in X_2$;
\item $\eta (0_x ,  \ceV{\vfx}(y))=
 \lambda^l (\ceV{\vfx}(xy) )-\lambda^l (\ceV{\vfx}(y))$, $\forall (x, y)\in X_2$;
\item $\lambda^r (\Vec{\vfx})(u)=\lambda^l (\ceV{\vfx})(u)=0$, $\forall u\in X_0$; and
\item $\eta (\Vec{\vfx}(x),  -\ceV{\vfx}(x^{-1}))= \lambda^l (\ceV{\vfx} )(x^{-1})-
\lambda^r (\Vec{\vfx})(x)$, $\forall x\in X_1 $.
\end{enumerate}
\end{lem}
\begin{pf} Consider the curve $ r(t)=(\exp{t\Vec{\vfx}}, x, y)$ in $X_3$
 through
the point $(s(x), x , y)$. By definition, we have
$$\dot{r}(0)\per \partial  \eta=
\eta (0_x , 0_y )-\eta (\Vec{\vfx}(x),  0_y )+\eta (\Vec{\vfx}(s(x)),  0_{xy})
-\eta (\Vec{\vfx}(s(x)), 0_{x}).$$
Thus (i) follows immediately since $\partial  \eta=0$.
Similarly (ii) can be proved by considering the curve $(x, y, \exp{t\ceV{\vfx}})$
through the point $(x, y , t(y))$.
(iii) follows from (i) and (ii) by taking $x=y=u\in X_0$.  Finally, 
using (i)-(iii), we have

\be
\eta (\Vec{\vfx}(x) , -\ceV{\vfx}(x^{-1}))&=&\eta (\Vec{\vfx}(x),  0_{x^{-1}})
-\eta (0_{x},  \ceV{\vfx}(x^{-1}) )\\
&=&[\lambda^r (\Vec{\vfx})(xx^{-1})-\lambda^r (\Vec{\vfx})(x)]
-[\lambda^l  (\ceV{\vfx})(xx^{-1}) -\lambda^l (\ceV{\vfx}) (x^{-1})]\\
&=& \lambda^l (\ceV{\vfx} )(x^{-1})-\lambda^r (\Vec{\vfx})(x).
\ee
Thus (iv) follows.
\end{pf}

For any  $\vfx\in \gm (A)$, by  $\vfx^{r}$ and $\vfx^l$ we denote  the vector fields
  on $X_2$ given by
$\vfx^{r}(x, y)= (\Vec{\vfx}(x),  0_y)$ and $ \vfx^l (x, y)=(0_{x} , \ceV{\vfx}(y))$,
$\forall (x, y)\in X_2$.
It is clear that the flows of  $\vfx^{r}$ and $ \vfx^l $ are, respectively, given
by
\begin{equation}
\label{eq:XY}
\phi_{t}(x, y)=(\exp{(t\Vec{\vfx})}x, y), \  \ \
\psi_{t}(x, y)=(x , y \exp{(t\ceV{\vfx})}).
\end{equation}

\begin{lem} 
\label{lem:eta}
For any $\vfx, \vfy\in \gm (A)$,
\begin{enumerate}
\item $(d\eta ) ( \vfx^{r}, \vfy^l )(x, y)=\Vec{\vfx} (\lambda^l (\ceV{\vfy} ) )(xy)
-\ceV{\vfy} (\lambda^r (\Vec{\vfx } ))(xy)$. 
\item $(d\eta ) ( \vfx^{r}, \vfy^r )(x, y)=(d\lambda^r )(\Vec{\vfx}, \Vec{\vfy})(xy)
-(d\lambda^r  )(\Vec{\vfx}, \Vec{\vfy})(x)$.
\item  $(d\eta ) ( \vfx^{l}, \vfy^l )(x, y)=  (d\lambda^l )(\ceV{\vfx}, \ceV{\vfy} )(xy)
-(d\lambda^l )(\ceV{\vfx}, \ceV{\vfy})(y)$.
\end{enumerate}
\end{lem}
\begin{pf} (i) From Eq. \eqref{eq:XY},
 one easily sees that the vector fields $\vfx^{r}$ and $ \vfy^l$ 
commute with each other: $[\vfx^{r}, \vfy^l ]=0$.

According to Lemma \ref{lem:3.1} (ii), $\eta ( \vfy^l )(x, y)=
\lambda^l (\ceV{\vfy}(xy)) -\lambda^l (\ceV{\vfy}(y))$. It thus
follows that 
$$\vfx^{r}  (\eta ( \vfy^l ) ) (x, y)
=\difft|_{t=0}[ \lambda^l (\ceV{\vfy}(\exp{(t\Vec{\vfx})}xy)-
\lambda^l (\ceV{\vfy}(y)) ]
=\Vec{\vfx} (\lambda^l (\ceV{\vfy} ) )(xy).$$
Similarly, one shows that $ \vfy^l (\eta (\vfx^{r} ))(x, y)=
\ceV{\vfy} (\lambda^r (\Vec{\vfx })(xy)$. (i) thus
follows.

(ii) We have
\be
\vfx^{r} (  \eta ( \vfy^r ) ) (x, y)&=&\vfx^{r}  [ \lambda^r (\Vec{\vfy}(xy))-
\lambda^r (\Vec{\vfy}(x))]\\
&=&\difft|_{t=0} [ \lambda^r (\Vec{\vfy}(\exp{(t\Vec{\vfx})}xy))-
\lambda^r (\Vec{\vfy}( \exp{(t\Vec{\vfx})}x))]\\
&=& \Vec{\vfx}(\lambda^r (\Vec{\vfy} ))(xy)- \Vec{\vfx}(\lambda^r (\Vec{\vfy} ))(x).
\ee
Hence
\be
(d\eta ) ( \vfx^{r}, \vfy^r )(x, y)&=& \vfx^{r}(\eta ( \vfy^r ))(x, y)-
\vfy^{r}(\eta ( \vfx^r )) (x, y) -\eta ([\vfx^r , \vfy^r ])(x, y)\\
&=& (d\lambda^r )(\Vec{\vfx}, \Vec{\vfy})(xy)- (d\lambda^r )(\Vec{\vfx}, \Vec{\vfy})(x).
\ee

(iii) can be proved similarly.\end{pf}

\begin{prop}
Assume that $\eta +\omega \in Z^3_{DR}(X\lcom )$ is a 3-cocycle.
\label{cor:3.3}
\begin{enumerate}
\item $\epsilon^* \omega =-d\eta_0$;
\item For any $\vfx, \vfy\in \gm (A)$,
$\omega (\Vec{\vfx}, \ceV{\vfy})=\Vec{\vfx}\lambda^l (\ceV{\vfy})-\ceV{\vfy}\lambda^r (\Vec{\vfx}).$
\item $\omega -d\lambda^r\in \Omega^2 (X^t)$ is
a right invariant  (leafwise) closed $2$-form, and
therefore induces a  Lie algebroid $2$-cocycle $\omega^r
 \in \gm (\wedge^2 A^* )$.
\item $\omega -d\lambda^l\in \Omega^2 (X^s)$ 
is a left invariant  (leafwise) closed $2$-form,  and therefore 
induces  a Lie algebroid $2$-cocycle 
$\omega^l \in \gm (\wedge^2 A^* )$.
\item $\omega^r$ and $\omega^l$ are related by 
$$\omega^r +\omega^l +\rho^*d \eta_{0}=0, $$
i.e. $\omega^r $ and $-\omega^l $ are cohomologous 
Lie algebroid $2$-cocycles. Here $\rho: A\to TX_0$ is the anchor of
the Lie algebroid $A$.
\end{enumerate}
\end{prop}
\begin{pf} It is not difficult to see that $m_{*}\vfx^{r}(x, y)=\Vec{\vfx}(xy)$ 
and $m_{*}\vfy^l (x, y)=\ceV{\vfy}(xy)$.
 Thus
$$(\partial   \omega )( \vfx^{r}, \vfy^l )(x, y)=
-\omega (\Vec{\vfx}, \ceV{\vfy} )(xy).$$
On the other hand, according to Lemma \ref{lem:eta}(i),
we have
$$(d\eta ) ( \vfx^{r}, \vfy^l )(x, y)=\Vec{\vfx} (\lambda^l (\ceV{\vfy} ) )(xy)-
\ceV{\vfy} (\lambda^r (\Vec{\vfx }))(xy).$$
Since $\partial  \omega +d\eta =0$, it thus follows that
$$\omega (\Vec{\vfx}, \ceV{\vfy})(xy)=
\Vec{\vfx}\lambda^l (\ceV{\vfy} )(xy)-\ceV{\vfy}\lambda^r (\Vec{\vfx })(xy).$$
(ii) thus follows by letting $y=t(x)$.

For (iii), we note that
$$(\partial  \omega )( \vfx^{r}, \vfy^r )(x, y)=
\omega (\Vec{\vfx}, \Vec{\vfy})(x)-\omega (\Vec{\vfx}, \Vec{\vfy})(xy).$$
The conclusion thus follows from Lemma \ref{lem:eta}.
(iv) can be proved similarly.

Finally,  consider the map $ \kappa : X_1    \to  X_2$ 
as in Eq. \eqref{eq:alp}.
 It is clear that 
\begin{equation}
\label{eq:ax}
 \kappa_{*}\Vec{\vfx}(x)=(\Vec{\vfx}(x),  -\ceV{\vfx}(x^{-1}) ).
\end{equation}
Thus 
\be
(\kappa^* \eta ) (\Vec{\vfx})(x)&=&\eta  (\Vec{\vfx}(x),  -\ceV{\vfx}(x^{-1}))
\ \ \ \mbox{(by Lemma \ref{lem:3.1} (iv))}\\
&=&- \lambda^r (\Vec{\vfx})(x)+\lambda^l (\ceV{\vfx} )(x^{-1})\\
&=&  (-\lambda^r -\iota^* \lambda^l ) (\Vec{\vfx})(x).
\ee
It follows that $\kappa^* \eta =-\lambda^r -\iota^* \lambda^l $.
Here $\kappa^* \eta$ is considered as  a fiberwise one-form on
$X^t_1$ by restriction.

For any $\vfx\in \gm (A)$, write 
$\overline{\vfx}(x, x^{-1})=(\Vec{\vfx}(x),  -\ceV{\vfx}(x^{-1}) )
\in T_{(x, x^{-1})} X_2$. Eq. \eqref{eq:ax} means that 
$\kappa_{*}\Vec{\vfx}(x)=\overline{\vfx}(x, x^{-1})$.
Hence for any $\vfx, \vfy\in \gm (A)$,
\be
(d\eta ) (\overline{\vfx}(x, x^{-1}), \overline{\vfy}(x, x^{-1}))
&=& d(\kappa^* \eta) (\Vec{\vfx}(x), \Vec{\vfy}(x))\\
&=&-(d\lambda^r )(\Vec{\vfx}(x), \Vec{\vfy}(x))-
(d \lambda^l ) (\iota_{*}\Vec{\vfx}(x), \iota_{*} \Vec{\vfy}(x))\\
&=&-(d\lambda^r )(\Vec{\vfx}(x), \Vec{\vfy}(x))-(d \lambda^l )
(\ceV{\vfx}(x^{-1}), \ceV{\vfy}(x^{-1})).
\ee
 
On the other hand, we have $m_{*}\overline{\vfx}(x, x^{-1})= \epsilon_*s_{*}
\Vec{\vfx}(x)$
and $m_{*}\overline{\vfy}(x, x^{-1})=\epsilon_* s_{*}\Vec{\vfy}(x)$.
To see this, note that $(\exp{t\Vec{\vfx}} \cdot x, (\exp{t\Vec{\vfx}} \cdot x)^{-1})$
is the flow generated by  $\overline{\vfx}(x, x^{-1})$  on $X_2$.
Thus we have 
\be
(\partial  \omega ) (\overline{\vfx}(x, x^{-1}), \overline{\vfy}(x, x^{-1}) )
&=& \omega (\Vec{\vfx}, \Vec{\vfy} )(x)+\omega (\ceV{\vfx}, \ceV{\vfy} )(x^{-1})
-\omega (s_{*}\Vec{\vfx}(x), s_{*}\Vec{\vfy}(x)).
\ee
(v) thus follows immediately. \end{pf}.

\subsection{Lie algebroid central extensions}

As in the last subsection, let 
$\eta +\omega \in Z^3_{DR}(X\lcom )$ be  a  de-Rham 3-cocycle
of  a Lie groupoid $X_1\toto X_0$ and $\omega $  represents
an integer  cohomology class in $H^2_{DR}(X_1 )$.
Let   $\pi: R_1\to X_1$ be  an $S^1$-bundle
and $\theta \in \Omega^1 (R_1 )$
 a connection 1-form  of the bundle so that $d\theta =-\pi^*\omega$.

Recall that, for a given Lie algebroid $A\to M$,
 any  Lie algebroid 2-cocycle $\gamma \in \gm  (\wedge^2 A^*)$
induces a Lie algebroid central extension $\tilde{A}=A\oplus (M\times \rr)$ 
as follows.  The anchor map  $\tilde{\rho} (\vfx+f)=\rho (\vfx)$, $\forall
\vfx\in \gm (A)$ and $f\in C^\infty (M)$, and the bracket is
$$[\vfx+f, \vfy+g]=[\vfx, \vfy]+(\rho(\vfx )(g)-\rho (\vfy )(f)+\gamma (\vfx, \vfy)), $$
$\forall \vfx, \vfy\in \gm (A)$ and $f, g\in  C^\infty (M)$.
Denote by $\tilde{A}^r$ and $\tilde{A}^l$ the Lie algebroid central
 extensions of
$A$ by the  2-cocycles $\omega^r$ and $-\omega^l$, respectively.
Then $\tilde{A}^r$ and $\tilde{A}^l$ are isomorphic, while the
isomorphism is given by
\begin{equation}
\tilde{A}^r \lon \tilde{A}^l ,  \ \ \ \vfx+ f\to \vfx+ ( f+\rho(\vfx)\per \eta_{0}), \ \ 
\forall \vfx\in \gm (A).
\end{equation}

Let $\theta^r =\theta + \pi^* \lambda^r $ and $\theta^l
=\theta + \pi^* \lambda^l$ be the fiberwise one-forms on $R^\tlt_1$
and $R_1^\ts$, respectively. Then $\theta^r$ is
a fiberwise connection one-form on  the fiberwise $S^1$-principal
bundle $ R^\tlt_1 \to  X_1^t$ with curvature being $-\omega+d\lambda^r$,
while $\theta^l $ is a fiberwise connection
 one-form on the fiberwise $S^1$-principal
bundle $R^\ts_1\to  X_1^s$ with curvature  being 
$-\omega+d\lambda^l$.
  For any $\vfx\in \gm (A)$, denote by $\Hat{\Vec{\vfx}}\in \calx (R_1^\tlt ) $ 
  the horizontal lift of 
$\Vec{\vfx} $ with respect to $\theta^r$, and $\Hat{\ceV{\vfx}}\in \calx (R^\ts_1 ) $ 
 the horizontal lift of $\ceV{\vfx}$ with respect to 
$\theta^l$. I.e., 
\begin{equation} \label{eq:Xhatr}
\left\{
\begin{array}{ll}
& \Hat{\Vec{\vfx}}\per (\theta +\pi^* \lambda^r )=0,\\
&\pi_{*} \Hat{\Vec{\vfx}}=\Vec{\vfx}
\end{array}
\right.
\end{equation}

and
\begin{equation} \label{eq:Xhatl}
\left\{
\begin{array}{ll}
& \Hat{\ceV{\vfx}}\per (\theta + \pi^* \lambda^l )=0, \\
&\pi_{*} \Hat{\ceV{\vfx}}=\ceV{\vfx}.
\end{array}
\right.
\end{equation}

Introduce  linear  maps $\phi: \gm (\tilde{A}^r )\lon \calx (R_1)$ 
and $\psi: \gm (\tilde{A}^l )\lon \calx (R_1)$, respectively,
 by 
\begin{equation}
\phi: \ \vfx+f \lon \Hat{\Vec{\vfx}}+(\pi^* s^* f)\xi
\end{equation}
and
\begin{equation}
\psi: \ \vfx+f \lon \Hat{\ceV{\vfx}}+(\pi^* t^* f)\xi,
\end{equation} 
$\forall \vfx\in \gm (A)$ and $f\in C^{\infty}(X_0 )$. Set
\begin{equation}
{\mathcal D}_s= \phi  \gm (\tilde{A}^r ), \ \ \mbox{and } \  {\mathcal D}_t=\psi
  \gm (\tilde{A}^l )\subset \calx (R_1)
\end{equation}

\begin{prop}
\label{pro:3.10}
\begin{enumerate}
\item   Both   $\phi$ and $\psi$ are Lie algebra homomorphisms.
\item Vector fields in ${\mathcal D}_s  $ and ${\mathcal D}_t$
mutually  commute.
\end{enumerate}
\end{prop}
\begin{pf} 
$\forall \vfx, \vfy\in \gm (A) $ and $f, g\in C^{\infty}(X_0)$,
we have
\be
&&[ \Hat{\Vec{\vfx}}+(\pi^* s^* f)\xi, \  \Hat{\Vec{\vfy}}+(\pi^* s^* g)\xi ]\\
&=&[\Hat{\Vec{\vfx}}, \Hat{\Vec{\vfy}}]+[\Hat{\Vec{\vfx}}, (\pi^* s^* g)\xi ]
+[(\pi^* s^* f)\xi , \Hat{\Vec{\vfy}} ]
\ee
Since the vector field $\Hat{\Vec{\vfx}}$ is $S^1$-invariant, we have
$[\Hat{\Vec{\vfx}}, \xi ] =0$.  On the other hand,  since
$s_{*}\pi_* \Hat{\Vec{\vfx}} =s_{*} \Vec{\vfx}=\rho(\vfx)$, it follows
that $\Hat{\Vec{\vfx}} (\pi^* s^* g)=\rho(\vfx)g$. Therefore,
$$[\Hat{\Vec{\vfx}}, (\pi^* s^* g)\xi ] =(\pi^* s^* g)[\Hat{\Vec{\vfx}}, \xi ] +(\Hat{\Vec{\vfx}} (\pi^* s^* g))\xi
=(\rho(\vfx)g) \xi .$$
Similarly,  one proves that $[(\pi^* s^* f)\xi , \Hat{\Vec{\vfy}} ] =-(\rho(\vfy)f) \xi$.
Finally note that
$$[\Hat{\Vec{\vfx}}, \Hat{\Vec{\vfy}}]=\Hat{\Vec{[\vfx, \vfy]}}+
\pi^* s^* \omega^r (\vfx, \vfy)\xi .$$
Hence it follows that $\phi $ is indeed a Lie algebra homomorphism.
Similarly, one proves that $\psi$ is  also a Lie algebra  homomorphism.

For the second part, for any $\vfx, \vfy\in \gm (A)$ and $f, g\in
C^{\infty}(X_0 )$,
we have
\be
[\phi (\vfx+f), \psi (\vfy+g)]&=&
[\Hat{\Vec{\vfx}}+(\pi^* s^* f)\xi,\ \Hat{\ceV{\vfy}}+(\pi^* t^* g)\xi]\\
&=&[\Hat{\Vec{\vfx}},  \Hat{\ceV{\vfy}}]+[\Hat{\Vec{\vfx}}, (\pi^* t^* g)\xi]
+[(\pi^* s^* f)\xi , \Hat{\ceV{\vfy}}].
\ee

Now
$[\Hat{\Vec{\vfx}}, (\pi^* t^* g)\xi]=[(\pi^* s^* f)\xi , \Hat{\ceV{\vfy}}]=0$
since $t_{*}\pi_{*}\Hat{\Vec{\vfx}}=t_{*}\Vec{\vfx}=0$ and
$s_{*}\pi_{*} \Hat{\ceV{\vfy}}=s_{*}\ceV{\vfy}=0$.

It remains to show that $[\Hat{\Vec{\vfx}},  \Hat{\ceV{\vfy}}] =0$.
For this, first of all, note that $\pi_* [\Hat{\Vec{\vfx}},  \Hat{\ceV{\vfy}}]
=[\Vec{\vfx}, \ceV{\vfy}]=0$.
By Eqs. \eqref{eq:Xhatr}, \eqref{eq:Xhatl}, we have  
 $\Hat{\Vec{\vfx}}\per \theta =
-\pi^* (\lambda^r (\Vec{\vfx}))$ and $\Hat{\ceV{\vfy}}\per \theta 
=-\pi^*  (\lambda^l (\ceV{\vfy} ))$. It thus follows
that
\be
[\Hat{\Vec{\vfx}},  \Hat{\ceV{\vfy}}] \per \theta&=&
\Hat{\Vec{\vfx}}\theta (\Hat{\ceV{\vfy}})-\Hat{\ceV{\vfy}}\theta (\Hat{\Vec{\vfx}})
-(d\theta )(\Hat{\Vec{\vfx}},\Hat{\ceV{\vfy}})\\
&=&\pi^* [-\Vec{\vfx}\lambda^l (\ceV{\vfy} )+\ceV{\vfy} \lambda^r (\Vec{\vfx})
+ \omega ( \Vec{\vfx}, \ceV{\vfy}) ] \ \ \mbox{(by Proposition \ref{cor:3.3} (ii))}\\
&=&0.
\ee
This concludes the proof of the proposition.\end{pf}

Introduce  distributions $\D_{s}$ and
$\D_{t}$  on  $\pi^{-1}(\Lambda )$ as follows.
For any $(\x, \y ,\z )\in \pi^{-1}(\Lambda )$, 
\begin{eqnarray}
\D_{s}|_{(\x, \y ,\z )}&=&\{(\phi (\vfx+f)(\x), 0, \phi (\vfx+f))(\z))|\forall \vfx\in   \gm (A), \ f\in
C^{\infty}(X_0 )\}\\
\D_{t}|_{(\x, \y ,\z )}&=&
\{(0, \psi (\vfx+f) (\y ), \psi (\vfx+f)) (\z) ) |\forall \vfx\in   \gm (A), \ f\in
C^{\infty}(X_0 )\}.
\end{eqnarray}

It is clear, from Proposition \ref{pro:3.10}, that both $\D_{s}$ and
$\D_{t}$  are integrable distributions. By $\D$, we denote the
distribution on $\pi^{-1}(\Lambda )$ defined by the equation
$\Tilde{\Theta}=0$. According to Eq. \eqref{eq:dTheta}, $\D$ is an integrable
distribution. 

\begin{prop}
\label{pro:3.13}
We have $\D_{s}\subseteq \D$ and $\D_{t}\subseteq \D$.
\end{prop}
\begin{pf} Let $v=(\phi (\vfx+f)(\x), 0, \phi (\vfx+f))(\z))$. Then
\be
v\per (\theta ,\theta , -\theta )&=& (\Hat{\Vec{\vfx}}\per \theta)(\x )
+f(s\smalcirc \pi (\x ))
-(\Hat{\Vec{\vfx}}\per \theta)(\z )-f(s\smalcirc \pi (\z )) \\
&=&-\lambda^r (\Vec{\vfx} )(x)+\lambda^r (\Vec{\vfx} )(z) .
\ee

On the other hand, 
$$v\per \pi^* \eta =(\Vec{\vfx}(x), 0 )\per \eta=
\eta (\Vec{\vfx}(x), 0_{y}).$$
Thus $$v\per \Tilde{\Theta}=v\per \Theta-v\per \pi^* \eta =0, $$
according to Lemma \ref{lem:3.1} (i).
Hence we have proved that $\D_s \subseteq \D$.
Similarly, one shows that $\D_t\subseteq \D$.
\end{pf} 

We are now ready to prove the main theorem of this section.

\begin{them}
Let $X_1 \toto X_0$ be an $s$-connected Lie groupoid, and $\eta +\omega
\in Z^3_{DR} (X\lcom )$ a de-Rham  3-cocycle, where $\eta \in \Omega^1(X_2)$ 
and $\omega \in \Omega^2 (X_1 )$.  Assume that $\omega $ represents
an integer cohomology class in $H^2_{DR}(X_1)$, so that there exists
an $S^1$-bundle $\pi:R_1\to X_1$ with a connection $\theta \in
\Omega^1 (R_1 )$, whose curvature is $-\omega$.
Assume that $\epsilon \upst R_1$ endowed with the flat connection
$\epsilon\upst\theta-\pi\upst \epsilon_2\upst\eta$ is holonomy free. 
 (Here $\epsilon :X_0 \to X_1$ and $\epsilon_2:
 X_0 \to X_2$ are the respective identity morphisms.)
Then $R_1\toto R_0$, where $R_0 =X_0$,
 admits in a natural way the structure of a Lie groupoid,
such that it   becomes an $S^1$-central extension of 
$X_1\toto X_0 $
and $\eta +\omega$ the pseudo-curvature of $\theta$.
In particular, $\eta +\omega$
 is of integer class in $H^3 (X\lcom , \zz )$.
\end{them}
\begin{pf}  Take a horizontal section $\epsilon' $
of the bundle $R_1|_{X_0}\to  \epsilon (X_0 )$:
$\epsilon (u)\to \tu, \forall u\in X_0 $.
Consider the foliation in $R_1\times R_1\times R_1$
 defined by ${\mathcal F}_{s }+{\mathcal F}_{t}$.
Let $$I=\{(\tu,\tu,\tu)| \ \ \forall  u\in X_0\}. $$
Then $I$ is transversal to the foliation 
 ${\mathcal F}_{s }+{\mathcal F}_{t }$. By the 
 method of characteristics \cite{CDW},
there is   a   minimal ${\mathcal F}_{s }+{\mathcal F}_{t }$-invariant
 submanifold $\TLambda$ containing $I$
which is  immersed in $R_1\times R_1\times R_1$. 
Proposition \ref{pro:3.13} implies that $\tilde{\Theta}=0$ when
being restricted to $\TLambda$.

It is clear that $\TLambda$ is $T^2$-invariant
since the $T^2$-generating vector fields
 $\xi_1 \in {\mathcal F}_{\alpha } $
and $\xi_2 \in {\mathcal F}_{\beta }$.
Now we need to show that $\TLambda$ is a graph over $R_{2}$.  
 Let $\pr_{12} :\ \ R_1\times R_1\times R_1\to R_1\times R_1$
 be the natural projection onto the
first two coordinates:
$\pr_{12} (x,y,z)=(x,y)$.  First, we  show that
$\pr_{12}  (\TLambda)=R_{2}$. 

Note that  $(\x,\y,\z)\in \TLambda$ 
if and only if $\x=\phi^{\alpha }\tu, \; \y=\phi^{\beta }\tu$,
  and $z=\phi^{\alpha }\phi^{\beta }\tu$
 for some $u\in X_0 $, where $\phi^{\alpha }$ is a product of 
 flows in  ${\mathcal D}_{s }$ and
 $\phi^{\beta }$ is a product of  flows in   ${\mathcal D}_{t }$.

Since  $\tlt_* {\mathcal D}_{s }=0$ and
 $\ts_*  {\mathcal D}_{t }=0$,    the  flow of ${\mathcal D}_{s }$
 preserves $\tlt $-fibres and similarly   the flow of
 ${\mathcal D}_{t }$ preserves $\ts $-fibres;  thus
$$\tilde{t} (\x)=\tilde{t} (\phi^{\alpha }\tu)=u,$$ 
and
$$\ts  (\y)=\ts (\phi^{\beta }\tu)=u. $$ 
I.e.,  $\tilde{t} (\x)=\ts (\y)$, namely, $(\x,\y)\in R_{2}$. 
Therefore,
$$\pr_{12} (\TLambda)\subseteq R_{2}. $$

Conversely, for any $(\x,\y)\in R_{2}$, assume that
 $\tilde{t} (\x)=\ts (\y) =u\in X_0 $.
 Since  $X_1\toto X_0$ is $t$-connected and
$(t\smalcirc \pi ) (\x)=u$, there exists a product  $\phi^{\alpha_{0}}$
of flows  generated by   vector
fields of the form $\Vec{\vfx}$ for $\vfx\in \gm (A)$,  such
that $\phi^{\alpha_{0}}(u)=\pi (\x)$. For each $\vfx\in
\gm (A)$, we  denote the flow of the
vector field $\Hat{\Vec{\vfx}}\in {\calx }(R_1)$ by
 $\Phi_{t}^{\alpha }$.
Since $\pi_* \Hat{\Vec{\vfx}}=\Vec{\vfx}$, then
$$\pi \smalcirc \Phi_{t}^{\alpha }=\phi_{t}^{\alpha_{0}}\smalcirc \pi . $$
As each  fibre of the $S^1$-bundle $R_1\to X_1$
 is  compact,  $\Phi_{t}^{\alpha }$
is defined provided that $\phi_{t}^{\alpha_{0}}$ is defined. 
Let $\Phi^{\alpha }$ denote the product of flows  corresponding 
to $\phi^{\alpha_{0}}$.  Then we have 
$$\pi \circ \Phi^{\alpha }=\phi^{\alpha_{0}}\circ \pi . $$
Hence 
$$\pi (\Phi^{\alpha }(\tu ))=
\phi^{\alpha_{0}} (u)  =\pi (\x).$$
 Therefore
$\x=\lambda \cdot \Phi^{\alpha }(\tu)$ for some $\lambda\in S^1$.
Note that the flow $\psi_{t}(\x)=t\cdot \x$ on $R_1$
 is generated by the   standard  Euler vector field $\xi $, which
 is also in ${\mathcal D}_{s }$.
 Hence   we conclude that 
there exists a product of flows  $\hat {\Phi }^{\alpha }$ generated by  
the vector fields in  ${\mathcal D}_{s }$ such that
$\x= \hat{\Phi }^{\alpha }(\tu)$.
 Similarly, we can  find a product of flows  $\hat {\Phi }^{\beta }$ 
 generated by the 
  vector fields in ${\mathcal D}_{t }$ such that
$\y= \hat{\Phi }^{\beta }(\tu )$.
So $(\x, \y, \hat{\Phi }^{\alpha }\hat{\Phi }^{\beta }(\tu ))
\in \TLambda$, i.e.  
$(\x, \y)=\pr_{12} (\x, \y,\hat {\Phi }^{\alpha }\hat {\phi }^{\beta }
(\tu )) \in \pr_{12} (\TLambda )$. 
Thus  we have proved  that $\pr_{12} (\TLambda)=R_{2}$.

Finally,  note that if 
$\x=\phi^{\alpha }\tu =\phi^{\alpha }_1 \tu$ and $ \y=\phi^{\beta }\tu$,
then $\z=\phi^{\alpha }\phi^{\beta }\tu$ and 
$z_1 =\phi_1^{\alpha }\phi^{\beta }u$.  
Thus $\z=\phi^{\beta } \phi^{\alpha } \tu=
\phi^{\beta } (\x)
=\z_1$. Similarly, one shows that $\z$ is also independent
of the choice of the flows $\phi^{\beta }$.
This shows that 
$\TLambda$ is indeed a graph over $R_{2}$.
Now the conclusion follows from Proposition
\ref{pro:4.8}.
 \end{pf}    

\begin{rmk}
It would be interesting to investigate how the integrability condition of 
Crainic-Fernandes \cite{CF} is related to the theorem above.
\end{rmk}


\begin{thebibliography}{00}
\bibitem{sga4}
Artin, M.,  Grothendieck, A., and   Verdier, J.-L.,
Th\'eorie des topos et cohomologie \'etale des sch\'emas,
  S\'eminaire de g\'eom\'etrie alg\'ebrique du Bois-Marie (SGA4) 
 \emph{Lecture Notes in Mathematics} {\bf 269, 270, 305}.
 Springer-Verlag, Berlin-New York, 1972-1973.

\bibitem{art.maz}
	 Artin, M., and Mazur, B.,
 Etale homotopy, {\em Lecture Notes in Mathematics} {\bf 100}
 Springer-Verlag, Berlin-New York 1969.


\bibitem{BX} Behrend, K., and Xu, P.,
$S^1$-bundles and gerbes over differentiable stacks,
\emph{C. R. Acad. Sci. Paris, Serie I}  {\bf 336}
(2003), 163-168.



\bibitem{BXZ}
Behrend, K., Xu, P. and Zhang, B.,
Equivariant gerbes over compact simple Lie groups.
\emph{C. R. Math. Acad. Sci. Paris, Serie I} {\bf 336} (2003), 251--256.



  
\bibitem{Brylinski:book}
Brylinski, J.-L., Loop spaces, characteristic classes and 
geometric quantization, Progress in Mathematics, {\bf 107} Birkh\"{a}user,
   1993.

\bibitem{Brylinski:gerbe}
Brylinski, J.-L.,  Gerbes on complex reductive Lie groups,
math.DG/0002158.

\bibitem{BSS}
Bunke, U.,   Schick, T.,   Spitzweck, M.,
Sheaf theory for stacks in manifolds and twisted cohomology for $S^1$-gerbes,
{\em Algebraic \& Geometric Topology} {\bf 7} (2007), 1007-1062.


\bibitem{Chatterjee}
Chatterjee, D.,
On the construction of abelian gerbes, PhD thesis (Cambridge), 1998.

\bibitem{CDW}
Coste, A., Dazord, P., and Weinstein, A,
 Groupo{\"\i}des symplectiques,
{\em  Publications du D{\'e}partement}
 Nouvelle Serie. A {\bf 2}
1--62,  Univ. Claude-Bernard, Lyon, (1987) 1--62.
 	
 	
									 
\bibitem{Crainic}
Crainic, M.,
Differentiable and algebroid cohomology, van Est isomorphisms, and
    characteristic classes, 
{\em Comment. Math. Helv.} {\bf 78} (2003),  681--721

\bibitem{CF}
Crainic, M., and Fernandes, R.,
 Integrability of Lie brackets, 
{\em Ann. of  Math.}  {\bf 157} (2003), 575--620.


\bibitem{deligne}
Deligne, P.,
Th\'eorie de {H}odge. {III},
{\em Inst. Hautes \'Etudes Sci. Publ. Math.} {\bf 44} (1974), 5--77.



\bibitem{Dijkgraaf}  Dijkgraaf, R., The mathematics of fivebranes,
{\em Documenta Math. Extra Volume ICM}  (1998).

\bibitem{Dupont}
Dupont, J., Curvature and characteristic classes,
   {\em Lecture Notes in Mathematics}  {\bf 640}
Springer-Verlag, Berlin-New York, 1978.

\bibitem{Duskin}
Duskin, J.,
An outline of nonabelian cohomology in a topos. I.
 The theory of bouquets and gerbes,
{\em Cahiers Topologie G\'eom. Diff\'erentielle}{\bf 23} (1982),  
165--191.

\bibitem{Freed}  Freed, D.,
 Higher algebraic structures and
quantization, {\em Commun. Math. Phys.} {\bf 159} (1994), 343-398.



\bibitem{GN}
Gawedzki, K., and  Reis, N.,
WZW branes and gerbes,
{\em Rev. Math. Phys.} {\bf 14} (2002),  1281--1334.

\bibitem{GS}
Ginot, G., and Stienon, M.,
G-gerbes, principal 2-group bundles and characteristic classes,
arXiv:0801.1238



\bibitem{Giraud}
Giraud, J., 
Cohomologie non ab\'elienne,
Die Grundlehren der mathematischen Wissenschaften {\bf 179},
 Springer-Verlag 1971.

\bibitem{sga1}
Grothendieck, A., Revtements \'etales et groupe fondamental,
S\'eminaire de g\'eom\'etrie alg\'ebrique du Bois-Marie 1960--1961 (SGA 1),
 \emph{Lecture Notes in Mathematics} {\bf 224}
 Springer-Verlag, Berlin-New York, 1971.



\bibitem{Hitchin}
Hitchin, N., Lectures on special Lagrangian submanifolds, 
Winter School on Mirror Symmetry, Vector Bundles and Lagrangian Submanifolds
(Cambridge, MA, 1999), \emph{AMS/IP Stud. Adv. Math.}  {\bf 23},  151--182.                             


\bibitem{hilsum-skandalis87}
Hilsum, M.,  and  Skandalis, G.,  Morphismes $K$-orient\'{e}s
d'espaces de feuilles et fonctorialit\'{e} en th\'{e}orie de Kasparov
(d'apr\`{e}s une conjecture d'A. Connes).
\emph{Ann. Sci. \'{E}cole Norm. Sup.} {\bf   20} (1987),  325--390. 


\bibitem{Kalkkinen}  Kalkkinen, J.,
 Gerbes of massive type II
configurations, {\em J. High Energy Physics} {\bf 9907} (1999),
 2-22.

\bibitem{Kostant}
 Kostant, B.,
 Quantization and unitary representations. I. Prequantization.
 Lectures in modern analysis and applications, III, 
{\em Lecture Notes in Math.} {\bf 170}  87--208.

\bibitem{LL}
Laumon, G., and Moret-Bailly, L.,
 Champs alg\'ebrique,,
 Results in Mathematics and Related Areas. 3rd Series. A Series of Modern Surveys in Mathematics, {\bf 39} Springer-Verlag, (2000).
 
\bibitem{LTX} Laurent-Gengoux, C., Tu, J.-L., and Xu, P.,
Chern-Weil map for principal bundles over groupoids.
{\em Math Z.} 
 {\bf 255} (2007), 451-491.

	

\bibitem{LU}
Lupercio, E., and Uribe, B.,
     Gerbes over orbifolds and twisted $K$-theory,
{\em Comm. Math. Phys.} {\bf  245} (2004),  449--489.


\bibitem{Mackenzie}
  Mackenzie, K.,
 General theory of Lie groupoids and Lie algebroids,
{\em  London Mathematical Society Lecture Note Series} {\bf 213}
 Cambridge University Press, Cambridge, 2005.



\bibitem{Meinrenken}
Meinrenken, E.,
The basic gerbe over a compact simple Lie group,
 {\em Enseign. Math. (2)}{\bf 49} (2003), 307--333.

\bibitem{MW}
 Mikami, K., and Weinstein, A., Moments and reduction for symplectic
groupoid actions, {\em Publ. RIMS Kyoto Univ.} {\bf  24} (1988), 121-140.


\bibitem{muhly-renault-williams87}
Muhly, P.,  Renault, J.  Williams, D.,
Equivalence and isomorphism for groupoid $C\sp *$-algebras,
\emph{J. Operator Theory} {\bf  17} (1987),  3--22.


\bibitem{Murray}
Murray, M. K.,
Bundle gerbes,
{\em J. London Math. Soc.} {\bf   54} (1996), 403--416.

\bibitem{MurrayS}
Murray, M. K.,  Stevenson, D.,
 Bundle gerbes: stable isomorphism and local theory,
{\em  J. London Math. Soc.} {\bf  62}  (2000), 925--937. 


\bibitem{Metzler} 
 Metzler, D., Topological and smooth stacks,
math.DG/0306176.

\bibitem{Moer}
Moerdijk, I.,
Orbifolds as groupoids: an introduction,
Orbifolds in mathematics and physics (Madison, WI, 2001),
{\em Contemp. Math.} {\bf 310}, 205--222,  2002.

 

\bibitem{Moer:regular}
Moerdijk, I.
Lie groupoids, gerbes, and non-abelian cohomology,
\emph{$K$-Theory} {\bf 28} (2003), no. 3, 207--258.

\bibitem{Moerdijk}
Moerdijk, I.,  Introduction to the language of stacks and gerbes, 
math.AT/0212266.

\bibitem{MoerdijkM}
Moerdijk, I.,   Mrcun, Introduction to Foliations and Lie Groupoids,
 Cambridge Studies in Advanced Mathematics, {\bf 91}, Cambridge University Press, Cambridge, 2003.

\bibitem{Noohi1}
Noohi, B.,
Fundamental groups of algebraic stacks,
{\em  J. Inst. Math. Jussieu} {\bf  3}  (2004),  69--103.

\bibitem{Noohi2}
Noohi, B.,
Foundations of topological stacks I, 
math.AG/0503247.


\bibitem{Stienon}
Stienon, M.,
Equivariant Dixmier-Douady classes, arXiv:0709.2720.


\bibitem{tu:05} Tu, J.-L.,  Groupoid cohomology and extensions,
{\em Trans. AMS} 
{\bf  358}  (2006),   4721--4747.



\bibitem{TXL}  Tu, J.-L., Xu, P., and Laurent-Gengoux, C., Twisted
  $K$-theory of differentiable stacks.  \emph{Ann. Sci. \'Ecole
  Norm. Sup.} (4)  {\bf 37}  (2004),  no. 6, 841--910.

\bibitem{TX1}
Tu, J.-L., and Xu, P., Chern character for twisted K-theory of orbifolds,
 {\em Adv Math}  {\bf 207}  (2006),   455--483. 

\bibitem{TX2}
Tu, J.-L., and Xu, P., The ring structure for equivariant twisted K-theory,
{\em J. Reine Angew. Math.} (to appear);
math.KT/0604160.


\bibitem{Sharpe}
Sharpe, E.,  Discrete torsion, {\em Phys. Rev.} {\bf 68} (2003), 7- 

\bibitem{Souriau}
 Souriau, J.-M.,
 Structure of dynamical systems, A symplectic view of physics,
 Translated from the French by C. H. Cushman-de Vries,
{\em Progress in Mathematics} {\bf  149} Birkhuser Boston, Inc., Boston, MA,
 1997.

\bibitem{Weil}
 Weil, A.,
 Sur les th\'eor\`emes de de Rham, {\em Comment. Math. Helv.} {\bf  26}, (1952). 119--145.
 	


\bibitem{WeinsteinXu}
Weinstein, A. and Xu, P.,
Extensions of symplectic groupoids and quantization.
\emph{J. Reine Angew. Math.} {\bf 417} (1991), 159--189.


\bibitem{WS}
Wirth, J., and Stasheff, J.,
 Homotopy transition cocycles,
math.AT/0609220.





\end{thebibliography}
\end{document}